\documentclass[oneside, 12pt]{amsart}
\usepackage{amscd, amssymb, amsmath, mathrsfs}
\usepackage[english]{babel}
\usepackage{booktabs}
\usepackage{tikz-cd}
\usepackage{url}
\usepackage[pdftex, colorlinks=true,  citecolor=blue, linkcolor=blue, linktocpage=true]{hyperref}

\DeclareUrlCommand\arXiv{\urlstyle{same}}

\newcommand\mylabel[1]{\label{#1}\marginpar{\vspace{-1ex}\medskip\medskip\footnotesize \tt #1}}
\renewcommand\mylabel[1]{\label{#1}}

\newcommand{\mydate}{
\number\day\space
\ifcase\month \or January\or February\or March\or April\or May\or June\or July\or August\or September\or October\or November\or December\fi 
\space\number\year}

\DeclareUrlCommand\arXiv{\urlstyle{same}}

\newtheorem{theorem}{Theorem}[section]
\newtheorem*{maintheorem}{Theorem}
\newtheorem{lemma}[theorem]{Lemma}
\newtheorem{proposition}[theorem]{Proposition}
\newtheorem{corollary}[theorem]{Corollary}

\theoremstyle{definition}

\newtheorem*{acknowledgement}{Acknowledgement}

\theoremstyle{remark}

\newcommand{\ZZ}{\mathbb{Z}}
\newcommand{\QQ}{\mathbb{Q}}
\newcommand{\RR}{\mathbb{R}}
\newcommand{\CC}{\mathbb{C}}
\newcommand{\FF}{\mathbb{F}}

\newcommand{\PP}{\mathbb{P}}

\newcommand{\GG}{\mathbb{G}}

\newcommand{\shA}{\mathscr{A}}

\newcommand{\shE}{\mathscr{E}}
\newcommand{\shF}{\mathscr{F}}

\newcommand{\shM}{\mathscr{M}}

\newcommand{\shN}{\mathscr{N}}
\newcommand{\shL}{\mathscr{L}}

\newcommand{\foX}{\mathfrak{X}}
\newcommand{\foY}{\mathfrak{Y}}

\newcommand{\aff}{\text{\rm aff}}
\newcommand{\Aff}{\text{\rm Aff}}
\newcommand{\alg}{\text{\rm alg}}

\newcommand{\Ass}{\operatorname{Ass}}

\newcommand{\Br}{\operatorname{Br}}

\newcommand{\Card}{\operatorname{Card}}

\newcommand{\CH}{\operatorname{CH}}

\newcommand{\contr}{\operatorname{contr}}

\newcommand{\Frac}{\operatorname{Frac}}

\newcommand{\Gal}{\operatorname{Gal}}
\newcommand{\GL}{\operatorname{GL}}

\newcommand{\gr}{\operatorname{gr}}

\newcommand{\Hom}{\operatorname{Hom}}
\newcommand{\Hilb}{\operatorname{Hilb}}

\newcommand{\id}{{\operatorname{id}}}

\newcommand{\I}{\text{\rm I}}
\newcommand{\II}{\text{\rm II}}
\newcommand{\III}{\text{\rm III}}
\newcommand{\IV}{\text{\rm IV}}
\newcommand{\tI}{\text{\rm \~{I}}}

\newcommand{\Kernel}{\operatorname{Ker}}

\newcommand{\invlim}{\varprojlim}

\newcommand{\lra}{\longrightarrow}

\newcommand{\Mat}{\operatorname{Mat}}

\newcommand{\maxid}{\mathfrak{m}}

\newcommand{\MW}{\operatorname{MW}}

\newcommand{\Num}{\operatorname{Num}}

\newcommand{\primid}{\mathfrak{p}}
\renewcommand{\O}{\mathscr{O}}

\newcommand{\Pic}{\operatorname{Pic}}

\newcommand{\pr}{\operatorname{pr}}
\newcommand{\Proj}{\operatorname{Proj}}

\newcommand{\quadand}{\quad\text{and}\quad}

\newcommand{\ra}{\rightarrow}

\newcommand{\rank}{\operatorname{rank}}
\newcommand{\red}{{\operatorname{red}}}

\newcommand{\sep}{{\operatorname{sep}}}

\newcommand{\Sing}{\operatorname{Sing}}

\newcommand{\Spec}{\operatorname{Spec}}

\newcommand{\Sym}{\operatorname{Sym}}

\newcommand{\uHom}{\underline{\operatorname{Hom}}}

\newcommand{\ind}{\text{\rm ind}}
\newcommand{\bk}{{\bar{k}}}
\newcommand{\bY}{{\bar{Y}}}
\newcommand{\bX}{{\bar{X}}}
\newcommand{\bE}{{\bar{E}}}

\newcommand{\bC}{{\bar{C}}}

\newcommand{\bJ}{{\bar{J}}}
\newcommand{\ba}{{\bar{a}}}

\newcommand{\Enr}{{\text{\rm Enr}}}
\newcommand{\inv}{{-1}}

\begin{document}

\title[Enriques surfaces over the integers]
      {There is no  Enriques surface over the integers}

\author[Stefan Schr\"oer]{Stefan Schr\"oer}
\address{Mathematisches Institut, Heinrich-Heine-Universit\"at,
40204 D\"usseldorf, Germany}
\curraddr{}
\email{schroeer@math.uni-duesseldorf.de}


\dedicatory{Third revised version, 19 July 2022}

\begin{abstract}
We show that there is no    family of Enriques surfaces over the ring of integers.
This extends   non-existence results of Minkowski for   families of finite \'etale schemes,
of Tate and Ogg for families of elliptic curves, and of Fontaine and Abrashkin for families of abelian varieties
and more general smooth proper schemes with certain  restrictions on   Hodge numbers.
Our main idea is to study the local system of numerical classes of invertible sheaves.
Among other things, our result also hinges on counting rational points, Lang's classification of   
rational elliptic surfaces in characteristic two,  the theory of exceptional Enriques surfaces due to 
Ekedahl and Shepherd-Barron, some recent results on the base of their versal deformation, Shioda's theory of Mordell--Weil lattices,
and an extensive  combinatorial study for  the pairwise interaction of    genus-one fibrations.
\end{abstract}

\maketitle
\tableofcontents

\section*{Introduction}
\mylabel{Introduction}

Which   smooth proper morphism $X\ra\Spec(\ZZ)$ do exist?  This tantalizing question seems to go  back to 
Grothendieck (\cite{Mazur 1986}, page 242), but it is rooted in algebraic number theory:
By   Minkowski's discriminant bound, for each number field $K\neq \QQ$ the corresponding
number ring  is not   \'etale over the ring $R=\ZZ$, so the only examples
in relative dimension $n=0$ are finite sums of $\Spec(\ZZ)$. Moreover, this reduces the general question  to  the
case $\Gamma(X,\O_X)=\ZZ$, such that the structure morphism  $X\ra\Spec(\ZZ)$ is a   contraction of fiber type.

A result of Ogg \cite{Ogg 1966}, which he himself attributes to Tate, asserts that each Weierstra\ss{} equation with integral 
coefficients has discriminant $\Delta\neq \pm 1$. 
Thus there is no family of elliptic curves
$E\ra\Spec(\ZZ)$. In other words, the Deligne--Mumford stack $\shM_{1,1}$ of smooth pointed curves of genus one
has  fiber category $\shM_{1,1}(\ZZ)=\varnothing$.
This was generalized by Fontaine \cite{Fontaine 1985},
who showed that there are no families $A\ra\Spec(\ZZ)$  of non-zero abelian varieties.
In other words, the  fiber category $\shM_{\text{Ab}}(\ZZ)$ for the  stack  
of   abelian varieties contains only  trivial objects.
This was independently established by Abrashkin \cite{Abrashkin 1985}, who also settled the case of abelian surfaces and   threefolds earlier
(\cite{Abrashkin 1976}, \cite{Abrashkin 1977}).

The theory   of relative Picard schemes \cite{FGA V} then shows that there is no family of smooth curves $C\ra\Spec(\ZZ)$
of genus $g\geq 2$. In other words, the Deligne--Mumford stack $\shM_g$ of smooth curves of genus $g$ has fiber category
$\shM_g(\ZZ)=\varnothing$.
Since the ring of integers has trivial Picard and Brauer groups, the only example in relative dimension $n=1$
are finite sums of the projective line $\PP^1_\ZZ$.
Abrashkin \cite{Abrashkin 1990} and Fontaine \cite{Fontaine 1993}   
also established that for each smooth proper $X\ra\Spec(\ZZ)$, the Hodge numbers of the complex fiber $V=X_\CC$
are very restricted, namely
$H^j(V,\Omega^i_{V/\CC}) = 0$ for $i+j\leq 3$ and $i\neq j$.
Note that this ensures that the   stack $\shM_{\text{HK}}$ of hyperk\"ahler varieties has fiber category 
$\shM_{\text{HK}}(\ZZ)=\varnothing$.

The only examples for smooth proper schemes over the integers that easily come to mind
are  the projective space $\PP^n$ and schemes   stemming from this,
for example sums, products,   blowing ups with respect to linear centers, or the Hilbert scheme 
$X=\Hilb^n_{\PP^2/\ZZ}$.
Similar examples arise from $G/P$, where $G$ is a Chevalley group  and $P$ is a  parabolic subgroup,
and from   torus embeddings $\operatorname{Temb}_\Delta$, stemming from a  fan of regular cones $\sigma\in\Delta$.
I am not aware of any other construction. It would be interesting to understand the case of   elliptic surfaces over $\PP^1$.

In light of the Enriques classification of algebraic surfaces, it is natural to consider  families
$Y\ra \Spec(\ZZ)$ of smooth proper surfaces with $c_1=0$. 
These are the abelian surfaces, bielliptic surfaces, K3 surfaces and Enriques surfaces.
They can be distinguished by their second Betti number, which takes the   values $b_2=6,2,22,10$.
The former two are ruled out by the Hodge numbers $h^{0,1}\geq 1$ of the complex fiber. Similarly, K3 surfaces do not 
occur, by $h^{2,0}=1$.
So only the Enriques surfaces remain, which over the complex numbers indeed have $h^{i,j}=0$ for $i+j\leq 3$ and $i\neq j$.

At first glance, the family of  K3 coverings $X\ra Y$ seems to show that there is no   family of Enriques surfaces.
There is, however, a  crucial loophole: At prime $p=2$, the K3 cover can degenerate into a \emph{K3-like covering}, and then
the fiber $X\otimes\FF_2$ becomes singular. In light of this, it actually seems    possible that
such families of Enriques surfaces exist. The main result of this paper asserts that it does not happen:

\begin{maintheorem}
{\rm (see Thm.\ \ref{no enriques over integers})}
The   stack $\shM_\Enr$ of Enriques surfaces
has fiber category $\shM_\Enr(\ZZ)=\varnothing$. In other words, 
there is no smooth proper family  of Enriques surfaces over the ring $R=\ZZ$. 
\end{maintheorem}

Note that the finiteness of the fiber categories $\shM_\Enr(R)$ for rings of integers $R$ in   number fields  was recently
established by Takamatsu \cite{Takamatsu 2019}. A similar finiteness holds for abelian varieties
by Faltings \cite{Faltings 1983} and Zarhin \cite{Zarhin 1985}, and also for K3 surfaces by Andr\'e \cite{Andre 1996}.
Results of this type are often refereed to as the \emph{Shafarevich Conjecture} \cite{Shafarevich 1963}

The proof of our non-existence result is long and indirect, and is given in the final Section \ref{Proof}.
It is actually a  consequence of the following:

\begin{maintheorem}
{\rm (see Thm.\ \ref{no non-exceptional})}
There is no Enriques surface $Y$ over   $k=\FF_2$  that is non-exceptional and has
constant Picard scheme $\Pic_{Y/k}=(\ZZ^{\oplus 10}\oplus\ZZ/2\ZZ)_k$.
\end{maintheorem}

The  \emph{exceptional Enriques surfaces} where introduced by Ekedahl and Shepherd-Barron \cite{Ekedahl; Shepherd-Barron 2004}.
They indeed can be discarded from our considerations, as I recently showed in \cite{Schroeer 2021b} that they do not admit a lifting
to the ring $W_2(\FF_2)=\ZZ/4\ZZ$.
The above result will be deduced from an explicit classification of \emph{geometrically rational elliptic surfaces} $\phi:J\ra\PP^1$
over $k=\FF_2$ that have constant Picard scheme, and satisfy 
certain additional technical conditions stemming from the theory of Enriques surfaces:

\begin{maintheorem}
{\rm (see Thm.\ \ref{classification non-zero j} and Thm.\ \ref{classification zero j})}
Up to isomorphisms, there are exactly eleven Weierstra\ss{}
equations $y^2+a_1xy + \ldots =x^3+a_2x^2+\ldots$ with coefficients $a_i\in \FF_2[t]$
that define a geometrically rational elliptic surface $\phi:J\ra\PP^1$  
with constant Picard scheme $\Pic_{J/\FF_2}$  having at most one rational point $a\in \PP^1$ where $J_a$ is semistable or supersingular.
\end{maintheorem}

For our Enriques surfaces $Y$ over the prime field $k=\FF_2$ with constant Picard scheme,
this narrows down the possible configurations of fibers over the three rational points $t=0,1,\infty$
in genus-one fibrations $\varphi:Y\ra\PP^1$, by passing to the jacobian fibration $\phi:J\ra\PP^1$ and using 
a result of Liu, Lorenzini and Raynaud \cite{Liu; Lorenzini; Raynaud 2005}.

It turns out that there are fifteen possible configurations, 
where the  additional four cases come from quasielliptic fibrations. In all but one situation there is precisely one ``non-reduced''
Kodaira symbol, from the list $\II^*,\III^*,\IV^*,\I^*_4, \I_2^*, \I_1^*$.
We then make an extensive combinatorial analysis for the possible interaction of pairs of genus-one fibrations
$\varphi,\psi:Y\ra\PP^1$ whose intersection number takes the minimal value  $\varphi^\inv(\infty)\cdot\psi^\inv(\infty)=4$,
in the spirit of Cossec and Dolgachev \cite{Cossec; Dolgachev 1989}, Chapter III, \S 5.
This only becomes feasible after  discarding the   exceptional Enriques surfaces,
which roughly speaking have ``too few'' fibrations and ``too large'' degenerate fibers.

The paper is organized as follows:
In Section \ref{Local system} we treat the group of numerical classes of invertible
sheaves as a local system $\Num_{X/k}$, that is, a sheaf on the \'etale site of the ground field $k$. 
Here the main result is that if it is constant,
then every projective contraction of the base-change to $k^\alg$ descends to a contraction of
$X$.
Geometric consequences for smooth surfaces $X=S$ are examined in more detail in Section \ref{Contractions}.
Then we analyze curves of canonical type in Section \ref{Curves}. They 
occur in genus-one fibrations, and we stress   arithmetical aspects.  
In Section \ref{Enriques surfaces} we review the theory of Enriques surfaces over algebraically closed fields, including
the notion of exceptional Enriques surfaces. Section \ref{Families} contains a detailed study of families  of Enriques surfaces,
in particular over the ring of integers.
In Section \ref{Geometrically rational} 
we collect basic facts on geometrically rational elliptic surfaces, again from an arithmetical point of view.
In Section \ref{Counting points}   we turn to counting points over finite fields with the Weil Conjectures,
and give some formulas that apply to Enriques surfaces and geometrically rational elliptic surfaces
whose local system $\Num_{X/k}$ is constant.
Section \ref{Multiple fibers} contains an observation on the behavior of reduced fibers in passing from
an elliptic fibration to its jacobian fibration.
Section \ref{Passage to jacobian} contains a main technical result: If $Y$ is an   Enriques surface over $k=\FF_2$
with constant Picard scheme, then the resulting geometrically rational surfaces $J$, arising as jacobian for
elliptic fibrations, also has constant Picard scheme.
Using this, we classify in Section \ref{Classification jacobian} all Weierstra\ss{} equations with coefficients from $\FF_2[t]$ that
possibly could describe such $J$. In Section \ref{Possible configurations}, we make a list of the possible fiber configurations
over the rational points $t=0,1,\infty$ that could arise for the Enriques surface $Y$,
where we also take into account quasielliptic fibrations.
This list is narrowed down to two configurations in Sections  \ref{Elimination I4*} and \ref{Elimination III*}.
In Section \ref{Restricted configurations} we   work again over arbitrary algebraically closed ground fields $k$, and show that
Enriques surfaces all whose genus-one fibrations are subject to certain severe combinatorial restrictions do not exist.
In the final Section \ref{Proof} we combine all these observations and deduce  our main results.

\begin{acknowledgement}
The research was supported by the Deutsche Forschungsgemeinschaft by the grant  \emph{SCHR 671/6-1 Enriques-Mannigfaltigkeiten}.
It was also conducted  in the framework of the   research training group
\emph{GRK 2240: Algebro-geometric Methods in Algebra, Arithmetic and Topology}.
I wish to thank Will Sawin for useful remarks, and Matthias Sch\"utt for helpful comments, in particular for 
pointing out  \cite{Shioda 1992} and \cite{Schuett; Shioda 2019},
which correct two entries in the classification of Mordell--Weil lattices for rational elliptic surfaces.
I would also like to thank the referees for very  careful and thorough reading and many valuable suggestions, which helped to improve the paper
and remove mistakes. In particular, for pointing out a defective argument in the first version concerning   two-sections in
what are  now Propositions \ref{no I4*+R} and   \ref{no I4*+R second}.
\end{acknowledgement}

\section{The local system of numerical classes}
\mylabel{Local system}

Throughout, $k$ denotes  a ground field of characteristic $p\geq 0$. Let $X$ be a proper scheme,   $\Pic_{X/k}$ the
Picard scheme, and $\Pic^\tau_{X/k}$ be the open subgroup scheme parameterizing numerically trivial
invertible sheaves  (see for example \cite{SGA 6}, Expos\'e XII, Corollary 1.5 and Expos\'e XIII, Theorem 4.7, compare also
\cite{Laurent; Schroeer 2021}, Section 2). 
The latter is   \emph{algebraic}, which means that the structure morphism is of finite type and separated,
and the resulting quotient
$$
\Num_{X/k}=(\Pic_{X/k})/(\Pic^\tau_{X/k})
$$
is an \'etale group scheme.
As such, it corresponds to the Galois representation on the stalk
$$
\Num_{X/k}(k^\alg)=\Num_{X/k}(k^\sep)=\ZZ^{\oplus\rho},
$$
where  $k^\sep\subset k^\alg$ are chosen separable and algebraic closures,
and  $\rho\geq 0$ denotes the Picard number for the base-change $\bX=X\otimes k^\alg$. 
The case $\rho=0$ is also allowed, it   may occur  for   schemes without ample sheaves (\cite{Schroeer 1999}, Section 3).
One may call  $\Num_{X/k}$  the \emph{local system  of numerical classes}.
Note that all these group schemes exist, even if $X$ fails to be reduced, irreducible,
or equidimensional.

We say that   $\Num_{X/k}$  is a  \emph{constant local system} if the corresponding Galois representation $\Gal(k^\sep/k)\ra\GL_\rho(\ZZ)$
is trivial. In other words, the group scheme  $\Num_{X/k}$  is isomorphic to  
$(\ZZ^{\oplus\rho})_k=\bigcup\Spec(k)$, where the disjoint union runs over all elements $l\in\ZZ^{\oplus\rho}$. 

In this section, we collect some noteworthy consequences of this condition. It can be solely characterized in terms
of invertible sheaves as follows:

\begin{lemma}
\mylabel{invertible sheaves}
The group scheme $\Num_{X/k}$ is constant if and only if there is an integer $d\geq 0$ with the following property:
For every invertible sheaf $\shL$ on $\bX$  there is a numerically trivial invertible sheaf $\shN$  on $\bX$
such that $\shL^{\otimes d}\otimes\shN$ is the pullback of some invertible sheaf on $X$.
\end{lemma}
 
\proof
The condition is sufficient: It ensures that the Galois action on the stalk  $\Num_{X/k}(k^\alg) =\ZZ^{\oplus\rho}$ becomes
trivial on some subgroup of finite index. Hence it is trivial after tensoring with $\QQ$, and thus
on the torsion-free group $\ZZ^{\oplus\rho}$.

The condition is also necessary:
Each rational point $l\in \Num_{X/k}$ can be seen as a representable  torsor
$T\subset\Pic_{X/k}$ for the algebraic group scheme $P=\Pic^\tau_{X/k}$.
Fix a closed point $a\in T$. The torsor becomes trivial after base-changing to the finite extension field $E=\kappa(a)$.
So we may regard the isomorphism class of $T$ as an element in $l\in H^1(k,P)$, where    the cohomology is taken
with respect to the finite flat topology. 

Such classes have  finite order, which follows from general principles:
Set $X'=X\otimes E$ and $P'=\Pic_{X'/E}$, and view them as $k$-schemes. First note that for each $k$-algebra $A$,
the projection $\pr:X'\otimes A\ra X\otimes A$ is locally free of rank $n=[E:k]$,
so  besides the pull-back   $\pr^*:\Pic(X\otimes A)\ra\Pic(X'\otimes A)$ we also have a \emph{norm map} $N$ in the other direction,
where $N\circ\pr^*$ is multiplication by $n$ (for details see \cite{EGA II}, Section 6.5).
Now consider the    diagram 
$$
\begin{tikzcd} 
 	& 			& H^1(k,P)\ar[d,"\pr^*"]\ar[dl,"\pr^*"]\\
0\ar[r]	& H^1(k,\pr_*P')\ar[r]\ar[ur,bend left=20,"N"]	& H^1(E,P')\ar[r]		& H^0(k,R^1\pr_*P').
\end{tikzcd}
$$
Here the lower   sequence is exact, coming  from the Leray--Serre spectral sequence.
Note that the term on the right is zero, because  the coefficient sheaf 
$R^1\pr_*P'$  vanishes, by the arguments for Lemma \ref{brauer group} below. 
Our class $l\in H^1(k,P)$ becomes trivial
in $H^1(E,P')$, hence $n \cdot l=N(\pr^*(l))=0$.
 
Summing up,  there is some $n\geq 1$ such that $nl$ comes from a rational point $l'\in \Pic_{X/k}$.
We next show that some multiple of $l'$ comes from an invertible sheaf. Set $S=\Spec(k)$.
The Leray--Serre spectral sequence for the structure morphism
$h:X\ra S$ yields an exact sequence
$$
H^1(X,\GG_m)\lra H^0(S,R^1h_*(\GG_m))\lra H^2(S,h_*(\GG_{m,X})).
$$
The arrow on the left coincides with the canonical map $\Pic(X)\ra \Pic_{X/k}(k)$. The 
term on the right is a torsion group by Lemma \ref{brauer group} below.
Consequently there is some integer $n'\geq 1$ so that the point  $n'l'$ comes
from some invertible sheaf $\shL$ on $X$. 
Applying the above reasoning for the members of  some generating set $l_1,\ldots, l_r\in\Num_{X/k}(k)$,
we conclude that the composite map $\Pic(X)\ra\Num_{X/k}(k^\alg)$ has finite cokernel, and the assertion follows.
\qed

\medskip
Attached to any noetherian scheme $Z$ is the \emph{dual graph} $\Gamma=\Gamma(Z)$. Its 
finitely many vertices $v_i\in\Gamma$ correspond to the irreducible components $Z_i\subset Z$, and two vertices $v_i\neq v_j$ are joined by an edge if $Z_i\cap Z_j\neq\varnothing$. Any morphism $Z'\ra Z$ that sends irreducible components to irreducible components,
in particular flat maps,  induces a map $\Gamma(Z')\ra\Gamma(Z)$ between vertex sets that sends edges to edges.
Such maps are termed \emph{graph morphisms}. They are called  \emph{graph isomorphism}
if the maps on vertex and edge sets are bijective.

\begin{proposition}
\mylabel{graph isomorphism}
Suppose   $X$ is one-dimensional. Then the local system  $\Num_{X/k}$ is constant if and only if 
 $\Gamma(\bX)\ra\Gamma(X)$ is a graph isomorphism.
\end{proposition}

\proof
Let $C_1,\ldots,C_r\subset X$ be the irreducible components, endowed with reduced scheme structure.
First note that an  invertible sheaf $\shL$ is numerically trivial if and only if $(\shL\cdot C_i)=0$ for $1\leq i\leq r$.

We now show that our condition is sufficient: If $\Gamma(\bX)\ra\Gamma(X)$ is a graph isomorphism,
the stalk $\Num_{X/k}(k^\alg)$ has rank $\rho=r$.
For each $1\leq i\leq r$, choose a closed point $a_i\in X$ not contained in $\Ass(\O_X)$
and the union of the $C_j$, $j\neq i$, and some effective Cartier divisor $D_i\subset X$ supported by $a_i$.
Recall that $\Ass(\O_X)$ is the set of all points where the local ring $\O_{X,\zeta}$ admits a copy of the
residue field $\kappa(\zeta)$ as a submodule. Here it  comprises the generic points, together with closed points that give
embedded components.
The   invertible sheaves $\shL_i=\O_X(D_i)$ define elements $l_i\in\Num_{X/k}(k^\alg)$ that
generate a subgroup of finite index. In turn, the Galois representation must be trivial.

The condition is also necessary: Let $E_1,\ldots,E_s\subset\bX$ be the irreducible components.  
Fix some $1\leq i_0\leq s$, and choose some effective 
Cartier divisor  on $X$ supported by $E_{i_0}$ but disjoint from the other components.
Using  Lemma \ref{invertible sheaves}, we find  some invertible sheaf $\shL$ on $X$ with $(\shL_\bX\cdot E_i)\geq 0$,
with inequality if and only if $i=i_0$.
It follows that $(\shL\cdot  C_j)\geq 0$, with inequality for precisely one index $j=j_0$,
such that $C_{j_0}\otimes k^\alg=E_{i_0}$ as closed sets.
Consequently $\Gamma(\bX)\ra\Gamma(X)$ is bijective, whence a graph isomorphism.
\qed

\medskip
Let us use the following terminology: A \emph{contraction} of $X$ is a proper scheme $Z$ together with a   morphism
$f:X\ra Z$ with $\O_Z=f_*(\O_X)$.

\begin{theorem}
\mylabel{projective contractions}
Suppose  the local system $\Num_{X/k}$ is constant.  
Then every contraction $\bX\ra \bY$  to some projective scheme $\bY$ is isomorphic to the base-change of a contraction $X\ra Y$.
\end{theorem}

\proof
First note that each  invertible sheaf $\shL$ on $X$ yields the graded ring 
$$
R(X,\shL)=\bigoplus_{t=0}^\infty H^0(X,\shL^{\otimes t}),
$$
and the resulting homogeneous spectrum $P(X,\shL)=\Proj R(X,\shL)$. If $\shL$ is semiample,
then $Z=P(X,\shL)$ is projective and yields a contraction, and each   contraction with $Z$ projective is of this form.
The sheaf  $\shL$ is unique up to preimages of semiample sheaves on $Z$.

According to \cite{EGA IVc}, Theorem 8.8.2 there is a finite   extension $k\subset k'$
such that the given contraction $\bX\ra \bY$ is the base-change of some contraction $f':X'\ra Y'$ 
with $Y'$ projective, where $X'=X\otimes k'$. Choose some semiample invertible sheaf $\shL'$
on $X'$ defining the contraction, as discussed in the preceding paragraph.
Enlarging $k'$, we may assume that it is a splitting field. Then it is the composition
of a Galois extension by some purely inseparable extension, and it suffices to treat these cases
individually.

Suppose first that $k'$ is Galois, with   Galois group $G=\Gal(k'/k)$.  
For each $\sigma\in G$, the pullback $\shL'_\sigma=\sigma^*(\shL')$ is semiample, and also numerically equivalent to $\shL'$
because $\Num_{X/k}$ is constant.
Replacing $\shL'$ by the tensor product $\bigotimes_\sigma\shL'_\sigma$, we may assume that the $k'$-valued
point $\Pic_{X/k}(k')$ corresponding to $\shL'$ is Galois-invariant, hence defines a rational point $l\in\Pic_{X/k}$.
Passing to a multiple, it comes from an invertible sheaf $\shL$ on $X$, and the assertion follows.

It remains to treat the case that $k'$ is purely inseparable.
In this setting we can prove a    stronger statement:
Any contraction $X'\ra Y'$ to a \emph{proper}  scheme $Y'$ is the base-change of some contraction $X\ra Y$.
To see this, recall that our scheme $X=(|X|,\O_X)$ comprises an underlying topological space and a structure sheaf,
and that the given contraction is a pair $(f',\varphi')$ where $f:|X'|\ra |Y'| $ is a continuous map  and
$\varphi':\O_{Y'}\ra\O_{X'}$ is an $f'$-homomorphism (\cite{EGA I}, Chapter 0, Section 3.5.1).
Since $\pr:X'\ra X$ is a homeomorphism, we are forced to set  $|Y|=|Y'|$ and $f=f'$  and   $\O_Y=f_*(\O_X)$.
Define $Y=(|Y|,\O_Y)$ and write  $\varphi:\O_Y\ra\O_X$ for the canonical $f$-homomorphism.
This gives a morphism   of ringed spaces $(f,\varphi):X\ra Y$.

It remains to verify that this is the desired contraction. Obviously, the topological space $|Y|$ is quasicompact, and $f_*(\O_X)=\O_Y$.
The rest of the  problem is essentially local: Let $V'\subset Y'$ be an affine open subscheme,  $U'=f'^{-1}(V)$ the preimage, 
and $U\subset X$ the corresponding open subscheme. 
The condition $f'_*(\O_{X'})=\O_{Y'}$ ensures that    
$V'=\Spec\Gamma(U',\O_{X'})$. In other words, $U'\ra V'$ is the \emph{affine hull}, and this morphism is proper.
Consider the affine hull $V=\Spec\Gamma(U,\O_X)$. Then $U'\ra V'$ is the base-change of $U\ra V$, according to \cite{EGA IIIa}, Proposition 1.4.15.
Using that $V'\ra V$ is a homeomorphism, we get an identification $(|V|,\O_Y|_{|V|})=V$, compatible with the morphisms from $U$.
It follows that the ringed space $Y=(|Y|,\O_Y)$
is a scheme and that $f:X\ra Y$ is a morphisms of schemes.  Moreover, we infer $Y'=Y\otimes k'$ and $f'=f\otimes k'$.
Applying \cite{EGA IVb}, Proposition 2.7.1 to the structure morphism $Y\ra \Spec(k)$, we see that $Y$ is proper.
\qed

\medskip
Note that the    part in the  previous proof that  deals with purely inseparable extensions $k\subset k'$ 
works without any assumptions on $\Num_{X/k}$. 
The  proof of Lemma \ref{invertible sheaves} relies on  the following observation:

\begin{lemma}
\mylabel{brauer group}
Write $S=\Spec(k)$, and let $h:X\ra S$ be the structure morphism. Then we have an identification $H^2(S,h_*(\GG_{m,X}))=\Br(X^\aff)$,
and this is a torsion group.
\end{lemma}

\proof
For any $k$-algebra $A$, we have $\Gamma(X\otimes A,\GG_a)=\Gamma(X^\aff\otimes A,\GG_a)$,
because $X$ is   quasicompact and separated, and $k\ra A$ is flat.  It follows that the multiplicative groups
on $X$ and $X^\aff$ have the same direct image on $S$, and it suffices to treat the case that $X=\Spec(R)$ is finite.
One easily checks that the restriction map $H^2(X,\GG_m)\ra H^2(X_\red,\GG_m)$ is bijective.
The term on the right is the sum of Brauer groups for  fields, which is torsion.

Consider the Leray--Serre spectral sequence $H^i(S,R^jh_*(\GG_{m}))\Rightarrow H^{i+j}(X,\GG_m)$.
We show that the edge map $H^2(X,\GG_m)\ra H^2(S,h_*(\GG_{m}))$ is bijective by checking that 
  $R^ih_*(\GG_{m})=0$ for $i=1,2$.
Let $A$ be a finite $k$-algebra. Then $X\otimes A$ is semilocal, hence its Picard group vanishes, and we immediately see the  
vanishing for $i=1$. Now suppose we have a  class $\alpha\in H^2(X\otimes A,\GG_m)$.
Let $k^\sep$ be some separable closure. The base-change $A\otimes_kk^\sep$ is the product of finitely many 
strictly henselian local Artin rings, so the pullback of $\alpha$ vanishes. Then there is already
a finite extension $k\subset k'$ on which the pullback vanishes. This shows $R^2h_*(\GG_m)=0$.
\qed

\section{Contractions of surfaces}
\mylabel{Contractions}

We now examine in dimension two the results of the preceding section.
Let $k$ be a ground field of characteristic $p\geq 0$, and 
 $S$ be a smooth proper surface with $h^0(\O_S)=1$.
Let $C\subset X$ be a curve. Decompose it into irreducible components  $C=m_1C_1+\ldots+m_rC_r$ with multiplicities $m_i\geq 1$, 
and let 
$$
N=N(C)=(C_i\cdot C_j)_{1\leq i,j\leq r}\in\Mat_r(\ZZ)
$$
be the resulting intersection matrix.  
We say that the curve $C$ is \emph{negative-definite}  or \emph{negative-semidefinite} if the
intersection matrix $N$ has the respective property. 
One may use the intersection numbers  to endow $\Gamma=\Gamma(C)$ with the
structure of an \emph{edge-labeled graph}, where the   labels attached to the edges are the intersection numbers
$(C_i\cdot C_j)>0$ for $i\neq j$

Now let  $f:S\ra Z$ be a
contraction. Then  either $Z$ is a geometrically normal   surface and the morphism is birational,
or $Z$ is a smooth curve and the morphism is of fiber type.
Let $\ba:\Spec(\Omega)\ra Z$ be a geometric point, for some algebraically closed field $\Omega$,
with image point $a\in Z$.
Write $S_\ba=X\times_Z\Spec(\Omega)=f^{-1}(\ba)$ and $S_a=S\times_Z\Spec\kappa(a)=f^{-1}(a)$ for the resulting geometric 
and schematic fibers, respectively. 
The goal of this section is to establish the following:

\begin{theorem}
\mylabel{geometric and schematic fiber}
Suppose the geometric fiber $S_\ba$   is one-dimensional.
If the group scheme $\Num_{S/k}$ is constant,   the following holds:
\begin{enumerate}
\item 
The map $\Gamma(S_\ba)\ra\Gamma(S_a)$ is a graph isomorphism.
\item 
The above graph isomorphism respects edge labels, provided that the  reduced scheme  $(S_a)_\red$ 
is geometrically reduced  over the residue field $\kappa(a)$.
\item 
If the geometric fiber $S_\ba$ is reducible, the residue field extension $k\subset\kappa(a)$ is purely inseparable.
\end{enumerate}
\end{theorem}

\proof
Assertion (ii) is an immediate consequence from (i) and the definitions, and left to the reader.
To verify (i) and (iii) we write  $E=\kappa(a)$. Without restriction, it suffices to treat the case that $\Omega= E^\alg$.
The geometric fiber $S_\ba=S_a\otimes_E\Omega\subset S_a\otimes_k\Omega$ is a   curve on the base-change $S\otimes_k\Omega$.
This curve is connected, by the condition $\O_Z=f_*(\O_S)$.
Let  $C_i\subset S_\ba$, $1\leq i\leq r$ be the irreducible components, and write $N=(C_i\cdot C_j)_{1\leq i,j\leq r}$
for the resulting intersection matrix.

We first consider   the case that $f:X\ra Z$ is birational. Then $N$ is negative-definite, in particular
the canonical map $\bigoplus\ZZ C_i\ra \Num_{S/k}(\Omega)$ is injective. It follows that the Galois action
on the set $\{C_1,\ldots,C_r\}$ is trivial. In turn, $\Gamma(S_\ba)\ra\Gamma(S_a)$ is a graph isomorphism.
Now suppose that $r\geq 2$. Seeking a contradiction, we assume that $k\subset E$ is not purely inseparable.
Then $S_a\otimes_k\Omega$ is disconnected. On the other hand, the fiber $S_a$ is connected,
 by the condition $\O_Z=f_*(\O_S)$.
It follows that the Galois action on the   irreducible components $E_i\subset S_a\otimes_k\Omega$
is non-trivial.
As above, the canonical map $\bigoplus\ZZ E_i\ra\Num_{S/k}(\Omega)$ is injective, so the Galois action on 
the set $\{E_1,\ldots,E_s\}$ must be trivial, contradiction.

It remains to treat the case that  $f:X\ra Z $ is of fiber type. Then the matrix $N$ is negative-semidefinite.
For $r=1$ our graphs have but one vertex, and the map $\Gamma(S_\ba)\ra\Gamma(S_a)$ is obviously a graph isomorphism.
Now suppose  $r\geq 2$. Then  the curves $C_1+\ldots+C_{r-1}$ and $C_2+\ldots+C_r$ are negative-definite. Using the 
previous paragraph,
we infer that $\Gamma(S_\ba)\ra\Gamma(S_a)$ is again a graph isomorphism. Likewise, one shows that
 $k\subset E$ is purely inseparable.
\qed

\section{Curves of canonical  type}
\mylabel{Curves}

Let $k$ be a ground field of characteristic $p\geq 0$, and $S$ be a smooth proper surface with $h^0(\O_S)=1$.
In this section we apply Theorem \ref{geometric and schematic fiber} to fibers in genus-one fibrations,
and collect additional facts for the situation that $k=\FF_q$ is finite.

We start with  some useful general terminology:
Let $C\subset S$ be a curve, and  $C=m_0C_0+\ldots+m_rC_r$ be the decomposition into irreducible components, together with
their multiplicities.
We say that $C$ is of \emph{fiber type} if  
$(C\cdot C_i)= 0$ for all indices  $0\leq i\leq r$.
Obviously, this condition transfers to multiples and connected components of the curve.
If $C$ is connected and $\gcd(m_0,\ldots,m_r)=1$, we say that
$C$ is an \emph{indecomposable curve of fiber type}. 
Then the intersection matrix $N=(C_i\cdot C_j)_{1\leq i,j\leq r}$ is negative-semidefinite,
and   the radical of
the intersection form on the lattice $L=\bigoplus_{i=0}^r \ZZ C_i$ is generated by $m_0C_0+\ldots+m_rC_r$.

Given a connected curve  of fiber type $C=\sum m_iC_i$, we call $m=\gcd(m_0,\ldots,m_r)$ its \emph{multiplicity}.
We say that $C$ is \emph{simple} if $m=1$, and \emph{multiple} otherwise.
Setting $n_i=m_i/m$, we call  $C_\ind=\sum n_iC_i$ the \emph{underlying indecomposable curve of fiber type}.

Suppose now that we have a contraction $f:S\ra B$ of fiber type onto some curve $B$.
For each geometric point $\ba:\Spec(\Omega)\ra B$, the geometric fiber $S_\ba=f^{-1}(\ba)$ is a connected 
curve of fiber type on the base-change $S\otimes_k\Omega$. It is useful to write $f^\inv_\ind(\ba)$ for the underlying
indecomposable curve of fiber type, and $f^\inv_\red(\ba)$ for its reduction.
An analogous notation will be used for  the schematic fiber $S_a=f^\inv(a)$
over the image point $a\in B$.
 
We are mainly interested in \emph{curves of canonical type}. By definition, this means 
$(C\cdot C_i)=(K_S\cdot C_i)=0$ for all indices $0\leq i\leq r$. The above locutions for curves of fiber type
are used in an analogous way for  curves of canonical types.

Suppose now that $f:S\ra B$ is a genus-one fibration that is relatively minimal.
Then every geometric fiber $ f^{-1}(\ba)$ is a connected curve of canonical type.
The underlying indecomposable curves of canonical type $f^\inv_\ind(\ba)$ were classified by N\'eron \cite{Neron 1964} and Kodaira \cite{Kodaira 1963}.
Throughout, we denote the \emph{Kodaira symbols} 
\begin{equation}
\label{kodaira symbols}
\I_m,\;\I_n^*,\;\II,\;\III,\;\IV,\;\IV^*,\;\III^*,\; \II^*
\end{equation}
to designate the geometric fibers.  The symbols $\I_m$, $m\geq 1$ are also called \emph{semistable} or \emph{multiplicative},
whereas $\I_n^*,\II, \ldots,\II^*$,  are referred to as \emph{unstable} or \emph{additive}. Here $n\geq 0$.
As customary, the symbol $\I_0$ indicates that $f^\inv_\ind(\ba)$ is an elliptic curve; we then distinguish the cases
that the  smooth curve is \emph{ordinary} or \emph{supersingular}.
By abuse of terminology 
we say that a geometric fiber is \emph{reducible or singular} if the scheme $f^\inv_\ind(\ba)$ has the respective property.
If this curve is  singular with  $f$  elliptic, or reducible with $f$ quasielliptic, or multiple, 
 we say that the geometric fiber is \emph{degenerate}.

The   classification of geometric fibers $f^\inv(\ba)$ was extended to schematic fibers $ f^{-1}(a)$ by Szydlo \cite{Szydlo 2004},
where the situation becomes considerably more involved.
If all irreducible components   of $S_a=f^\inv(a)$ 
are birational to $\PP^1$, the map $\Gamma(S_\ba)\ra\Gamma(S_a)$ is a graph isomorphism respecting edge labels,
and we also use the Kodaira symbols in \eqref{kodaira symbols} to designate   $f^{-1}_\ind(a)=m_1C_1+\ldots+m_rC_r$.
We need, however, the following modifications for $r\leq 2$, which play an important role in our later applications:

For $r=1$ the curve  $C=C_1$ has $h^0(\O_C)=h^1(\O_C)=1$, and the normalization $\nu:\PP^1\ra C$ is described by 
a  \emph{conductor square}
$$
\begin{CD}
\Spec(R)	@>>>	\PP^1\\
@VVV		@VV\nu V\\
\Spec(k)	@>>>	C,
\end{CD}
$$
which is both cartesian and cocartesian, see \cite{Fanelli; Schroeer 2020}, Appendix A for more details. 
Geometrically speaking, the closed subscheme $\Spec(R)\subset\PP^1$
is contracted to a rational point $\Spec(k)\subset C$.
Here $R$ is a finite $k$-algebra of length two. If $R=k\times k$ or $R=k[\epsilon]$ with $\epsilon^2=0$ 
then $C$ is the rational nodal curve or the rational cuspidal curve and we use the standard Kodaira symbol $\I_1$ and $\II$,
respectively. If  $R=E$ is a field, then $C$ is a  twisted form, 
and we use   \emph{twisted Kodaira symbols} $\tI_1$ and $\widetilde{\II}$ instead.
Note that in the twisted situation, a non-rational point on $\PP^1$ gets identified
with a rational point on $C$, with consequences for the size of $C(k)$.

A similar situation appears for $r=2$: Then   $C=C_1+C_2$ sits in the conductor square
$$
\begin{CD}
\Spec(R\times R)	@>>>	\PP^1\amalg\PP^1\\
@VVV		@VV\nu V\\
\Spec(R)	@>>>	C,
\end{CD}
$$
where $R$ is a $k$-algebra of length two. For $R=k\times k$ or $R=k[\epsilon]$, where $\epsilon^2=0$, we again use the standard Kodaira symbols
$\I_2$ and $\III$, respectively. If $R=E$ is a field extension, we have twisted forms and  use $\tI_2$ and $\widetilde{\III}$ instead.
Again, the twisted situation has   consequence for the size of $C(k)$.

Note that over perfect ground fields, and in particular over finite ground fields, for $r\leq 2$ the cases $\widetilde{\II}$ and $\widetilde{\III}$ do not occur,
and   we only have $\tI_1$ and $\tI_2$ as twisted Kodaira symbols. 
One then says that $C$ is a \emph{non-split semistable fiber}.  
For $r\geq 3$ there are further   possible twists, but then  $\Gamma(S_\ba)\ra\Gamma(S_a)$ is not a graph isomorphism.
We   record the following consequence of Theorem \ref{geometric and schematic fiber}:

\begin{proposition}
\mylabel{constant num over finite field}
Suppose the ground field $k=\FF_q$ is finite and
that  the  group scheme $\Num_{S/k}$ is constant.
Then $\Gamma(S_\ba)\ra\Gamma(S_a)$ is a graph isomorphism respecting edge labels, and 
the implications (i)$\Rightarrow$(ii)$\Leftrightarrow$(iii) hold  among the following three conditions:
\begin{enumerate}
\item The geometric fiber $S_\ba$ is reducible.
\item The closed point $a\in B$ is a rational point.
\item The irreducible components of the schematic fiber $S_a$ are birational to $\PP^1$.
\end{enumerate}
\end{proposition}

\proof
The   assertion on $\Gamma(S_\ba)\ra\Gamma(S_a)$ immediately follows from Theorem \ref{geometric and schematic fiber},
which also ensures that   $k=\kappa(a)$ provided that $S_a$ is reducible.
This already gives (i)$\Rightarrow$(ii). The implication (ii)$\Leftarrow$(iii) is trivial.
For the converse,  let $C\subset S_a$ be an irreducible component, and $C'\ra C$ be its normalization.
Then $C\otimes\Omega$ is reduced, with normalization $C'\otimes\Omega$, because the field extension $k\subset\Omega$
must be separable. The curve $C\otimes\Omega$ is  birational to $\PP^1_\Omega$, because 
$\Gamma(S_\ba)\ra\Gamma(S_a)$ is a graph isomorphism. Hence the same holds for $C'\otimes\Omega$.
In other words, $C'$ is a one-dimensional Brauer--Severi variety.  
By Wedderburn's Theorem, we must have $C'\simeq \PP^1$.
\qed

\begin{proposition}
\mylabel{canonical type rational points}
Assumptions as in Proposition \ref{constant num over finite field}.
Let $C=f^\inv(a)$ be a sche\-matic  fiber over a rational point $a\in\PP^1$, and $r\geq 1$ be the number
of irreducible components. Suppose the reduced scheme $C_\red$ is singular. Then the number  of rational points in the fiber 
is  given by the following table:
$$
\begin{array}[b]{@{}lllll}
\toprule
\text{\rm Kodaira symbol}	& \I_r	& \tI_1	&\tI_2	& \text{\rm unstable}\\
\midrule
\Card C(\FF_q)				& rq	& q+2 	& 2q+2	& rq+1 \\
\bottomrule
\end{array}
$$
For the field $k=\FF_2$ with two elements, the  Kodaira symbol $\I_0^*$ is impossible.
\end{proposition}

\proof
The computation of $n\geq 0$ follows directly by comparing $C=f^\inv_\red(a)$ with its normalization $C'=\PP^1\cup\ldots\cup\PP^1$, 
and is left to the  reader. If the Kodaira symbol is $\I_0^*$, then in $C=C_0+C_1+\ldots+C_4$
the terminal components $C_1,\ldots,C_4$ intersect the central component $C_0=\PP^1$ in four rational points.
This implies $4\leq \Card\PP^1(\FF_q)=q+1$, which means $q\neq 2$.
\qed

\medskip
We now   examine the situation that $k=\FF_2$  and  the scheme $ f^\inv_\ind(a)$ is smooth. 
This genus-one curve contains a rational point (\cite{Lang 1956}, Theorem 2), hence
can be regarded as an elliptic curve $E=f^\inv_\ind(a)$. 
Let us recall from \cite{Husemoller 1987}, Chapter 3, \S6 
that up to isomorphism, there are five elliptic curves $E_1,\ldots,E_5$ over the prime field $k=\FF_2$.
The groups $E_i(k)$ are cyclic, and all   $1\leq n\leq 5$ occur as orders. Deviating from loc.\ cit., we
choose indices in a canonical way according to the order of the group of rational points, and record:

\begin{proposition}
\mylabel{elliptic curves}
Up to isomorphisms, the elliptic curves $E_i$ over the field $k=\FF_2$ with $\Card E_i(\FF_2)=i$
are given by the following table:
$$
\begin{array}[b]{@{}lll}
\toprule
\text{elliptic curve}	& \text{Weierstra\ss{} equation}	& \text{invariant}\\
\toprule
E_1			& y^2+y=x^3+x^2+1 			& \text{supersingular with $j=0$}\\
E_3			& y^2+y=x^3\\
E_5			& y^2+y=x^3+x^2\\
\midrule
E_2			& y^2+xy=x^3+x^2+x  			& \text{ordinary with $j=1$}\\
E_4			& y^2+xy=x^3+x\\
\bottomrule
\end{array}
$$
\end{proposition}

\section{Enriques surfaces and exceptional Enriques surfaces}
\mylabel{Enriques surfaces}

In this section we collect the relevant facts on Enriques surfaces over algebraically closed fields. For general information,
we refer to the monograph of Cossec and Dolgachev \cite{Cossec; Dolgachev 1989}. Here we
are particularly interested in the so-called exceptional Enriques surfaces, which were
introduced and studied by Ekedahl and Shepherd-Barron \cite{Ekedahl; Shepherd-Barron 2004}.

Throughout, we fix an algebraically closed ground field $k$ of characteristic $p\geq 0$.
An \emph{Enriques surface}  is a smooth proper  scheme $Y$ with 
$$
h^0(\O_Y)=1\quadand  c_1=0\quadand b_2=10.
$$
The first condition means that $Y$ is connected, and  the second condition
 signifies that the dualizing sheaf $\omega_Y$ is numerically trivial.
In the Enriques classification of algebraic surfaces, Bombieri and Mumford  showed in  \cite{Bombieri; Mumford 1976}, \S 3    that
the   Picard scheme $\Pic^\tau_{Y/k}$ of numerically trivial sheaves is a finite group scheme
of order two, whose group of rational points is generated by the dualizing sheaf $\omega_Y$.
Moreover, $\Num(Y)$ is a  free abelian group of rank $\rho=10$, the Betti numbers satisfies $b_1=b_3=0$ and $b_2=\rho$. 
Furthermore, $h^1(\O_Y)=h^2(\O_Y)$, and this number is at most one.
The intersection form on $\Num(Y)$  has signature $(1,9)$ by the Hodge Index Theorem, must be  even by Riemann--Roch, and is actually 
unimodular (\cite{Illusie 1979}, Corollary 7.3.7). By the classification of indefinite unimodular forms (\cite{Milnor; Husemoller 1973}, Chapter II, Theorem 5.3), 
we have $\Num(Y)\simeq U\oplus E_8$, where the first summand has Gram matrix 
$(\begin{smallmatrix}0&1\\1&0\end{smallmatrix})$  and the second summand is the negative-definite $E_8$-lattice.
Each curve of canonical type $C\subset Y$ induces a genus-one fibration $\varphi:Y\ra\PP^1$, and
this  leads to the following fact (\cite{Bombieri; Mumford 1977}, Theorem 3 and Proposition 11, combined with
 \cite{Lang 1983}, Theorem 2.2):

\begin{proposition}
There is at least one genus-one fibration $\varphi:Y\ra\PP^1$.
Any such fibration has a multiple fiber $mF$ of multiplicity $m=2$,
and there is an integral curve $C\subset Y$ with $C\cdot \varphi^\inv(\infty)=2$.
Such a two-section can be chosen to be  a $(-2)$-curve, or a half-fiber in some other genus-one fibration $\psi:Y\ra\PP^1$.
\end{proposition}

Here a curve $F\subset Y$ is called a \emph{half-fiber} if $2F$ is a fiber in some
genus-one fibration $\varphi:Y\ra\PP^1$ over some rational point.
Another important fact (\cite{Bombieri; Mumford 1977},  Proposition 11):

\begin{proposition}
\mylabel{omega order two}
The following four conditions are equivalent:
\begin{enumerate}
\item The dualizing sheaf $\omega_Y$ has order two in the Picard group.
\item The cohomology group $H^1(Y,\O_Y)$ vanishes.
\item There is a genus-one fibration $\varphi:Y\ra\PP^1$ without wild fibers.
\item There is a genus-one fibration with two multiple fibers $2F_1$ and $2F_2$.
\end{enumerate}
Under these conditions,   (iii) and (iv) are valid for all genus-one fibrations, and the dualizing sheaf
is given by $\omega_Y\simeq\O_Y(F_1-F_2)$. Moreover, the above conditions are automatic for $p\neq 2$.
\end{proposition}

\begin{lemma}
\mylabel{multiple for p=2}
Let $\varphi:Y\ra\PP^1$ be a genus-one fibration with two multiple fibers $2F_1\neq 2F_2$
in characteristic $p=2$. Then each $ F_i$ is either an ordinary elliptic curve  or unstable.
\end{lemma}

\proof
Write $C=F_i$.
The   multiple fibers   must be tame, hence the   sheaf $\shN=\O_C(C)$ has order $m=2$ in
 $\Pic(C)$. This Picard group has no element of order two if $C$ is a supersingular elliptic curve.
If $C$ is semistable, then $\Pic^0(C)=\GG_m(k)=k^\times $ likewise has no elements of order two. 
\qed
 
\medskip
Suppose now that we are in characteristic $p=2$. Up to isomorphism, there are three group schemes of order two,
namely $\mu_2$ and $\ZZ/2\ZZ$ and $\alpha_2$. In turn, there are three types of Enriques surfaces $Y$, which are called
\emph{ordinary}, \emph{classical} or \emph{supersingular}, according to the following table:
$$
\begin{array}[c]{@{}lll l}
\toprule	
\text{\rm group scheme $P=\Pic^\tau_{Y/k}$}		& \mu_2					& \ZZ/2\ZZ					& \alpha_2\\
\midrule
\text{\rm Cartier dual $G=\uHom(P,\GG_m)$}		& \ZZ/2\ZZ				& \mu_2						& \alpha_2\\
\midrule
\text{\rm fundamental group $\pi_1(Y,\Omega)$}	& \ZZ/2\ZZ				& {\rm trivial}				& {\rm trivial}	\\
\midrule
\text{\rm dualizing sheaf $\omega_Y$}			& \O_Y					& \O_Y(F_1-F_2) 			& \O_Y\\
\midrule
\text{designation of $Y$}						& \text{\rm ordinary}	& \text{\rm classical} 	& \text{\rm supersingular}\\
\bottomrule
\end{array}
$$

In case that $P=\Pic^\tau_{Y/k}$ is unipotent, in other words, the Cartier dual $G=\uHom(P,\GG_m)$ is local,
one says that $Y$ is \emph{simply-connected}.
The canonical inclusion of $P$ into the Picard scheme corresponds to a non-trivial $G$-torsor $\epsilon:X\ra Y$,
compare the discussion in \cite{Schroeer 2021a}, Section 4.
If the Enriques surface $Y$ is simply-connected  then $\epsilon:X\ra Y$ is a universal homeomorphism, and  $X$ is called
the \emph{K3-like covering}. Indeed, the scheme $X$ is an integral proper surface that is Cohen--Macaulay 
with $\omega_X=\O_X$ and $h^0(\O_X)=h^2(\O_X)=1$ and $h^1(\O_X)=0$.
However, the singular locus $\Sing(X)$ is non-empty.

The simply-connected Enriques surface can be divided into two   subclasses depending on
properties of the normalization $\nu:X'\ra X$. The induced projection $X'\ra Y$ is purely inseparable and flat of degree two,
hence the ramification locus for $\nu$ is an effective Cartier divisor $R'\subset X'$, giving $\omega_{X'}=\O_{X'}(-R')$.
Ekedahl and Shepherd-Barron \cite{Ekedahl; Shepherd-Barron 2004} 
observed that $R'$ is the preimage of some effective Cartier divisor $C\subset Y$
referred to as the  \emph{conductrix}, and  $2C$ is called the \emph{biconductrix}.
They  called an Enriques surface \emph{exceptional} if it is simply-connected and 
the biconductrix has $h^1(\O_{2C})\neq 0$. From this amazing definition, they obtained a complete classification
of exceptional Enriques surfaces. Moreover, Salomonsson \cite{Salomonsson 2003} described birational models by explicit equations.
The following property of  \emph{non-exceptional Enriques surfaces} will play an important role later:

\begin{proposition}
\mylabel{properties non-exceptional}
Suppose $p=2$, and that the Enriques surface $Y$ is non-except\-ional. Then the following holds:
\begin{enumerate}
\item 
There is no quasielliptic fibration $\varphi:Y\ra\PP^1$ having a simple fiber  
with Kodaira symbol $\III^*$ or $\II^*$.
\item 
There is no genus-one fibration $\varphi:Y\ra\PP^1$ having a multiple fiber   with Kodaira symbol
$\IV^*$, $\III^*$ or $\II^*$ that admits a $(-2)$-curve   as two-section.
\item
For every genus-one fibration $\varphi:Y\ra\PP^1$ there is another genus-one fibration $\psi:Y\ra\PP^1$
with intersection number $\varphi^\inv(\infty)\cdot \psi^\inv(\infty)=4$.
\end{enumerate}
\end{proposition}

\proof
The first two assertions  are \cite{Ekedahl; Shepherd-Barron 2004}, Theorem A.
The third result is \cite{Cossec; Dolgachev 1989}, Theorem 3.4.1.
\qed

\section{Families of Enriques surfaces}
\mylabel{Families}

In this section we examine    arithmetic properties of Enriques surfaces $Y$
over non-closed ground fields. In fact, we treat them as a special case of the following general notion:
A \emph{family of Enriques surfaces} over an arbitrary ground ring $R$ is an algebraic space $\foY$,
together with a morphism  $f:\foY\ra\Spec(R)$
that is flat,  proper and of finite presentation such that for each geometric point $\Spec(\Omega)\ra\Spec(R)$,
the base-change $Y$ is an Enriques surfaces over $k=\Omega$.  
Let $\shM_\Enr$ be the resulting stack in groupoids 
whose fiber categories $\shM_\Enr(R)$ consists 
of the families of Enriques surfaces $f:\foY\ra\Spec(R)$.  This stack lies over the site $(\Aff/\ZZ)$ of
all affine schemes, endowed with the fppf topology. 
We already   can formulate the main result of this paper:

\begin{theorem}
\mylabel{no enriques over integers}
The fiber category $\shM_\Enr(\ZZ)$ is empty. In other words, 
there is no family of Enriques surfaces over the ring $R=\ZZ$.
\end{theorem}

The proof is long and indirect, and will  finally be achieved  in Section \ref{Proof}. 
The following result  (\cite{Schroeer 2021b}, Theorem 7.2) tells us that we do not have to worry
about 
exceptional Enriques surface:

\begin{theorem}
\mylabel{no exceptional}
There is no family of Enriques surfaces over the ring $R=\ZZ/4\ZZ$ whose
geometric fiber is exceptional.
\end{theorem}

We now collect further preliminary facts.
Throughout, $f:\foY\ra \Spec(R)$ denotes a family of Enriques surfaces over an arbitrary ground ring $R$.
To simplify notation, we set $S=\Spec(R)$.
According to \cite{Artin 1969}, Theorem 7.3 
the fppf sheafification of the naive Picard functor $A\mapsto\Pic(\foY\otimes_RA)$ is representable by a 
group object in the category of algebraic spaces over $R$.  
The following observation, which immediately follows from \cite{Ekedahl; Hyland; Shepherd-Barron 2012}, Corollary 4.3,
will be important:

\begin{proposition}
The algebraic space $\Pic_{\foY/R}$ is a scheme.
The sheaf of subgroups $\Pic^\tau_{\foY/R}$ is representable by an open embedding,
and its structure morphism  is  locally free of rank  two. 
Moreover, the     quotient sheaf $\Num_{\foY/R}$ is representable by a local system of free abelian groups of rank $\rho=10$.
\end{proposition}
 
Each triple $(L,a,b)$ where $L$ is an invertible $R$-module and $a\in L$, $b\in L^{\otimes -1}$ are
elements with $a\otimes b=2_R$ yields a finite flat group scheme $G^L_{a,b}$, characterized by
$$
G^L_{a,b}(A)=\{f\in L\otimes_RA\mid f^{\otimes 2}=a\otimes f \},
$$
with group law $f_1\star f_2=f_1+f_2+ b\otimes  f_1\otimes   f_2$.  
According to \cite{Tate; Oort 1970}, Theorem 2 each 
group scheme that is locally free of rank two is isomorphic to some $G^L_{a,b}$. Cartier duality corresponds
to the involution $(L,a,b)\mapsto (L^{\otimes -1},b,a)$.
The special case $L=R$, $a=1$, $b=2$ yields the constant group scheme $(\ZZ/2\ZZ)_R$,
whereas $a=2$, $b=1$ defines the multiplicative group scheme $\mu_{2,R}$.

Write $P=\Pic^\tau_{\foY/R}$,   and   $G=\uHom(P,\GG_{m,R})$ for the Cartier dual.
According to \cite{Raynaud 1970}, Proposition 6.2.1, the inclusion $P\subset\Pic_{\foY/R}$ corresponds to a global section
of $R^1f_*(G_\foY)$, which we denote by a formal symbol $[\foX]$. This  notation may  be explained as follows:
 The Leray--Serre spectral sequence for the structure morphism $f:\foY\ra S=\Spec(R)$
yields an exact sequence
\begin{equation}
\label{5-term sequence}
0\ra H^1(S, G)\ra H^1(\foY,G_\foY)\ra H^0(S,R^1f_*(G_\foY))\stackrel{d}{\ra} H^2(S,G)\ra H^2(\foY,G_\foY).
\end{equation}
If $[\foX]\in H^0(S,R^1f_*(G_\foY))$ maps to zero under the differential $d$, the section comes
from a $G_\foY$-torsor $\epsilon:\foX\ra \foY$, which we call a \emph{family of canonical coverings}. 

\begin{proposition}
\mylabel{canonical covering}
A family of canonical coverings $\epsilon:\foX\ra\foY$ exists 
under any of the following three conditions:
\begin{enumerate}
\item The pullback map $H^2(S,G)\ra H^2(\foY,G_\foY)$ is injective.
\item The   structure morphism $f:\foY\ra S$ admits a section.
\item The group scheme $P\ra S$ is \'etale and $\Pic(S)=0$.
\end{enumerate}
If a family of canonical coverings exists,  it is unique up to isomorphism and twisting by  pullbacks of $G$-torsors over $S$.
\end{proposition}

\proof
The first assertion and the statement about uniqueness immediately follow from
the five-term exact sequence \eqref{5-term sequence}.
If the structure morphism admits a section, the map $H^2(S,G)\ra H^2(\foY,G_\foY)$ has a retraction and is thus injective,
which gives (ii).

Finally, suppose that $P\ra S$ is \'etale. Suppose first that $R$ is noetherian.
Consider the invertible sheaf $\shL=\Omega^2_{\foY/R}$.
For each $a\in S$, the restriction $\shL|Y$ to the fiber $Y=f^{-1}(a)$ is the dualizing sheaf, which has order two
in $\Pic(Y)$. Moreover, we have $h^i(\shL^{\otimes 2}|Y)=h^i(\O_Y)=0$ for all $i\geq 1$. In turn, the formation of the direct image
$\shN=f_*(\shL^{\otimes 2})$ commutes with base-change, and we conclude that $\shN$ is invertible.
With Nakayama, we infer that $f^*(\shN)\ra \shL^{\otimes 2}$ is bijective.
Using $\Pic(R)=0$ we obtain a trivialization $\varphi:\O_\foY\ra\shL^{\otimes 2}$.
We endow the locally free sheaf
$\shA=\O_\foY\oplus\shL^{\otimes -1}$
with an  algebra structure, by declaring $(a,b)\cdot (a',b') = (aa'+\varphi^{-1}(bb'), ab'+a'b)$.
Then the relative spectrum $\foX=\Spec(\shA)$ is the desired canonical covering.

It remains to treat general rings $R$. Since   $f:\foX\ra\Spec(R)$ is of finite presentation, it comes
from a family of Enriques surfaces over some noetherian subring $R'$. We then can argue as in the previous paragraph,
after enlarging the noetherian ring $R'$ to make $\shN$ trivial. 
\qed

\medskip
We say that our family of Enriques surfaces has \emph{constant local system} if 
the local system $\Num_{\foY/R} $ is isomorphic  to  
$(\ZZ^{\oplus 10})_R$. Likewise, we say that it has \emph{constant Picard scheme} if 
the group space $\Pic_{\foY/R}$ is isomorphic to $(\ZZ/2\ZZ\times\ZZ^{\oplus 10})_R$.
Moreover, we say that it has \emph{split Picard scheme} if the short exact sequence of group spaces
\begin{equation}
0\lra \Pic^\tau_{\foY/R}\lra \Pic_{\foY/R}\lra \Num_{\foY/R}\lra 0
\end{equation}
splits. Our interest in these notions stems from the following fact: 

\begin{proposition}
\mylabel{constant picard scheme}
If the base ring is $R=\ZZ$, then our family of Enriques surfaces $f:\foY\ra\Spec(R)$ has constant Picard scheme.
\end{proposition}

\proof
Choose  a geometric point $\Spec(\Omega)\ra S$.
The local system $\Num_{\foY/R}$ comes from a representation of the algebraic fundamental group
$\pi_1(S,\Omega)$ on the stalk $\Num_{Y/k}(\Omega)=\ZZ^{\oplus 10}$. According to Minkowski's theorem  
(see for example \cite{Neukirch 1999}, Theorem 2.18), the scheme $S=\Spec(\ZZ)$
is simply-connected, so our family of Enriques surfaces has constant local system.

For the ring $R=\ZZ$, up to isomorphism the possible triples describing the group scheme $P=\Pic^\tau_{\foY/R}$
are either $(R,1,2)$ or $(R,2,1)$. In other words, we either have $P=(\ZZ/2\ZZ)_R$ or $P=\mu_{2,R}$.
Seeking a contradiction, we assume the latter holds.
The Picard group $\Pic(R)$ vanishes, because the ring  $R=\ZZ$ is a principal ideal domain.
Furthermore, we have $\Br(R)=0$. Indeed, the Hasse principle (see for example \cite{Neukirch; Schmidt; Wingberg 2008}, Theorem 8.1.17)
gives a short exact sequence
$$
0\lra\Br(\QQ)\lra \bigoplus\Br(\QQ_\primid)\stackrel{\delta}{\lra} \QQ/\ZZ\lra 0,
$$
where the sum runs over all places $\primid$ including the Archimedean one. Since we have $\Br(\ZZ_p)=\Br(\FF_p)=0$ for all primes $p>0$ by Wedderburn,
one sees that $\Br(\ZZ)$ is contained in the Brauer group $\Br(\RR)\cap\Kernel(\delta)=0$.
 
Over the localization $R'=R[1/2]$, the group scheme $\mu_{2,R}$ becomes \'etale, hence there is a canonical covering
$\foX\otimes R'\ra\foY\otimes R'$, by Proposition \ref{canonical covering}.
Now consider the generic fiber $X=\foX\otimes\QQ$, which is a K3 surface. In particular, the Hodge number
$h^{2,0}(X)=h^0(\Omega_X^2)$ is non-zero.
By construction, $X$ has good reduction over the localization $\ZZ_{(p)}$ for all primes $p\geq 3$.
We now consider the ring of Witt vectors 
$A=W(\FF_2^\sep)$, which is the maximal unramified extension of the local ring $\ZZ_2$.

Recall that Fontaine showed that for each proper smooth scheme over $\QQ$ that has good reduction over
$W(\FF_p^\sep)$ for all primes $p>0$, the Hodge numbers $h^{i,j}$ must vanish for $i\neq j$ and $i+j\leq 3$
(\cite{Fontaine 1993}, Theorem 1 on page 44, and first remark afterwards).
Essentially the same  result was independently obtained by Abrashkin, who showed that for each smooth proper
scheme over $\ZZ$, the Hodge numbers $h^{i,j}$ of  the complex fiber vanish for $i\neq j$ and $i+j\leq 3$
(\cite{Abrashkin 1990}, \S7,  Section 6, Theorem on page 516. Note that the formulation given there contains an obvious misprint,
confer the Theorem on page 514).

Since our K3 surface $X=\foX\otimes\QQ$ has Hodge number $h^{2,0}=1$ and   good reduction over   $W(\FF_p^\sep)$ for all $p\neq 2$, it follows that 
the base-change $X_F$ to the field of fractions $F=\Frac(A)$ 
must have bad reduction over $A$. 

One the other hand, using Hensel's Lemma, any $\FF_2^\alg$-valued point for $\foY\ra S$ can be extended to an $A$-valued point.
Again with Proposition \ref{canonical covering} we conclude that a family of canonical covering  $\foX\otimes A\ra\foY\otimes A$ exists,
whence $X_F$ has good reduction over $A$, contradiction.

Summing up, the Picard scheme sits in  a short exact sequence
\begin{equation}
\label{extension local systems}
0\lra (\ZZ/2\ZZ)_R\lra \Pic_{\foY/R}\stackrel{\pr}{\lra} (\ZZ^{\oplus 10})_R\lra 0.
\end{equation}
For each global section $l$ on the right, the preimage $T=\pr^{-1}(l)$ is a representable torsor for
the group scheme $H=(\ZZ/2\ZZ)_R$. In particular, the structure morphism $T\ra \Spec(R)$ is finite
\'etale, of degree two. Since $\pi_1(S,\Omega)=0$, the total space $T$ is the disjoint union of two copies
of $S=\Spec(\ZZ)$. From this we infer that the short exact sequence \eqref{extension local systems} splits.
In turn, our family of Enriques surfaces has constant Picard scheme.
\qed

\medskip
When combined with suitable assumptions
on the Picard group or the Picard group of the base ring $R$, our constancy condition
on   the Picard scheme have remarkable consequence for the fibers of the
morphism $f:\foY\ra\Spec(R)$.
In what follows, let $R\ra k$ be a homomorphism
to a field, choose   an algebraic closure $\bk=k^\alg$,
and write  $Y=\foY\otimes_Rk$  and $\bY=\foY\otimes_R \bk$ for the resulting Enriques surfaces.

\begin{theorem}
\mylabel{family with constant pic}
Suppose that the family of Enriques surfaces  $f:\foY\ra\Spec(R)$ has constant Picard scheme.
Assume furthermore  that  $\Pic(R)$  and  $\Br(R)$ vanish.
Then there is at least one canonical covering $\epsilon:\foX\ra \foY$, and 
the set of isomorphism classes  of such coverings  
is a principal homogeneous space for the  group
$R^\times/R^{\times 2}$. Furthermore, the induced Enriques surfaces $Y$ and $\bY$ enjoy the following properties:
\begin{enumerate}
\item 
The dualizing sheaf $\omega_Y$   has order two in the Picard group.
\item 
Every class $l\in \Num_{Y/k}(\bk)$ comes from an invertible sheaf $\shL$ on $Y$,
which is unique up to twisting by $\omega_Y$.
\item 
Every $(-2)$-curve $\bE\subset\bY$ is the base-change of a $(-2)$-curve $E\subset Y$, which is isomorphic  
to the projective line $\PP^1_k$.
\item 
Every genus-one fibration $ \bY\ra\PP^1_{\bk}$ is the base-change of a genus-one fibration
$ Y\ra\PP^1_k$. The latter has exactly two multiple fibers, both of which are tame and lie over   $k$-rational points.
\item 
There is at least one genus-one fibration $\varphi:Y\ra\PP^1$.
\end{enumerate}
\end{theorem}

\proof
By assumption we have $\Pic^\tau_{\foY/R}=(\ZZ/2\ZZ)_R$. 
Write $S=\Spec(R)$. The Kummer sequence $0\ra \mu_2\ra\GG_m\stackrel{2}\ra\GG_m\ra 0$ gives an exact sequence
$$
\Pic(S)\stackrel{2}{\lra}\Pic(S)\lra H^2(S,\mu_{2})\lra \Br(S)\stackrel{2}{\lra}\Br(S).
$$
Our assumptions ensure that the term in the middle vanishes.
According to Proposition \ref{canonical covering},   there is at least one canonical covering $\epsilon:\foX\ra\foY$.
Moreover, the set of isomorphism classes is a principal homogeneous space for the group $H^1(S,\mu_2)$.
Again by the Kummer sequence we obtain an identification $H^1(S,\mu_2)=R^\times/R^{\times 2}$.
This establishes the statements on the family of canonical coverings.

From $\Pic^\tau_{\foY/R}=(\ZZ/2\ZZ)_R$ we immediately get (i).
To see (ii), we may regard the element $l$ as a section of $R^1f_*(\GG_{m,\foY})$.
The Leray--Serre spectral sequence for the structure morphism $f:\foY\ra S$ yields
$$
H^1(\foY,\GG_{m,\foY})\lra H^0(S, R^1f_*(\GG_{m,\foY}))\lra H^2(S,\GG_{m,S}).
$$
The arrow on the right is the zero map  by assumption, hence the section $l$ comes from an invertible sheaf on $\foY$.
Base-changing along $R\ra k$, we find the desired invertible sheaf $\shL$ on $Y$.


We now check (iii).  Let $\bE\subset\bY$ be a $(-2)$-curve, and consider the invertible sheaf $\bar{\shL}=\O_{\bY}(\bE)$.
We just saw that it is  the base-change of some invertible sheaf $\shL$ on $Y$.
Using $h^0(\bar{\shL})=1$ we infer that there is a unique curve $E\subset Y$ inducing $\bE\subset\bY$, and therefore $E^2=-2$.
This curve is  a twisted form of the projective line. Since its class in $\Num(Y)$ must be primitive in light of the self-intersection number,
and  the intersection form is unimodular,   there
is an invertible sheaf $\shN$ with $(\shN\cdot E)=1$, hence $E\simeq \PP^1$.

It remains to verify (iv). 
Let $\bY\ra\PP^1_{\bar{k}}$ be a genus-one fibration.
According to Theorem \ref{projective contractions}, it is the base-change of a fibration $Y\ra B$, where $B$ is a twisted form
of the projective line.  Consider the invertible sheaf  $\bar{\shL}=\O_\bY(\bar{F})$, where $\bar{F}\subset\bY$
is a half-fiber. It arises from some invertible sheaf $\shL$ on $Y$.   Using $h^0(\bar{\shL})=1$ we infer that there is a unique curve
$F\subset Y$ inducing $\bar{F}\subset \bY$. This curve is a half-fiber for $Y\ra B$. 
If $\bar{F}$ is reducible, it contains a $(-2)$-curve. By the previous paragraph, this gives an inclusion $\PP^1\subset F$.
In particular, $Y$ contains a rational point, hence $B=\PP^1$. Now suppose that $\bar{F}$ is irreducible.
Then $F$ is an integral curve with $h^0(\O_F)=h^1(\O_F)=1$. Since $\Num(Y)$ is unimodular, there
is an invertible sheaf $\shN$ on $Y$ with $(\shN\cdot F)=1$. With Riemann--Roch we conclude that there
is an effective Cartier divisor of degree one on $F$. Again there is a rational point, and $B=\PP^1$.
Summing up, we have shown that each half-fiber on $\bY$ induces a half-fiber on $Y$ that lies over a rational point.
Since there are exactly two half-fibers on $\bY$, both of which are tame, the same holds for $Y$. 
\qed

\section{Geometrically rational elliptic surfaces}
\mylabel{Geometrically rational}

In this section we discuss the relation between Enriques surfaces and geometrically rational surfaces.
The results are well-known, but I could not find suitable references in the required generality.
Let $k$ be an arbitrary ground field of arbitrary characteristic $p\geq 0$.
A proper surface $J$ is called a \emph{geometrically rational surface} if it is   geometrically integral and the
base-change to $k^\alg$ is birational to the projective plane. Clearly we have $h^0(\O_J)=1$.

Suppose $J$ is such a surface that is  endowed with a   genus-one fibration $\phi:J\ra \PP^1$.
We assume that the fibration is jacobian and relatively minimal, and that the total space is smooth.
Fix a section $E\subset J$, and let $J\ra Z$ be the contraction of all vertical curves disjoint from the zero-section $E\subset J$.
The normal proper surface $Z$ is called the \emph{Weierstra\ss{} model}.
Write $\PP^1=\Spec k[t]\cup \Spec k[t^{-1}]$ for some indeterminate $t$.

\begin{proposition}
\mylabel{global weierstrass equations}
There are polynomials $a_i\in k[t]$ of degree $\deg(a_i)\leq i$
such that 
\begin{equation}
\label{two equations}
\begin{gathered}
y^2+a_1xy+a_3y=x^3+a_2x^2+a_4x+a_6,\\
y^2+a'_1xy+a'_3y=x^3+a'_2x^2+a'_4x+a'_6\qquad \text{with $a_i'=a_i/t^i$}
\end{gathered}
\end{equation}
are  Weierstra\ss{} equations for  $Z$ over the affine open sets $\PP^1\smallsetminus\{\infty\}=\Spec k[t]$
and $\PP^1\smallsetminus \{0\}=\Spec k[t^\inv]$, respectively.  
\end{proposition}

\proof
To simplify notation we write $\O_X(n)$ for the pullbacks of the invertible sheaves $\O_{\PP^1}(n)$.
The Canonical Bundle Formula \cite{Bombieri; Mumford 1977} gives $\omega_J=\O_X(d)$ for some integer $d$.
Base-changing to the algebraic closure $k^\alg$ we can apply \cite{Cossec; Dolgachev 1989}, Proposition 5.6.1 and conclude $d=-1$.
The Adjunction formula for $E\subset J$ gives the self-intersection number $(E\cdot E)=-1$,
and we have $\omega_{J/\PP^1}=\O_X(1)$.

As explained in \cite{Deligne 1975}, Section 1 the sheaves $\shE_i=f_*\O_X(iE)$ are locally free,
with $\rank(\shE_i)=i$ for  $i\geq 1$, and $\shE_0=\O_{\PP^1}$. The canonical inclusions of $\shE_i$
define an increasing filtration on $\shE=\shE_3$ with $\gr^*(\shE)=\bigoplus_{i=0,2,3}\O_{\PP^1}(-i)$.
The invertible sheaf $\shL=\O_X(3E)$ is relatively very ample, and yields a closed embedding
$Z\subset\PP(\shE)$.
On the projective line, we actually may choose splittings $\shE=\O_{\PP^1} \oplus\O_{\PP^1}(-2)\oplus\O_{\PP^1}(-3)$.

Write $\PP^1=\Proj k[t_0,t_1]$, such that $t=t_0/t_1$, and consider the global section $\omega=t_1$
of the invertible sheaf
$\O_{\PP^1}(1)=\Omega^1_{J/\PP^1}|E$. Over the affine open set $U=\Spec k[t]$, we choose local sections $x=t_0^{-2}$
and $y=t_0^{-3}$ in the second and third summand of $\shE$. 
As explained in loc.\ cit., there are  polynomials $a_i\in k[t]$ such that 
the first equation in \eqref{two equations} holds in the sheaf of graded rings $\Sym^\bullet(\shE)$ over $U$.
The situation over $U'=\Spec k[t^{-1}]$ is similar: Here we choose $\omega'=t_0$ and $x'=t_1^{-2} $ and $y'=t_1^{-3}$.
This gives polynomials $a_i'\in k[t^{-1}]$ with 
$y'^2+a_1'x'y'+\ldots = x'^3+\ldots$. Now we use that the  coefficients are unique, once the local sections 
in the summands of $\shE$ are chosen.
Multiplying the equation over $U$ with $t^{-6}=(t_0/t_1)^6$ and comparing coefficients with the equation over $U'$
we get $a_i'=a_i/t^i$. This gives the second equation in \eqref{two equations}, and  ensures $\deg(a_i)\leq i$.
\qed

\medskip
Note that at least one of the  coefficients $a_i$ in  \eqref{two equations} is non-constant, and at least one is not divisible by $t^i$,
because otherwise the geometrically rational smooth surface $J$ would be isomorphic to the product of $\PP^1$ and some elliptic curve.
Furthermore, this property of the coefficients remains  true after any degree-preserving change of coordinates.

Now let $Y$ be an Enriques surface, and suppose there is a genus-one fibration $\varphi:Y\ra\PP^1$.
Let $\eta\in\PP^1$ be the generic point, with function field $K=k(t)$. 
If the fibration is elliptic, $J_\eta=\Pic^0_{Y_\eta/K}$ is an elliptic curve.
If the fibration is quasielliptic, $\Pic^0_{Y_\eta/K}$ is a twisted form of the additive group $\GG_{a,K}$,
and we write $J_\eta$ for its unique regular compactification.
In both cases, $J_\eta$ is a regular curve with $h^0(\O_{J_\eta})=h^1(\O_{J_\eta})=1$, and we write
$\phi:J\ra\PP^1$ for the resulting relatively minimal genus-one fibration.
It comes with a zero-section,  and is thus a jacobian fibration.

\begin{proposition}
\mylabel{jacobian geometrically rational}
The regular surface $J$ is geometrically rational.
\end{proposition}

\proof
Let $U\subset Y$ be the open set where the coherent sheaf $\Omega^1_{Y/\PP^1}$ is invertible.
It is dense, because the generic fiber $Y_\eta$ is geometrically reduced, hence the image
$V=\varphi(U)$ is a dense open set. 
For each point $b\in V$, the fibration $\varphi:Y\ra\PP^1$ acquires a section after base-changing
to the strict henselization  $\O_{\PP^1,b}\subset R$, and we get $Y\times_{\PP^1}\Spec(R)\simeq J\times_{\PP^1}\Spec(R)$.
It follows that the scheme $J$ is smooth on the open set $\psi^{-1}(V)$, and we infer that the base-change $\bJ$ remains integral.

Moreover, the jacobian for the base-change $\bY$ coincides with the base-change of the jacobian $\bJ$, at least over some
dense open set in $\PP^1_{k^\alg}$, so we may assume from the start that $k$ is algebraically closed.
The assertion now follows from \cite{Cossec; Dolgachev 1989}, Theorem 5.7.2.
\qed

\medskip
In turn, the genus-one fibration $\varphi:Y\ra \PP^1$ yields polynomials $a_i\in k[t]$ of degree $\deg(a_i)\leq i$
such that the two Weierstra\ss{} equations in \eqref{two equations} describe the Weierstra\ss{} model $Z\ra \PP^1$
of the jacobian fibration $\phi:Y\ra\PP^1$. We shall see that for $k=\FF_2$, the possibilities for these polynomials
impose
strong restrictions on the Kodaira symbols occurring for $\varphi$.

\section{Counting points over finite fields}
\mylabel{Counting points}

In this section we count the number of rational points  on Enriques surfaces  with constant local system over finite fields.
The observation has nothing in particular to do with Enriques surfaces, so we first work in a general setting.
Let $k=\FF_q$ be a finite ground field, for some prime power $q=p^\nu$, and   $X$ be a smooth proper scheme of dimension $n\geq 0$,
with $h^0(\O_X)=1$.

Set $N_s=\Card X(\FF_{q^s})$ for  the integers $s\geq 1$, such that $N_1$ is the number of rational points.
The \emph{Hasse--Weil zeta function} is 
\begin{equation}
\label{zeta function}
Z(X,t) = \prod_a \frac{1}{1-t^{\deg\kappa(a)}} =\exp\left(\sum_{s=1}^\infty \frac{N_s }{s}t^s\right),
\end{equation}
where the product runs over all closed points $a\in X$. This formal power series actually belongs to
the subring $\QQ(t)\subset\QQ[[t]]$, by Dwork's Rationality Theorem \cite{Dwork 1960}.

Choose some algebraic closure $\bk=k^\alg$, write $\bX=X\otimes \bk$, and fix a prime $\ell\neq p$. 
For each $i\geq 0$ and each integer $j$ we can form the  $\ell$-adic cohomology groups
\begin{equation}
\label{l-adic cohomology groups}
H^i(\bX,\QQ_\ell(j))= \left( \invlim_r H^i(\bX,\mu^{\otimes j}_{\ell^r})\right)\otimes_{\ZZ_\ell}\QQ_\ell.
\end{equation}
Here $j\in\ZZ$ refers to the \emph{Tate twist} in the coefficient sheaves. Note that $H^i(\bX,\QQ_\ell(j))$
arises from $H^i(\bX,\QQ_\ell)$ by tensoring with the one-dimensional  vector space
$\QQ_\ell(j)=(\invlim_r\mu_{l^r}(\bk))\otimes_{\ZZ_\ell}\QQ_\ell$, which comprises only information on the action of $\Gal(\bk/k)$
on the $l$-primary roots of unity in $\bk$.
In particular, the $\QQ_\ell$-dimensions $b_i\geq 0$ of the $\ell$-adic cohomology 
$H^i(\bX,\QQ_\ell(j))$ do not depend on the Tate twist, and are called  the \emph{Betti numbers} of $X$.

The   power $F^\nu_X:X\ra X$ of the absolute Frobenius is a $k$-morphism, and induces by base-change a 
$\bk$-morphism $\Phi=F^\nu_X\otimes \id_\bk$ on $\bX$. Let us write $f_{i,j}=\Phi^*$ for the induced 
$\QQ_\ell$-linear endomorphism on the    cohomology group $H^i(\bX,\QQ_\ell(j))$. Of particular interest are the $f_i=f_{i,0}$:
The zeta function takes the form 
$$
Z(X,t)=\prod \det(1-tf_i)^{(-1)^{i+1}},
$$
by the  Grothendieck--Lefschetz Trace Formula \cite{FL}.
Here the  polynomial factors can be seen as
``reciprocal'' characteristic polynomials 
$\det(1-f_it)=t^{b_i}\det(t^{-1}-f_i)$ of the endomorphisms $f_i$. 
The characteristic polynomials of $f_{i,j}$ have coefficients from $\QQ$, and the eigenvalues
$\alpha\in\QQ^\alg$ are algebraic integers, 
of  length $|\alpha|=q^{i/2-j}$ in all complex embeddings, according to the Riemann Hypothesis (\cite{Deligne 1980}, Corollary 3.3.9,
see   \cite{Deligne 1974}, Theorem 1.6 for the projective case).

The following observations explain the  effect of the Tate twist:
First of all, we have a canonical identification $H^i(\bX,\QQ_\ell(j))=H^i(\bX,\QQ_\ell)\otimes \QQ_\ell(j)$.
Moreover, the  composition $(F_X^\nu\otimes\id_\bk)\circ (\id_X\otimes F_\bk^\nu)$ is a power of the absolute Frobenius of $\bX$, which
acts trivially on \'etale cohomology   with respect to any coefficient sheaf.
Finally, the Frobenius power $F_\bk^\nu$ acts by multiplication with $q=p^\nu$ on $\QQ_\ell(1)$. Combining these three observations,
we see that $f_{i,j}=f_i\otimes q^{-j}$, which reveals that the eigenvalues of $f_i$ are obtained from the eigenvalues of $f_{i,j}$
by multiplication with $q^j$.

The Tate twists become   relevant in connection with algebraic cycles:
Let $\CH^j(\bX)$ be the \emph{Chow group} of $j$-codimensional cycles modulo linear equivalence, which comes with 
a \emph{cycle class map}  $\CH^j(\bX)\ra H^{2j}(\bX,\QQ_\ell(j))$. 
We likewise can form Chow groups and  $\ell$-adic cohomology on the scheme $X$, before base-changing to the algebraic closure,
and have a commutative diagram 
\begin{equation}
\label{cycle class maps}
\begin{CD}
\CH^j(X)	@>>> 	H^{2j}(X,\QQ_\ell(j))\\
@VVV			@VVV\\
\CH^j(\bX)	@>>> 	H^{2j}(\bX,\QQ_\ell(j))
\end{CD}
\end{equation}
of cycle class maps. For more details on cycle class maps see \cite{Deligne 1977}.

\begin{proposition}
\mylabel{generalization wedderburn}
Suppose for all integers $j\geq 0$ we have $H^{2j+1}(\bX,\QQ_\ell)=0$, and that the 
composite map $\CH^j(X)\otimes\QQ_\ell \ra H^{2j}(\bX,\QQ_\ell(j))$ is surjective.
Then the zeta function is
\begin{equation}
\label{simple zeta function}
Z(X,t) = \prod_{j=0}^{n}\left(\frac{1}{1-q^jt}\right)^{b_{2j}}.
\end{equation}
In particular, the number of rational points is $\Card X(\FF_q)= \sum_{j=0}^{n} b_{2j}q^j$, and  the set of rational points is non-empty.
\end{proposition}

\proof
Suppose  that the endomorphism $f_{2j}$ is the homothety with eigenvalue $q^j$. Then 
$$
t^{b_{2j}}\det(t^{-1}-f_{2j}) = t^{b_{2j}}(t^{-1}-q^j)^{b_{2j}} = (1- q^jt)^{b_{2j}},
$$
and the assertion on the zeta function follows. Expanding \eqref{zeta function} and \eqref{simple zeta function} 
as formal power series in $t$ and comparing linear terms, we get the statement on the number of rational points.
Since $b_0=1$ the set  $X(\FF_q)$ must be non-empty.

It remains to verify the assertion on the endomorphism $f_{2j}$. 
Set $r=b_{2j}$ and choose prime cycles $Z_1,\ldots,Z_r\subset X$ whose cycle classes in $H^{2j}(\bX,\QQ_\ell(j))$ form
a vector space basis. Obviously, the action of $\id_X\otimes F_\bk$ on the base changes $\bar{Z}_1,\ldots,\bar{Z}_r\subset\bX$ is trivial,
so  the same holds for the action on $H^{2j}(\bX,\QQ_\ell(j))$. The latter action is   inverse to the action of
$\Phi=F_X^\nu\otimes\id_\bk$, as observed above.  In other words, $f_{2j,j}$ is the identity map.
We also observed above that $f_{2j,j}=f_{2j}\otimes q^{-j}$, thus $f_{2j}$ is the homothety with eigenvalue $q^j$.
\qed

\medskip
The conditions obviously hold if $b_{2j+1}=0$ and $b_{2j}=1$ for all $j\geq 0$. One may see the above result
as a generalization of  Wedderburn's Theorem, which states that every Brauer--Severi variety $X$ contains a rational point, 
and is thus isomorphic to   projective  space $\PP^n$. 
Let us now specialize to surfaces. Recall that $\rho=\rank \Num(\bX)$ denotes the Picard number.

\begin{corollary}
Suppose that $X$ is a surface with $b_1=0$ and $b_2=\rho$.
If the local system $\Num_{X/k}$ is constant, we have $\Card X(\FF_p)= 1+b_2q+q^2$. 
\end{corollary}

\proof
We have $b_0=1$, and with Poincar\'e Duality get $b_4=1$ and $b_1=0$. The cycle class map $\CH^1(X)\ra H^2(\bX,\QQ_\ell(1))$
factors over $\Num(X)$, and  by the assumptions  the inclusion $\Num(X)\otimes\QQ_\ell\subset H^2(\bX,\QQ_\ell(1))$ is an equality.
Thus we may apply the  Proposition and  the formula on  $\Card X(\FF_p)$ follows.
\qed

\begin{corollary}
\mylabel{points on enriques}
Suppose that $X$ is either an Enriques surface, or a  geometrically rational   surface endowed
with a relatively minimal genus-one fibration $X\ra\PP^1$. If the local system $\Num_{X/k}$ is constant,
we have  $\Card X(\FF_q) = 1+10q+q^2$.
\end{corollary}

\proof
In both cases, we have $b_1=0$. 
In light of the previous corollary, we only have to show
that $\bX$ has Betti number $b_2=10$. For Enriques surfaces, this holds by definition.
In the other case, the induced  fibration over $\bk$ remains a relatively minimal genus one-fibration.
We have $K_{\bX}^2=0$ by \cite{Cossec; Dolgachev 1989}, Proposition 5.6.1.
Hence there is a morphism $\bX\ra\PP^2_\bk$ that can be factored as a sequence of nine blowing ups,
so we also have $b_2=10$.   
\qed

\medskip
We shall be particularly concerned with the case $q=p=2$. The above formula then gives $\Card X(\FF_2)=25$.

\section{Multiple fibers and isogenies}
\mylabel{Multiple fibers}

In this section we relate  multiple   fibers in elliptic fibrations with simple elliptic fibers
via isogenies. The general  set-up is as follows: Suppose $R$ is an   discrete valuation ring,
with field of fractions $F=\Frac(R)$ and residue field $k=R/\maxid_R$. For the sake of exposition,
we assume that  the local ring $R$ is excellent, which here means that the field extension $\Frac(R)\subset\Frac(\hat{R})$
is separable, compare the discussion in \cite{Schroeer 2020}, Section 4.
Let $E_F$ be an elliptic curve over $F$ and $X_F$ be a principal homogeneous space,
representing a cohomology class in $H^1(F,E_F)$ with respect to the   flat topology.
Suppose that $E_F$ has good reduction over $R$, that is,   arises as generic fiber for some
family of elliptic curves $E\ra\Spec(R)$.
Furthermore, we assume that the principal homogeneous space is the generic fiber of some relatively minimal regular model
$X\ra\Spec(R)$, and that the underlying reduced subscheme $X_{k,\red}$  of the closed
fiber $X_k$ is an elliptic curve. We then observe:

\begin{proposition}
\mylabel{reduction isogeneous}
Under the above assumptions, the two elliptic curves $E_k$ and $X_{k,\red}$ over the residue field $k$ are isogeneous.
\end{proposition}

\proof
Without restriction, we may assume that the discrete valuation ring $R$ is henselian.
Let $m\geq 1$ be the multiplicity of the closed fiber for $X\ra\Spec(R)$, such
that $X_k=mX_{k,\red}$. By assumption, the reduction $X_{k,\red}$ is an elliptic curve,  in particular contains a
$k$-rational point $d\in X_{k,\red}$. This is an effective Cartier divisor,
and can be extended to an effective Cartier divisor $D_k\subset X_k$ on the   fiber.
In light of \cite{EGA IVd}, Theorem 18.5.11 we may extend further, and regard it as the closed
fiber of some effective Cartier divisor $D\subset X$.
The scheme $D$ is regular, because the intersection of Cartier divisors $D\cap X_{k,\red}=\Spec \kappa(d)$
is  zero-dimensional and regular. Writing $D=\Spec(R')$, we get a finite extension of discrete valuation rings
$R\subset R'$. The residue field extension has degree one, hence $k'=R'/\maxid_{R'}$ is a copy of $k=R/\maxid_R$.
By construction, the base-change $X_{R'}$
acquires a section;  by abuse of notation, we use the same symbol for the
effective Cartier divisor $D\subset X$ and the resulting section $D\subset X_{R'}$.
Note, however, that the scheme $X_{R'}$   usually becomes non-normal, and contains $D$
as a Weil divisor rather than a Cartier divisor.
Let $X'\ra X_{R'}$ be the normalization. By the valuative criterion for properness,
the inclusion $D\subset X_{R'}$ lifts to an inclusion $D\subset X'$.

To proceed, consider the base-change $E'=E_{R'}$, which is a family of elliptic curves over $R'$.
The section $D\subset X'$ yields a  rational map $X'\dashrightarrow E'$, which is 
defined on the generic fiber. Choose  proper birational  maps $X'\leftarrow Y\ra E'$
that restrict to identities on the generic fiber,  where $Y$ is a normal scheme.
According to \cite{Raynaud 1970}, Theorem 8.2.1 the morphism $Y\ra\Spec(R')$ is cohomologically flat,
because it admits a section. In turn, the closed fiber $C=Y\otimes_{R'}k'$ has cohomological invariants
$h^0(\O_C)=h^1(\O_{C})=1$.

Let $C_1,\ldots,C_r\subset C$ be those irreducible components
that are strict transforms of irreducible components in $X'_{k'}$, endowed with reduced scheme structure.
Then each $C_i$ dominates the elliptic curve $X_{k,\red}$,
hence is a curve with cohomological invariants $h^0(\O_{C_i})\geq 1$ and  $h^1(\O_{C_i})\geq 1$,
by Lemma \ref{covering} below. We now claim that  $r=1$.
Seeking a contradiction,   suppose that $r\geq 2$. In the   exact sequence
$\O_C\ra \bigoplus_{i=1}^r \O_{C_i}\ra \shF\ra 0$,
the term on the right is a skyscraper sheaf, whereas the other sheaves are one-dimensional.
From this we infer that the map  $H^1(C,\O_C)\ra \bigoplus H^1(C_i,\O_{C_i})$
is surjective. Since each summand is a non-trivial vector space, we must have $r=1$, with $h^1(\O_{C_1})=1$.

Now write $C_1'\subset Y$ for the strict transform of  $E_{k'}$. This is an elliptic curve, because the birational map $C_1'\ra E_{k'}$ 
must be an  isomorphism.
Arguing as in the preceding paragraph, we see $C_1= C_1'$.
Summing up, the morphism $C_1\ra E_{k'}=E_k$ is an isomorphism. Using its inverse, we obtain
the desired isogeny   $E_k=C_1\ra   X_{k,\red} $.
\qed

\medskip 
It is easy to construct an example, even in arbitrary dimensions: Suppose for simplicity that $R$ contains a field of representatives,
that is, a subfield that surjects onto the residue field, for example $R=k[[t]]$.
Let $A$ be an abelian variety over $k$, endowed with a finite subgroup $G\subset A(k)$, 
and suppose that there is a finite Galois extension $F\subset F'$ with
Galois group $\Gal(F'/F)=G$ such that $R\subset R'$ has trivial residue field extension. 
Consider the diagonal action on $A\times\Spec(R')$, which is free,
and the resulting quotient $X=(A\times\Spec(R'))/G$. Then the generic fiber is a twisted form
of $A\otimes_k F$, whereas the closed fiber contains the abelian variety $A/G$ as reduction.
See the dissertation of Zimmermann  for more on this construction (\cite{Zimmermann 2019a}, \cite{Zimmermann 2021}).

In light of general results of Serre, Tate and Honda on isogeny classes of abelian varieties over finite fields,
the  above result has the following consequences for finite fields:

\begin{corollary}
\mylabel{reduction same zeta}
Assumptions as in the proposition. If the residue field $k=\FF_q$ is finite, 
then for each integer $r\geq 1$, the number of $\FF_{q^r}$-valued points
in $E_k$ and $X_k$ coincide.
\end{corollary}

\proof
According to  \cite{Tate 1982}, Theorem 1 any    two isogeneous elliptic curves over 
the finite field $\FF_q$ have the same    zeta function, in other words,
the same number of $\FF_{q^r}$-valued points, $r\geq 1$.
\qed

\begin{corollary}
\mylabel{reduction isomorphic}
Assumptions as in the proposition. If the residue field is $k=\FF_2$, then
the elliptic curve $E_k$ is isomorphic to the scheme $X_{k,\red}$.
\end{corollary}

\proof
We recalled in Proposition \ref{elliptic curves} that  there are five isomorphism classes of elliptic
curves over $k=\FF_2$. In each case, the group of rational points is cyclic, and each of the
groups $\ZZ/n\ZZ$, $1\leq n\leq 5$ does occur. Hence the number of rational points 
already determines the isomorphism class, and the assertion follows from Corollary \ref{reduction same zeta}.
\qed

\medskip
Note that according to the results of  Tate \cite{Tate 1982} and Honda \cite{Honda 1968},
the isogeny classes of simple abelian varieties $A$ over a finite field $k=\FF_q$, $q=p^\nu$
corresponds to the conjugacy classes of $q$-Weil numbers $\pi\in\CC$.

In the proof for Proposition \ref{reduction isogeneous}, we have used the following observation on algebraic curves:

\begin{lemma}
\mylabel{covering}
Let $k$ be a ground field, $E$ be an elliptic curve, $C$ be an integral proper curve, and $f:C\ra E$ be a dominant morphism.
Then  $h^1(\O_C)\geq 1$. Moreover,  if $h^1(\O_C)=1$ then the curve $C$ is regular.
\end{lemma}

\proof
Let $d\geq 1$ be the degree of the dominant morphism $f:C\ra E$.
The induced homomorphism $\Pic(E)\ra\Pic(C)$ comes with a norm map in the reverse direction (see \cite{EGA II} Section 6.5), 
with the property $N_{C/E}(\shL|C)=\shL^{\otimes d}$ for all invertible sheaves $\shL$ on $E$.

To check $h^1(\O_C)\geq 1$, we may assume that $k$ is separably closed.
Seeking a contradiction, we assume that $h^1(\O_C)=0$, such that $\Pic^0_{C/k}$ is trivial.
We conclude that the abelian group $\Pic^0(E)=E(k)$ is annihilated by $d$.
By Bertini's Theorem (for example \cite{EGA V}, Proposition 4.3 or \cite{Jouanolou 1983}, Chapter I, Theorem 6.3), we find some smooth divisor $D\subset E$ that is disjoint
from the scheme $E[d]$ of $d$-torsion. Each point $a\in D$ of this subscheme is rational,
because $k$ is separably closed. The corresponding invertible sheaf $\shL=\O_E(a-0_E)$
 has $\shL^{\otimes d}\neq \O_E$, contradiction.

Now suppose that $h^1(\O_C)=1$. 
It remains to check that $C$ is regular.
Let $C'\ra C$ be the normalization. The short exact sequence $0\ra\O_C\ra\O_{C'}\ra\shF\ra 0$
yields a long exact cohomology sequence
$$
0\ra H^0(C,\O_C)\ra H^0(C',\O_{C'})\ra H^0(C,\shF)\ra H^1(C,\O_C)\ra H^1(C',\O_{C'})\ra 0.
$$
The term on the right does not vanish, by the preceding paragraph, hence the map on the right is bijective.
Thus $H^1(C,\O_C)$ is a one-dimensional  vector space over $k=H^0(C,\O_C)$. It is also a vector space over the field extension $k'=H^0(C',\O_{C'})$,
and we infer $k=k'$.
Then the above exact sequence   gives $h^0(\shF)=0$,   whence the torsion sheaf $\shF$ vanishes.
Thus the one-dimensional scheme $C$ is normal, hence regular.
\qed

\section{Passage to   jacobian   fibrations}
\mylabel{Passage to jacobian}

Throughout this section, $Y$ is an Enriques surface with constant Picard scheme $\Pic_{Y/k}$ over the field $k=\FF_2$.
According to Theorem \ref{family with constant pic}, 
there is at least one  genus-one fibration $\varphi:Y\ra\PP^1$. 
We now assume that this fibration is elliptic,  and let
$\phi:J\ra\PP^1$ be the resulting jacobian fibration, where $J$ is a  geometrically rational surface.
The goal of this section is to establish the following result:

\begin{theorem}
\mylabel{jacobian constant pic}
The geometrically rational surface $J$ has constant Picard scheme.
\end{theorem}

This ensures that the relation between the elliptic fibration on the Enriques surface
and its jacobian is stronger that one might expect. This crucial observation will allow us in the next sections
to use Weierstra\ss{} equations to reduce the possibilities for $J$, and thus   also for $Y$.

Choose an an algebraic closure $k\subset k^\alg$, and write $\bY=Y\otimes k^\alg$ and $\bJ=J\otimes k^\alg$
for the   base-changes.
Throughout, $\ba:\Spec(\Omega)\ra \PP^1$ denotes a  geometric point   whose image
point $a\in \PP^1$ is closed.

\begin{proposition}
\mylabel{same symbols}
For each  such  $\ba:\Spec(\Omega)\ra \PP^1$, the geometric fibers $J_\ba$ and $Y_\ba$ have the same Kodaira symbols.
If  $Y_a$ is simple,  we actually have $J_a\simeq Y_a$.
\end{proposition}

\proof
The first assertion follows from a general result of Liu, Lorenzini and Raynaud 
(\cite{Liu; Lorenzini; Raynaud 2005}, Theorem 6.6). Now suppose that $Y_a$ is simple.
In order to check $J_a\simeq Y_a$ we may base-change to $\kappa(a)=\FF_{2^\nu}$ and assume that $a\in\PP^1$ is rational.
If the fiber contains a rational point $y\in Y_a$ in the regular locus, the resulting 
effective Cartier divisor gives an identification $J\otimes R\simeq Y\otimes R$,
where $R=\O_{\PP^1,a}^h$ is the henselization, and in particular $J_a\simeq Y_a$.
It remains  to verify the existence of such a  rational point. 
For smooth fibers, this follows from Lang's result (\cite{Lang 1956}, Theorem 2).
Suppose now that $C=Y_a$ is singular, with irreducible components $C_1,\ldots, C_r$.
Each $C_i$ is birational to $\PP^1$, according to Proposition \ref{constant num over finite field}.
In case $r=1$  the description before Proposition \ref{canonical type rational points} 
reveals that the desired rational point exists.
If $r\geq 2$, we actually have $C_i\simeq \PP^1$, and
the canonical map $\Gamma(\bC)\ra\Gamma(C)$ is a graph isomorphism respecting edge labels, 
again by Proposition \ref{constant num over finite field}.
Choose a component $C_i$ corresponding to a terminal vertex of $\Gamma(C)$. Then there is exactly one
component $C_j$ intersecting $C_i$, which has $(C_i\cdot C_j)\leq 2$.
In turn, from the three rational points in $C_i=\PP^1$ at least one is contained in
the regular locus of $C_\red$. 
\qed

\medskip
We now state   several facts about the fibers of $\phi:J\ra\PP^1$, which can be formulated without referring to the Enriques surface $Y$:
\newcounter{jac}\setcounter{jac}{0}
\renewcommand{\thejac}{\rm \roman{jac}}
\newcommand{\labeljac}[1]{\refstepcounter{jac}\mylabel{#1}}
\newcommand{\refjac}[1] {{\rm (\ref{#1})}}

\begin{proposition}
\mylabel{jacobian fiber properties}
Let $\ba:\Spec(\Omega)\ra \PP^1$ be a geometric point whose image point $a\in \PP^1$ is closed.
Then the following holds:
\begin{enumerate}
\item 
\labeljac{non-rational J}
If $a\in \PP^1$ is non-rational, then the geometric fiber $J_\ba$ is irreducible.
\item
\labeljac{I0* J}
The Kodaira symbol of  $J_\ba$ is different from  $\I_0^*$.
\item
\labeljac{semistable J}
If $J_\ba$ is semistable, then  the   map $\Gamma(J_\ba)\ra\Gamma(J_a)$ is a graph isomorphism. 
\item
\labeljac{rational J}
There is at most one rational point $b\in\PP^1$ whose fiber $J_b$ is semistable or supersingular. 
\item
\labeljac{mw J}
The canonical inclusion $\MW(J/\PP^1)\subset \MW(\bJ/\bar{\PP}^1)$ 
of Mordell--Weil groups is an equality.
\end{enumerate} 
\end{proposition}

\proof
Recall from Proposition \ref{same symbols} that for all closed points $a\in \PP^1$, the fibers $J_a$ and $Y_a$ have the same Kodaira symbol.
Suppose now that $k\subset\kappa(a)$ is not an equality. Then $Y_\ba$ is   irreducible  
by Theorem \ref{geometric and schematic fiber}, which  settles \refjac{non-rational J}.
From Proposition \ref{canonical type rational points} we get \refjac{I0* J}.
Suppose next  that $J_a$   is semistable. Then $Y_a$ must be simple, by  Lemma \ref{multiple for p=2}, hence
$J_a\simeq Y_a$, again from  Proposition \ref{same symbols}. Assertion \refjac{semistable J} 
thus follows from Theorem \ref{geometric and schematic fiber}.

By Theorem \ref{family with constant pic}  there are two rational points $a\in \PP^1$ where the fiber $Y_a$ is multiple.
The corresponding $\varphi^{-1}_\ind(a)$ are neither semistable nor supersingular, again by Lemma \ref{multiple for p=2}.
So from the three rational points, there is only one $b\in\PP^1$  where $\varphi^{-1}_\ind(b)$ and hence $J_b$
can be semistable or supersingular, which gives \refjac{rational J}.

It remains to verify the last assertion \refjac{mw J}. According to the definition of Mordell--Weil groups we have
$$
\MW(J/\PP^1)=J(F)=\Pic^0_{Y_F/F}(F)=\Pic^0(Y_F),
$$
where $F$ is the function field of the projective line. The latter equality holds because $\Br(F)=0$, by Tsen's Theorem.
Choose a two-section $R\subset Y$, and consider the  surjective homomorphism
$$
\Pic(Y)\lra \Pic^0(Y_F),\quad \shL\longmapsto \shL(-nR) | Y_F,
$$
where $2n=\deg(\shL|Y_F)$. The kernel $T(Y/\PP^1)$ of this map is generated by the irreducible components 
of the closed fibers for $\varphi:Y\ra\PP^1$, together with $R$.  Consider the commutative diagram
$$
\begin{CD}
T(Y/\PP^1)	@>>> 	\Pic(Y)\\
@VVV		@VVV\\
T(\bY/\bar{\PP}^1)	@>>>	\Pic(\bY)
\end{CD}
$$
induced by base-change. By assumption,  the vertical map on the right is bijective.
According to Theorem \ref{geometric and schematic fiber}, the same holds for the vertical map on the left.
In turn, the induced vertical map  between cokernels is bijective, thus \refjac{mw J} holds.
\qed

\medskip
The number of irreducible components   of the geometric fiber $J_\ba$ depends only on the image point $a\in \PP^1$;
we denote this integer by $r_a\geq 1$. It coincides with the number of irreducible components in both $Y_\ba$ and $Y_a$.
If $a\in\PP^1$ is rational, we furthermore set
$$
n_a=
\begin{cases}
i		& \text{if $J_a$ is smooth and isomorphic to $E_i$;}\\
2r_a 	& \text{if $J_a$ is split semistable;}\\
2r_a+2	& \text{if $J_a$ is non-split semistable;}\\
2r_a+1	& \text{if $J_a$ is unstable.}
\end{cases}
$$ 
Recall from Proposition \ref{elliptic curves} that the $E_i$ denote the five elliptic curves over the prime field
$k=\FF_2$,  indexed by the order of the group of  rational points.
Also note  that $n_a=| Y_a(k)|$, and thus
$\sum n_a=25$ by  Corollary \ref{points on enriques}, where the sum runs over the three rational points $a\in\PP^1$. 
The numbers $n_a\geq 1$ are defined,  however,  
in terms of the irreducible components of the schematic and geometric fibers $J_a$ and $J_\ba$, and it is a priori not clear 
how they  are related to $|J_a(k)|$.

We now forget that $J$ comes from an elliptic Enriques surface,
and only demand that it satisfies the properties  observed above.
Theorem \ref{jacobian constant pic} becomes a consequence of the following more general statement:

\begin{theorem}
\mylabel{constant pic}
Let $J$ be a geometrically rational surface, endowed with 
an elliptic fibration $\phi:J\ra\PP^1$ that is relatively minimal and jacobian.
Suppose that conditions {\rm (i)--(v)} of Proposition \ref{jacobian fiber properties} and the equation $\sum n_a=25$ hold.
Then $J$ has constant Picard scheme.
\end{theorem}

The elliptic surface $J$ is geometrically rational, so the Picard scheme $\Pic_{Y/k}$ is  a local system
of free abelian groups. Their  rank is $\rho=10$, by the minimality of the fibration.
We have to show that the Galois group $G=\Gal(k^\alg/k)$ acts trivially on $\Pic_{J/k}(k^\alg)=\Pic(\bJ)=\ZZ^{\oplus 10}$.
Since there is no torsion part, it suffices to find a subgroup of finite index for which this holds.
The group $\Pic(\bJ)$ is generated  by a fiber, the vertical curves $\Theta$ disjoint from the zero-section,
and the sections. The former is obviously defined over $k$, and the latter  by condition \refjac{mw J}.
According to  condition \refjac{non-rational J}, each $\Theta\subset\bJ$ projects to  a fiber $J_a$ over some rational point  $a\in\PP^1$. 
Thus it suffices to prove the following:

\begin{proposition}
\mylabel{bijection irreducible components}
Let $\ba:\Spec(\Omega)\ra \PP^1$ be a geometric point. Assume that the image $a\in \PP^1$ is rational
and that the fiber $J_a$ is singular. Then  the canonical map  $\Gamma(J_\ba)\ra \Gamma(J_a)$ 
is a graph isomorphism.
\end{proposition}

\proof
The assertion  is trivial if $J_\ba$ has Kodaira symbol $\I_1$ or $\II$.
Condition \refjac{semistable J} ensures that it also holds for $\I_n$ with $n\geq 2$.
Now suppose that  $J_\ba$ has one of the Kodaira symbol $\III,\IV,\ldots,\II^*$.
Our task is to verify that the  Galois action of $G=\Gal(\FF_2^\alg/\FF_2)$ on the dual graph $\Gamma=\Gamma(J_\ba)$
is trivial.
Since it is a tree, it suffices that check that all terminal vertices are fixed.
This obviously holds for the terminal vertex $v_0\in\Gamma$ that corresponds to 
the irreducible component   that intersects the base-change of the zero-section $O\subset J$.
If there is at most one further terminal vertex, it must be fixed as well.
This already settles the cases where the Kodaira symbol is   $\III$, $\III^*$ or $\II^*$.

If the Kodaira symbol is $\I_m^*$, we must have $m\geq 1$ by condition \refjac{I0* J}.
Besides  $v_0\in\Gamma$, there are three other terminal vertices. One has a shorter distance
to $v_0$ than the others, and thus must be fixed, and we have to deal with the remaining two 
terminal vertices.
For  Kodaira symbols $\IV$ and $\IV^*$ there are also two remaining terminal vertices besides $v_0\in\Gamma$.
In all these cases, we shall check  that there is a section $P\subset J$ 
whose base-change intersects an irreducible component   $\Theta\subset J_\ba$
corresponding to one of the two remaining terminal vertices.

To achieve this, first note that $\Num(\bJ)$
is unimodular, because  the  minimal model of $\bJ$ is the projective plane.
Consider the orthogonal complement $W\subset\Num(\bJ)$ of the fiber $J_\ba$ and the zero-section $O$.
These two curves have Gram matrix $(\begin{smallmatrix}0&1\\1&-1\end{smallmatrix})$, and we infer that
$W$ is unimodular. It is generated by the vertical curves disjoint from the zero-section,
together with the combinations $(P-O)-(1+(P\cdot O))J_\ba$ for $P\in \MW(\bJ/\bar{\PP}^1)$.
The latter have self-intersection $-1-2(P\cdot O)-1$, hence  the lattice $W$ is even.
Furthermore, we consider the sublattice $W_\ba\subset W$  generated by the components $\Theta_1,\ldots,\Theta_r\subset J_\ba$ disjoint
from the zero-section. Note that this  is a root lattice, endowed with a canonical basis.

Suppose now that $J_\ba$ has Kodaira symbol $\IV$ or $\IV^*$, and that 
all sections $P\subset J$ have a base change that
passes through the component $\Theta_0\subset J_\ba$ corresponding to $v_0\in\Gamma$.
Then $W_\ba\subset W$ is an orthogonal direct summand, hence unimodular. On the other hand, this root lattice has type $A_2$ or $E_6$,
with discriminant $d=3$, contradiction.

It remains to treat the symbols $\I_m^*$ with $1\leq m\leq 4$.
Now $W_\ba$ is a root lattice of type $D_{m+4}$, which again is  not unimodular.
Our lattices and their dual lattices are related by a commutative diagram
\begin{equation}
\label{various lattices}
\begin{CD}
W_\ba	@>>>	W	@>>>	\Num(\bJ)\\
@VVV		@VVV		@VVV\\
W_\ba^*	@<<<	W^*	@<<<	\Num(\bJ)^*,
\end{CD}
\end{equation}
where the vertical maps are given by $x\mapsto (y\mapsto (x\cdot y))$. 
Seeking a contradiction, we now assume that each section $P\subset J$ passes through
the component  $\Theta_0\subset J_\ba$ corresponding to $v_0\in\Gamma$, or the
component $\Theta_1\subset J_\ba$ corresponding to the terminal vertex  $v_1\in\Gamma$  that has shorter distance to $v_0$.
The latter indeed happens, because otherwise $W_\ba\subset W$ would be an orthogonal direct summand, hence unimodular,
contradiction. Now fix some $P$ passing through $\Theta_1$.  Under the maps in \eqref{various lattices}, the combination $(P-O)-(1+(P\cdot O))J_\ba\in W$
maps to the dual basis vector  $\Theta_1^*\in W_\ba^*$. Since the maps  respect the 
intersection pairing, we conclude that $(\Theta_1^*\cdot\Theta_1^*)\in \QQ$ actually belongs to $2\ZZ$.
On the other hand, a direct computation with the root lattice $W_\ba$ of type $D_{m+4}$ reveals $(\Theta_1^*\cdot\Theta_1^*)=\pm 1$,
contradiction.
\qed

\medskip
My original arguments for the cases $\IV,\IV^*$ and $\I_m^*$ in the above proof relied on Lang's Classification \cite{Lang 2000}, 
which gives the possible  configuration of singular fibers for
$\bJ$, together with the  Oguiso--Shioda Table \cite{Oguiso; Shioda 1991}  for Mordell--Weil lattices 
and the Shioda's explicit formula \cite{Shioda 1990} for the height pairing $\langle P,Q\rangle\in\QQ$.
This  method was already used in \cite{Schroeer 2021a}, Section 15. The above more elegant
and conceptual argument was suggested by one referee.

\section{Classification of certain jacobian fibrations}
\mylabel{Classification jacobian}

Throughout, we work over the ground field $k=\FF_2$.
Let $Y$ be a geometrically rational surface, and $\phi:Y\ra\PP^1$ be an elliptic fibration that
is relatively minimal and jacobian. It is then described by some Weierstra\ss{} equation 
$$
y^2+a_1xy+a_3y = x^3+a_2x^2+a_4x+a_6,
$$
where the coefficients $a_i\in \FF_2[t]$ are polynomials of degree $\deg(a_i)\leq i$, and the discriminant does not vanish.
The goal of this section is to classify those Weierstra\ss{} equations that might arise from
an Enriques surface with constant Picard scheme.
It turns out that there are exactly eleven such Weierstra\ss{} equations, up to change of coordinates and
automorphism of the projective line.
This is quite remarkable, because all in all there are $2^{21}=2097152$ possible Weierstra\ss{} equations.

We make the following assumptions on the surface $J$ and the jacobian fibration $\phi:J\ra\PP^1$:
\newcounter{J}\setcounter{J}{0}
\renewcommand{\theJ}{\rm \roman{J}}
\newcommand{\labelJ}[1]{\refstepcounter{J}\mylabel{#1}}
\newcommand{\refJ}[1] {{\rm (\ref{#1})}}

\begin{enumerate}
\item 
\labelJ{constant}
The proper scheme $J$ has constant   Picard scheme $\Pic_{J/k}$.
\item 
\labelJ{rational}
There is at most one rational point $a\in\PP^1$ whose fiber $J_a$   is semistable or supersingular.
\end{enumerate}

The first condition arises from Theorem \ref{jacobian constant pic}, and the second from Proposition \ref{jacobian fiber properties}.
With Propositions \ref{constant num over finite field}, \ref{canonical type rational points} and \ref{points on enriques} 
we see that   the following holds as well:

\begin{enumerate}
\setcounter{enumi}{2}
\item 
\labelJ{non-rational}
The fibers $J_\ba$ over non-rational points $a\in\PP^1$ have Kodaira $\I_0$, $\I_1$ or $\II$.
\item 
\labelJ{I0*}
There is no fiber $J_\ba$ with Kodaira symbol $\I_0^*$.
\item 
\labelJ{points}
The number of rational points is $\Card Y(\FF_2)=25$.
\end{enumerate}

For each rational point $a\in \PP^1$, let  
$r_a\geq 1$ and $n_a\geq 1$  be the respective numbers of irreducible components and rational points in $J_a$.
We start our analysis by observing that there is a large fiber, in the sense that there are many irreducible
components:

\begin{proposition}
\mylabel{six   components}
There is a rational point $a\in\PP^1$ whose fiber $J_a$ contains at least six irreducible components.
\end{proposition}

\proof
Seeking a contradiction, we   assume that $r_a\leq 5$ for the three rational points $a\in\PP^1$.
If $J_a$ is unstable we actually have $r_a\leq 3$, because the Kodaira symbol $\I_0^*$ is impossible,
and  get $n_a\leq 7$.
The same estimate holds if the fiber is smooth, because it is then isomorphic to one of the elliptic curves
$E_i$, with $1\leq i\leq 5$.
In the semistable case we get $n_a\leq 10$, and there is at most one rational point with
semistable fiber. Writing $a,b,c\in\PP^1$ for the three rational points,
we get the  contradiction $25=\Card Y(\FF_2)=n_a+n_b+n_c\leq 10 + 7 + 7= 24$.
\qed

\medskip
We shall see that this large fiber must be unstable. The next result is a  first step into this direction:

\begin{proposition}
\mylabel{unstable fiber}
There is a rational point $a\in\PP^1$ whose fiber $J_a$ is reducible and unstable.
\end{proposition}

\proof
Seeking a contradiction, we assume that all reducible fibers are semistable.
According to Proposition \ref{six components}, there is a reducible fiber with at least six components, 
say $J_a$.
It follows that $J_a$ has Kodaira symbol $\I_r$ for some $r\geq 6$,
contributing $n=2r\geq 12$ rational points. We also have the upper bound $r\leq 9$, because
of the Picard number $\rho(J)=10$.

For each of the  two other rational points $b,c\in \PP^1$, it follows that the fiber 
is either ordinary or with Kodaira symbol $\II$.
In the former case, the number of rational points is even, in the latter case it is odd.
Since $\Card Y(\FF_2)=25$ is odd, we may assume without restriction that $J_b=E_{2i}$
for $1\leq i\leq 2$, and  $J_c$ has Kodaira symbol $\II$.
This gives $25=\Card Y(\FF_2) = 2r +  2i+3$. The only possible solution is $r=9$ and $i=2$.
But by Lang's Classification \cite{Lang 2000}, the configuration of singular fibers $\II+\I_9$ does not occur.
\qed

\medskip
The generic fiber   of the jacobian fibration $\phi:J\ra\PP^1$
has a certain invariant $j(J/\PP^1)$ in the function field $ k(\PP^1)=\FF_2(t)$, which is also called
the \emph{functional $j$-invariant} and  can   be seen  as a morphism $j:\PP^1\ra\PP^1$.
It is now convenient to distinguish the cases that the invariant   is zero or non-zero.

\begin{proposition}
\mylabel{non-constant j-invariant}
Suppose the functional $j$-invariant   is non-zero.
Then there is  exactly one  unstable fiber.
If this fiber contains $r\leq 5$ irreducible components, the configuration of fibers over   rational points
must be $\III+E_4+ \I_8$. 
\end{proposition}

\proof 
First note that by Proposition \ref{unstable fiber} there is a rational point $a\in\PP^1$ whose fiber $J_a$ is unstable.
As observed by Lang, there is no other unstable fiber (\cite{Lang 2000}, beginning of Section 1).
Indeed,  a fiber $J_x$ is unstable if and only if $c_4\in \FF_2[T]$ vanishes at the point $x$ (\cite{Deligne 1975}, Proposition 5.1), and 
in our situation  $c_4=a_1^4$ with $\deg(a_1)\leq 1$.
Assume now that $J_a$   contains at most five irreducible components. Since the Kodaira symbol $\I_0^*$ does not occur,
$J_a$ has at most three irreducible components.
According to Proposition \ref{six components}, there must be a semistable fiber $J_b$ over some rational point $b\in\PP^1$,
with $6\leq m\leq 9$ irreducible components.
In light of Lang's Classification \cite{Lang 2000}, Section 2 the only possible configuration of singular fibers
is $\III+\I_6$ and $\III+\I_8$. Consequently, the fiber over the last rational point $c\in\PP^1$ is ordinary,
whence $J_c=E_{2n}$ with $1\leq n\leq 2$.
The respective number of rational points in the fiber is thus 
$n_a=5$ and $n_b=2m$ and $n_c=2n$. From $25=\Card Y(\FF_2)=n_a+n_b+n_c$ we get $m+n=10$.
The only solution is  $m=8$ and $n=2$.
\qed

\medskip
We now can state the first half our our classification result:

\begin{theorem}
\mylabel{classification non-zero j}
Suppose the functional $j$-invariant is non-zero. 
Up to coordinate changes and automorphisms of the projective line,
the fibration $\phi:J\ra\PP^1$ is given by exactly one of the following eight Weierstra\ss{} equations: 
$$
\begin{array}[c]{@{}lll}
\toprule
\text{\rm Weierstra\ss{} equation}      & \text{\rm fibers over $t=0,1,\infty$}	& \text{\rm $j$-invariant}\\
\toprule
y^2+txy+t^2y=x^3+ tx^2+t^4x+t^5(1+t)    & \I_1^*+E_4+\I_4                       & t^4\\
\midrule
y^2+txy=x^3+t^3x+t^5(1+t)               & \III^*+E_4+E_4                        & t^2/(t^2+t+1)\\
\midrule
y^2+txy=x^3+t^3x 			            & \III^*+E_4+\I_2                       & t^2\\
y^2+txy=x^3+t^2x^2+t^3x                 & \III^*+E_2+\tI_2\\
\midrule
y^2+txy=x^3+t^5                         & \II^*+E_4+\I_1                        & t\\
y^2+txy=x^3+t^2x^2+t^5                  & \II^*+E_2+\tI_1\\
\midrule
y^2+txy=x^3+tx^2+t^4x                   & \I_4^* +E_2+E_4                       & 1\\
\midrule
y^2+txy+ty=x^3+tx^2+tx                     & \III+E_4+\I_8                         & t^8\\
\bottomrule
\end{array}
$$
For the invariant $j=t^2/(t^2+t+1)$, there is an additional singular fiber, which occurs over the point $x\in\PP^1$ with $\kappa(x)\simeq\FF_4$
and has Kodaira symbol $\I_1$.
In all other cases, the fibers over the non-rational points are smooth.
\end{theorem}

\proof
Let $J_a$ be the unstable fiber, with $1\leq r_a\leq 9$ irreducible components and $n_a=2r_a+1$ rational points. 
Note that the latter is odd.
Without restriction, we may assume that $a$ is given by $t=0$. Write $b,c\in\PP^1$ for the remaining two rational points,
whence
\begin{equation}
\label{cardinality}
25=\Card Y(\FF_2) = n_a+n_b+n_c.
\end{equation}

First, assume that $r_a\geq 6$.
To start with, we consider the situation that both fibers $J_b$ and $J_c$ are smooth.
Applying an automorphism of the projective line, we may assume that the first is ordinary. Then $n_b$ is even.
By \eqref{cardinality} the number $n_c$ is also even, hence $J_c$ is also ordinary.
Write  $J_b=E_{2i}$ and $J_c=E_{2j}$. Applying another automorphism of the projective line,  we may assume
that $b$ is given by $t=1$. Recall that $r_a\geq 1$ is the number of irreducible components in $J_a$.
 From \eqref{cardinality} one infers $r_a+i+j=12$.
Using $i+j\leq 4$ we obtain $r_a\geq 8$.
Consequently, the possible configurations of   fibers over rational points are 
$$
\II^*+E_2+E_4,\quad \I_4^*+E_2+E_4,\quad \III^*+E_4+E_4,\quad \I_3^*+E_4+E_4.
$$
In the first and last cases, Lang's Classification \cite{Lang 2000} tells us that
there must be an additional singular geometric fiber, which is unique and has Kodaira symbol $\I_1$.
This uniqueness ensures that it comes from a  fiber over a rational point, contradiction.

We now analyze the case that $J_a$ has Kodaira symbol $\I_4^*$.
We are thus in Lang Case 5D from \cite{Lang 2000}, and the Weierstra\ss{} equation takes the form
$y^2+txy=x^3+td_1x^2+t^4x$, where $t\nmid d_1$.
\emph{In this context, symbols like $d_i,e_i,\ldots$ denote   polynomials in the variable $t$ of degree $\leq i$.}
Here we have the two cases $d_1=1$ and $d_1=1+t$. The former is transformed into the latter, 
by the automorphism $t\mapsto t/(1+t)$ of the projective line,  followed by the change of coordinates 
$x=u^2x'$ and $y=u^3y'$ with $u=1/(1+t)$.
For $d_1=1$ the fibers, the  Weierstra\ss{} equation and the ensuing $j$-invariant are as in the table.

Suppose now that $J_a$ has Kodaira symbol $\III^*$.
Here we are in Lang Case 7, 
and the configuration of singular geometric fibers  must
be $\III^*+\I_1+\I_1$. The Weierstra\ss{} equation is of the form
$y^2+txy+t^3c_0y=x^3+t^2d_0x^2+ t^3e_1x+t^5f_1$ with $t\nmid e_1$.
Making a change of variables we achieve $c_0=0$ and $e_1=1$.
The discriminant becomes $\Delta=t^{10}(tf_1+1)$.
Since the fiber over  the unique point $x\in\PP^1$ with residue field $\kappa(x)\simeq \FF_4$
is singular, we must have $f_1=1+t$.
The case $d_0=1$ yields $J_b=E_2$, contradiction. Thus $d_0=0$.
Again the Weierstra\ss{} equation and the ensuing $j$-invariant are as in the table.

We next consider the situation that one of the fibers over rational points is semistable.
Applying an automorphism of the projective line, we may assume that this happens over $c\in\PP^1$,
and that this point is given by $t=\infty$. Recall that $r_c\geq 1$ is the number of irreducible components
in $J_c$. The remaining fiber is ordinary,
and we write $J_b=E_{2i}$ for some $1\leq i\leq 2$, such that $n_b=2i$.
Suppose first that the semistable fiber is untwisted.
Then $n_c=2r_c$, and \eqref{cardinality} gives $r_a+i+r_c=12$. In particular, we have
$$
r_c = 12-i-r_a\geq 10-r_a.
$$
For $r_a=6$ this estimate gives $r_c\geq 4$. By Lang's Classification, the configuration of fibers over rational points
must be $\I_1^*+E_4+\I_4$. For $r_a=7$ we get  $r_a\geq 3$, and now the classification ensures that the configuration must be $\IV^*+E_4+\I_3$.
In this case, however, there is a unique additional degenerate fiber, which has Kodaira symbol $\I_1$.
The uniqueness ensures that it lies over a rational point, contradiction.
For $r_a=8$ we obtain $r_c\geq 2$, and the configuration is $\III^*+E_4+\I_2$.
Finally, for $r_a=9$ the possible configuration is $\II^*+E_4+\I_1$.
If the semistable fiber is twisted, we have  $r_c\leq 2$ and $n_c=2r_c+2$.
Now \eqref{cardinality} becomes $r_a+i+r_c=11$ and thus $r_c\geq 11-i-r_a\geq 7$, and the argument is similar.
Summing up, we have the following possibilities:
$$
\I_1^*+E_4+\I_4,\quad \III^*+E_4+\I_2,\quad\III^*+E_2+\tI_2, \quad \II^*+E_4+\I_1,\quad \II^*+E_2+\tI_1.
$$
Let me give the details for the first case: Here we are in Lang Case 5A,
with a Weierstra\ss{} equation of the form $y^2+txy+t^2c_1y=x^3+td_1x^2+t^3e_1x+t^4c_2$,
subject to the conditions  $t\nmid c_1$ and $t\nmid d_1$.
Making a change of coordinates, we may assume $c_1=1$ and $e_1=t$.
Now the discriminant becomes
$$
\Delta=t^8(t^2c_2+t^3+td_1+t^4+1+t).
$$
The second factor is a unit, hence equals one, because we assume that the semistable fiber is at $t=\infty$.
Comparing coefficients yields the solutions  $c_2=t+t^2$, $d_1=1$ and $c_2=1+t+t^2$, $d_1=1+t$.
For the latter solutions, the semistable fiber is twisted, as revealed by the Tate Algorithm at $t=\infty$.
Hence we are in the former case, and the Weierstra\ss{} equation and the resulting $j$-invariant is as in the Table.
The remaining cases are even simpler, and handled in the same fashion. We leave the details to the reader.
Summing up, this settles the case $r_a\geq 6$.

It remains to treat the case that $r_a\leq 5$. According to Proposition \ref{non-constant j-invariant}, the configuration over
the singular fibers is $\III+E_4+\I_8$, and we are in Lang Case 2A, with Weierstra\ss{} equation
$y^2+txy + tc_2y=x^3+tc_1x^2+tc_3x+t^2c_4$, with $t\nmid c_2$ and $t\nmid c_3$ and $\Delta=t^4$.
Making successive changes  of coordinates  $x=x'+R$ and $y=y'+T$ 
where the indeterminate $t$ divides $R$ and $T$, 
we may assume that $c_2=c_3=1$.
Then $1=\Delta/t^4=c_4t^4+t^3+c_1t^3+1$. It follows $c_1=1+\alpha t$ and $c_4=\alpha$
for some scalar $\alpha\in \FF_2$. In case $\alpha=1$ we have $J_b=E_2$, contradiction, hence $\alpha=0$.
\qed

\medskip
It remains to understand the situation with trivial functional $j$-invariant:

\begin{theorem}
\mylabel{classification  zero j}
Suppose the functional $j$-invariant is zero.
Up to some coordinate changes and automorphisms of the projective line,
the fibration $\phi:J\ra\PP^1$ is given by exactly one of the following three Weierstra\ss{} equations: 
$$
\begin{array}[c]{@{}ll}
\toprule
\text{\rm Weierstra\ss{} equation}		& \text{\rm fibers over $t=0,1,\infty$}\\
\toprule
y^2+t^2y=x^3+tx^2			& \I_1^*+E_5+\IV\\	
\midrule
y^2+t^2y=x^3+t^3x 			& \IV^*+E_5+\III\\
\midrule
y^2+t^2y=x^3 			& \IV^*+E_3+\IV\\
\bottomrule
\end{array}
$$
In all three cases, the fibers over non-rational points $x\in\PP^1$ are smooth.
\end{theorem}

\proof
There are neither ordinary nor semistable fibers, because the  map $j:\PP^1\ra\PP^1$ misses
the values $t=1$ and $t=\infty$. In particular, each smooth fiber  over a rational point 
is isomorphic to one of the elliptic curves $E_1$, $E_3$ or $E_5$. So there are
at least two rational points $a\neq c$ on $\PP^1$ such that $J_a$ and $J_c$ are unstable.
By Lang's Classification, there are at most three singular geometric fibers. If present, 
the third singular geometric fiber must be Galois invariant, whence occurs over
the remaining rational point $b\in\PP^1$. Thus all fibers over non-rational points are smooth.

Suppose first that all fibers of $\phi:J\ra\PP^1$ are reduced, such that the Kodaira symbols for
the singular fibers are $\II$, $\III$ or $\IV$. Then each fiber over  a rational point contains
at most seven rational points, which gives the contradiction $25=\Card Y(\FF_2)\leq 3\cdot 7$.
Thus there is a non-reduced rational fiber.

Recall that $J_a$ and $J_c$ are unstable. Let $r_a,r_c\geq 1$ be the respective number of irreducible components.
After applying an automorphism of the projective line, we may assume that $J_a $  is non-reduced. 
By Lang's Classification (\cite{Lang 2000}, page 5829), the remaining fiber $J_b$ must be smooth 
(one referee pointed out that this also follows from the fact that the valuation of the discriminant
in characteristic $p=2$ is at least four at an additive fiber, and at least eight at a non-reduced fiber).
Write $J_b=E_{2i-1}$ for some $1\leq i\leq 3$.
Then $25=\Card Y(\FF_2) = (2r_a+1) +(2i-1) + (2r_c+1) $, or equivalently $r_a+i+r_c=12$.
Since $J_a$ is non-reduced and the Kodaira symbol $\I_0^*$ does not occur, we have $6\leq r_a\leq 9$.
In particular, $r_c\leq 5$, so the unstable fiber $J_c$ is reduced, such that actually $1\leq r_c\leq 3$.
The case $r_a=6$ implies $r_c=i=3$. In case $r_a=7$ we get $r_c=3$, $i=2$ or $r_c=2$, $i=3$.
Going through Lang's Classification, we see that these are the only possibilities,
which already give the second column of the table.

Suppose the configuration of fibers over rational points is $\I_1^*+E_5+\IV$.
Then we are in Lang Case 13A, and the Weierstra\ss{} equation takes the form
$y^2+t^2c_1y=x^3+td_1x^2+t^3e_1x+t^4d_2$ with $t\nmid d_1$ and $t\nmid c_1$.
The discriminant is $\Delta=t^8c_1^4$, thus $c_1=1$. Making a change of coordinate $y'=y+sx$,
we achieve $d_1=1$. Making a further change of coordinates, we may assume
that $e_1$ has no linear term, and that $d_2$ has neither constant nor quadratic term.
Now the Weierstra\ss{} equation becomes $y^2+t^2y=x^3+tx^2+t^3e_0x+t^5d_0$.
One checks that only $e_0=d_0=0$ yields $J_b=E_5$.

Next, suppose the configuration is $\IV^*+E_5+ \III$.
Now we are in Lang Case 14, with $y^2+t^2c_1y=x^3+t^2d_0x^2+t^3e_1x+t^4d_2$ subject to
the condition $t\nmid c_1$. The discriminant is $\Delta=t^8c_1^4$, whence $c_1=1$.
Making a suitable change of variables $y=y'+sx'$ we   achieve $d_0=0$.
As in the preceding paragraph, we put the Weierstra\ss{} equation into
the form $y^2+t^2y=x^3+t^3e_0x+t^5d'_0$.
Only the case $e_0=1$ and $d'_0=0$ yields $J_b=E_5$.

Finally, consider the configuration $\IV^*+E_3+\IV$.
We argue as in the preceding paragraph to reduce to the very same   Weierstra\ss{} equation.
Now only $e_0=d'_0=0$ yields the fiber $J_b=E_3$. 
\qed

\section{Possible configurations for Enriques surfaces}
\mylabel{Possible configurations}

Let $Y$ be an Enriques surface over the prime field $k=\FF_2$ with constant Picard scheme. 
According to Theorem \ref{family with constant pic}, it admits at least one genus-one fibration.
Let us collect the information on   fibers we have gathered so far:

\begin{proposition}
\mylabel{possible configurations}
Let $\varphi:Y\ra\PP^1$ be a genus-one fibration.
Up to  automorphisms of the projective line,
the possible configurations of fibers over the the rational points $t=0,1,\infty$ are given by
the following table:
\begin{equation}
\label{first restriction}
\begin{array}[c]{@{}lll}
\toprule
\text{elliptic}		& 		& \text{quasielliptic}\\
\midrule
\II^*+E_4+\I_1		&		& \II^*+\II+\II\\
\II^*+E_2+\tI_1		&		& \III^*+\III+\II\\	 
\III^*+E_4+\I_2\\	 
\III^*+E_2+\tI_2\\
\III^*+E_4+E_4\\
\IV^*+E_3+\IV\\		 
\IV^*+E_5+\III\\
\midrule
\I_4^* +E_4+E_2		&		& \I_4^*+\II+\II\\
\I_1^*+E_5+\IV      &		& \I_2^*+\III+\III\\
\I_1^*+E_4+\I_4\\
\midrule
\III+E_4+\I_8\\
\bottomrule
\end{array}
\end{equation}
If moreover the Enriques surface $Y$ is non-exceptional, the following holds:
\begin{enumerate}
\item
Each fiber  with Kodaira symbol  $\IV$, $\IV^*$,   $\II^*$, $\I_1^*$ is multiple.
\item
If there is a simple fiber with Kodaira symbol $\III^*$, the fibration is elliptic, and the configuration
of fibers over rational points is $\III^*+E_4+E_4$.
\item
For each multiple fiber $\varphi^{-1}(a)$ with Kodaira symbol $\IV^*$, $\III^*$ or $\II^*$, there is  no $(-2)$-curve
$R\subset Y$ with intersection number  $\varphi^\inv_\ind(a)\cdot R=1$.
\end{enumerate}
\end{proposition}

\proof
The main task is to verify the table.
Let $\bar{a}:\Spec(\Omega)\ra \PP^1$ be a geometric point, with image $a\in \PP^1$.
If this image point is non-rational, the geometric fiber $Y_\ba$ is irreducible, 
according to Theorem \ref{geometric and schematic fiber}.
In any case,  the map $\Gamma(Y_\ba)\ra\Gamma(Y_a)$ between dual graphs is a graph isomorphism.
Moreover, the Kodaira symbol $\I_0^*$ is impossible, according to Proposition \ref{canonical type rational points}.
In particular, $Y_a$ and $Y_\ba$ have the same number of irreducible components. Write $r_a\geq 1$ for this number.
Then we have the inequality $2+\sum (r_a-1)\leq \rho(Y)=10$, and thus  $r_0+r_1+r_\infty\leq 11$.

Suppose first that the fibration is quasielliptic.  Since the Picard group $\Pic^0(Y_F)$ for the generic fiber  
is annihilated by $p=2$, we actually have  $r_0+r_1+r_\infty=11$. The possible Kodaira symbols
and
the corresponding number of irreducible components  are as follows (\cite{Cossec; Dolgachev 1989}, Theorems 5.7.4--5.7.6):
$$
\begin{array}{lllllll}
\toprule
\text{Kodaira symbol}			& \II	& \III	& \I_2^*	& \III^*	& \II^*	& \I_4^*\\
\midrule
\text{number of irr.\ comp.}	& 1		& 2		& 7			& 8			& 9		& 9\\
\bottomrule
\end{array}
$$
Since $9+1+1=8+2+1=7+2+2=11$ are the only possible solutions, the second column of table \eqref{first restriction} follows.

Now suppose that $\varphi:Y\ra\PP^1$ is elliptic, and consider the resulting jacobian fibration $\phi:J\ra\PP^1$.
Then the smooth proper surface $J$ is geometrically rational (Proposition \ref{jacobian geometrically rational}), 
with constant Picard scheme (Theorem \ref{jacobian constant pic}).
As above, the geometric fibers $J_\ba$ are irreducible if $a\in \PP^1$ is non-rational,
and in any case $\Gamma(J_\ba)\ra\Gamma(J_a)$ is a graph isomorphism.
According to Theorems \ref{classification non-zero j} and \ref{classification zero j}, 
the first column of table \eqref{first restriction} gives the possible configurations
of singular fibers over the  rational points for $\phi:J\ra\PP^1$. 
By Proposition \ref{same symbols} the fibers $Y_\ba$ and $J_\ba$ have the same Kodaira symbols. 
Moreover, we  actually have $J_a\simeq Y_a$ if the fiber $Y_a$ is simple. In particular, this holds for all semistable
fibers. According to Proposition \ref{reduction isomorphic}, we also have $J_a\simeq (Y_a)_\red$ provided that $Y_a$ is multiple 
with smooth reduction. In turn, the first column of table \eqref{first restriction} indeed describes the configuration of
degenerate fibers for $\varphi:Y\ra\PP^1$.

Now suppose that $Y$ is non-exceptional. If $\varphi:Y\ra\PP^1$ is elliptic 
having some $\varphi^{-1}(a)$ with Kodaira symbol $\IV,\IV^*,\II^*,\I_1^*$ there must   be 
some $\varphi^\inv_\ind(b)$ that is supersingular or semistable, in light of  the table,
so $\varphi^{-1}(a)$ is multiple by Lemma \ref{multiple for p=2}.
If $\varphi:Y\ra\PP^1$ is quasielliptic, fibers with Kodaira symbol $\II^*$ are multiple 
according to    Proposition \ref{properties non-exceptional}. This establishes (i).

Next assume that $\varphi^{-1}(a)$ is simple with Kodaira symbol $\III^*$.
Again by Proposition \ref{properties non-exceptional} the fibration must be elliptic. By the table, the configuration
of fibers over rational points is $\III^*+E_4+\I_2$ or  $\III^*+E_2+\tI_2$ or $\III^*+E_4+E_4$.
The former are impossible by Lemma \ref{multiple for p=2}, which gives (ii).
Finally, assertion (iii) is a consequence of    Proposition \ref{properties non-exceptional}.
\qed

\medskip
Starting from the above information, we rule out step by step each possible Kodaira symbol, 
mostly from combinatorial considerations.  In this section, we will reduce
from fifteen to eight possible configurations for $\varphi:Y\ra\PP^1$.

\begin{proposition}
\mylabel{no I8}
Assumption as in Proposition \ref{possible configurations}.
Then there are no fibers with   Kodaira symbol $\I_8$. 
\end{proposition}

\proof
Seeking a contradiction, we assume that $\varphi^{-1}(\infty)$ has Kodaira symbol $\I_8$, and $\varphi^{-1}(0)$ has symbol $\III$.
Then   the configuration of reducible fibers has dual graph:
%
%
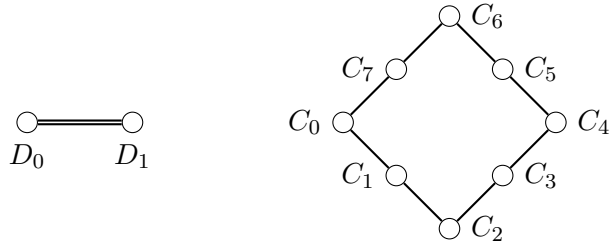
\begin{equation*}
\begin{gathered}
\begin{tikzpicture}
[node distance=1cm, font=\small]
\tikzstyle{vertex}=[circle, draw, fill=white, inner sep=0mm, minimum size=1.5ex]
\node[vertex]	(C0)  	at (0,0)[label=left:{$C_0$}] 		{};
\node[vertex]	(C1)			[below right of=C0, label=left:{$C_1$}]	{};
\node[vertex]	(C2)			[below right of=C1, label=right:{$C_2$}]	{};
\node[vertex]	(C3)			[above right of=C2, label=right:{$C_3$}]	{};
\node[vertex]	(C4)			[above right of=C3, label=right:{$C_4$}]	{};
\node[vertex]	(C5)			[above left of=C4, label=right:{$C_5$}]	{};
\node[vertex]	(C6)			[above left of=C5, label=right:{$C_6$}]	{};
\node[vertex]	(C7)			[below left of=C6, label=left:{$C_7$}]	{};
 
\node[]         (E)             [left of=C0,  xshift=-0.4cm,  label=below:{$$}]	{};
\node[vertex]	(D1)			[left of=E,  xshift=-0.4cm,  label=below:{$D_1$}]	{};
\node[vertex]	(D0)			[left of=D1,  xshift=-0.4cm, label=below:{$D_0$}]	{};
 
\draw[thick] (C0)--(C1)--(C2)--(C3)--(C4)--(C5)--(C6)--(C7)--(C0);
\draw[thick,double] (D0)--(D1);
\end{tikzpicture}
\end{gathered}
\end{equation*}
\emph{Here the double edge indicates that the   scheme $D_0\cap D_1$ has length two}.
The semistable fiber $\varphi^{-1}(\infty)$  is simple whereas the unstable fiber $\varphi^{-1}(0)$  is multiple.
Moreover,  there is a genus-one fibration $\psi:Y\ra\PP^1$ so that the multiple fibers $2F$ have  
$F\cdot \varphi^\inv_\ind(0)=1$. 
After renumeration, we may assume   $(F\cdot D_1)=1$, hence $D_0$ is vertical with respect to $\psi$.
Without restriction, we may assume that $\psi^\inv(\infty)$ is simple, whereas the $\psi$-fibers over $a=0,1$ are multiple.
Write $\psi^\inv_\ind(0)=\sum_{i=0}^r m_i \Theta_i$. After renumeration, we may assume $(D_1\cdot \Theta_0)>0$.
From
$$
1=\varphi^\inv_\ind(0)\cdot\psi^\inv_\ind(0) = D_1\cdot \sum_{i=0}^rm_i\Theta_i = m_0(D_1\cdot\Theta_0) + \sum_{i=1}^r m_i(D_1\cdot\Theta_i)
$$
we see that $(D_1\cdot\Theta_0)=m_0=1$, and that the curve $\Theta_1+\ldots+\Theta_r$  is disjoint from $D_1$.
Furthermore, $D_0$ does not belong to $\psi^\inv(0)$, because   $D_0=\Theta_j$ gives the contradiction
$$
1=\varphi^\inv_\ind(0)\cdot\psi^\inv_\ind(0) = D_1\cdot\psi^\inv_\ind(0)\geq D_1\cdot m_j\Theta_j\geq 2m_j\geq 2.
$$
Consequently $\Theta_1+\ldots+\Theta_r\subset\varphi^{-1}(\infty)$.
From $m_0=1$ we infer that  $\Theta_1+\ldots+\Theta_r$ is connected, hence is a chain. So the multiple fiber $\psi^\inv_\ind(0)$ is either
ordinary or has Kodaira symbol $\II$, $\III$ or $\IV$.
The same analysis applies to the other multiple fiber $\psi^\inv_\ind(1)$.
From Proposition \ref{possible configurations}, we see that the possible Kodaira symbols for  the simple fiber $\psi^\inv(\infty)$
are $\I_8$ or $\II^*, \III^*, \I_4^*, \I_2^*$.

To proceed write  $\psi^\inv(\infty)=\sum_{i=0}^s n_i\Upsilon_i$. Then $s\geq 6$ and  $(\Upsilon_i\cdot\Upsilon_j)\leq 1$. After renumeration $\Upsilon_0=D_0$. From
$$
2=\varphi^\inv_\ind(0)\cdot\psi^\inv(\infty) = D_1\cdot\sum_{i=0}^s n_i\Upsilon_i = 2n_0 + \sum_{i=1}^s n_i(D_1\cdot\Upsilon_i)
$$
we infer $n_0=1$, and that $\Upsilon_1+\ldots+\Upsilon_s$ is disjoint from $D_1$.
Seeking a contradiction, we now assume that $\psi^{-1}(\infty)$ is additive. By inspection, one sees that  $\Upsilon_0$ has exactly one neighbor in the dual graph,
say $\Upsilon_1$, and  that $\Upsilon=\Upsilon_2+\ldots+\Upsilon_s$ is not a disjoint union of chains.
Moreover, $\Upsilon$ is disjoint from $\varphi^{-1}(0)$, hence must be strictly 
contained in $\varphi^\inv(\infty)$. The latter is multiplicative, and it follows that $\Upsilon$ is a union of chains, contradiction.
We conclude that $\psi^\inv(\infty)$ has Kodaira symbol $\I_8$.

Summing up, the configuration of rational fibers for $\psi$
must be $\III+E_4+\I_8$, according to the   table in Proposition \ref{possible configurations}.
Without restriction we may assume that $\psi^{-1}_\ind(0)=\Theta_0+\Theta_1$ is reducible.
We saw above that $\Theta_1$ belongs to $\varphi^\inv(\infty)$, so after renumeration $\Theta_1=C_0$.
To simplify notation set $E=\Theta_0$, such that
$\psi^\inv_\ind(0)=E+C_0$   with $(E\cdot D_1)=1$ and $(E\cdot C_0)=2$. 
It follows that  $D_0$ and $C_2+\ldots+C_6$ belong  to $\psi^\inv(\infty)$. Let $E_1,E_7$ be the additional components.
The dual graph for our curves takes the following form:
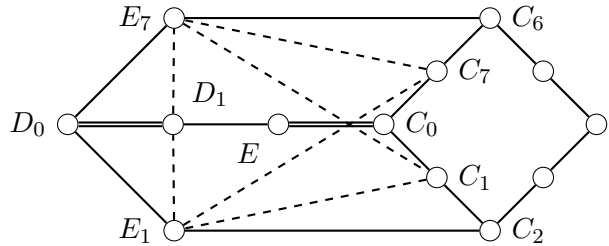
\begin{equation*}
\begin{gathered}
\begin{tikzpicture}
[node distance=1cm, font=\small]
\tikzstyle{vertex}=[circle, draw, fill=white, inner sep=0mm, minimum size=1.5ex]
\node[vertex]	(C0)  	at (0,0)[label=right:{$C_0$}] 		{};
\node[vertex]	(C1)			[below right of=C0, label=right:{$C_1$}]	{};
\node[vertex]	(C2)			[below right of=C1, label=right:{$C_2$}]	{};
\node[vertex]	(C3)			[above right of=C2, label=right:{}]	{};
\node[vertex]	(C4)			[above right of=C3, label=right:{}]	{};
\node[vertex]	(C5)			[above left of=C4, label=right:{}]	{};
\node[vertex]	(C6)			[above left of=C5, label=right:{$C_6$}]	{};
\node[vertex]	(C7)			[below left of=C6, label=right:{$C_7$}]	{};
 
\node[vertex]   (E)             [left of=C0,  xshift=-0.4cm,  label=below left:{$E$}]	{};
\node[vertex]	(D1)			[left of=E,  xshift=-0.4cm,  label=above right:{$D_1$}]	{};
\node[vertex]	(D0)			[left of=D1,  xshift=-0.4cm, label=left:{$D_0$}]	{};

\node[vertex]   (E7)            [left of=C6,  xshift=-3.2cm,  label=left:{$E_7$}]	{};
\node[vertex]   (E1)            [left of=C2,  xshift=-3.2cm,  label=left:{$E_1$}]	{};

\draw[thick] (C0)--(C1)--(C2)--(C3)--(C4)--(C5)--(C6)--(C7)--(C0);
\draw[thick,double] (D0)--(D1);
\draw[thick,double] (E)--(C0);
\draw[thick] (E)--(D1);
\draw[thick] (C2)--(E1)--(D0)--(E7)--(C6);
\draw[thick,dashed] (C7)--(E7)--(C1);
\draw[thick,dashed] (C1)--(E1)--(C7);
\draw[thick,dashed] (E1)--(D1)--(E7);
\end{tikzpicture}
\end{gathered}
\end{equation*}
\emph{Here the dashed edges indicate potential intersections}.
Using
$$
2=(D_0+D_1) \cdot\sum \psi^\inv(\infty) = (D_0+D_1)\cdot(E_1+E_7) = 1+1 + D_1\cdot(E_1+E_7)
$$
we conclude that $D_1$ is actually disjoint from $E_1\cup E_7$. Then
$$
2=2(D_0+D_1)\cdot E_1 = (C_0+\ldots+C_7)\cdot E_1 = (C_1+C_7)\cdot E_1 +1
$$
ensures that either $(E_1\cdot C_1)=1$, $(E_1\cdot C_7)=0$ or $(E_1\cdot C_7)=1$, $(E_1\cdot C_1)=0$.
In the former case we get a curve    $E_1+C_2+C_1$ of canonical type $\I_3$, in the letter 
$E_1+C_2+\ldots+C_7$ of canonical type $\I_7$, contradiction.
\qed

\medskip
For the next observation it is crucial to  exclude exceptional Enriques surfaces.

\begin{proposition}
\mylabel{no II*, IV*, IV}
Assumption as in Proposition \ref{possible configurations}. Assume  $Y$ is non-exceptional. 
Then there are no fibers with Kodaira symbol $\II^*$, $\IV^*$ or  $\IV$. Furthermore, there are no multiple fiber
with Kodaira symbol $\I_4^*$.
\end{proposition}

\proof
According to Proposition \ref{possible configurations}, all fibers with Kodaira symbol $\II^*$, $\IV^*$, $\IV$ are multiple.
Seeking a contradiction, we thus assume that $\varphi^{-1}(\infty)$ is a multiple with  Kodaira symbol $\II^*$, $\IV^*$, $\IV$ 
or  $\I_4^*$.
 
By  Proposition \ref{properties non-exceptional}, there is another genus-one fibrations $\psi:Y\ra\PP^1$ 
such that  $\varphi^{-1}(\infty)\cdot \psi^{-1}(\infty)=4$.
Without restriction, we may assume that $\psi^{-1}(\infty)$ is the   fiber  having
the maximal number of irreducible components.
Let  $2F$  a  multiple fiber  for $\psi$, such that  $F\cdot \varphi_\ind^{-1}(\infty)=1$.
 
We first treat the case   that  
$\varphi^{-1}(\infty)$ is multiple with Kodaira symbol $\II^*$. With enumeration as in the Bourbaki tables
\cite{LIE 4-6}, the dual graph for such a  fiber is:
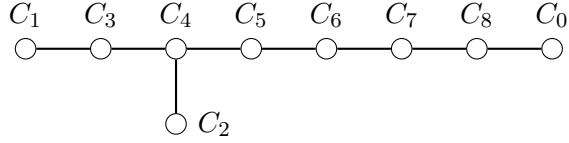
\begin{equation*}
\begin{gathered}
\begin{tikzpicture}
[node distance=1cm, font=\small]
\tikzstyle{vertex}=[circle, draw, fill=white, inner sep=0mm, minimum size=1.5ex]
\node[vertex]	(C1)  	at (1,0) 	[label=above:{$C_1$}] {};
\node[vertex]	(C3)			[right of=C1,  label=above:{$C_3$}]	{};
\node[vertex]	(C4)			[right of=C3,  label=above:{$C_4$}]	{};
\node[vertex]	(C2)			[below of=C4,  label=right:{$C_2$}]	{};
\node[vertex]	(C5)			[right of=C4,  label=above:{$C_5$}]	{};
\node[vertex]	(C6)			[right of=C5,  label=above:{$C_6$}]	{};
\node[vertex]	(C7)			[right of=C6,  label=above:{$C_7$}]	{};
\node[vertex]	(C8)			[right of=C7,  label=above:{$C_8$}]	{};
\node[vertex]	(C0)			[right of=C8,  label=above:{$C_0$}]	{};

\draw[thick]    (C1)--(C3)--(C4)-- (C5)--(C6)--(C7)--(C8)--(C0);
\draw[thick]    (C2)--(C4);
\end{tikzpicture}
\end{gathered}
\end{equation*}
We also denote this dual graph by $T_{2,3,6}$. \emph{In general, the symbol $T_{l,m,n}$ refers to a star-shaped graph
with three terminal chains of length $l,m,n$, where the central vertex belongs to each of the terminal chains}.

We have $\varphi^{-1}_\ind(\infty)=C_0+2C_1+3C_2+\ldots +2C_8$, and $C_0$ is indeed the only component  with multiplicity $m=1$.
Thus  $(C_0\cdot F)=1$.  
In turn, the connected curve $C_1+\ldots+C_8$ has dual graph $T_{2,3,5}$ and is necessarily contained in
the fiber $\psi^{-1}(\infty)$, which therefore must have Kodaira symbol $\II^*$. This is multiple by 
Proposition \ref{possible configurations}.
Let $R\subset\psi^{-1}(\infty)$ be the additional component, with $(R\cdot C_8)=1$.
Then 
$$
1=\varphi^{-1}_\ind(\infty)\cdot\psi^{-1}_\ind(\infty) =  \varphi^{-1}_\ind(\infty)\cdot R= (C_0+2C_8)\cdot R\geq 2(C_8\cdot R)=2,
$$
contradiction. Summing up, there are no fibers with Kodaira symbol $\II^*$.

Suppose next that $\varphi^{-1}(\infty)$ is multiple with Kodaira symbol $\I_4^*$.
The dual graph is:
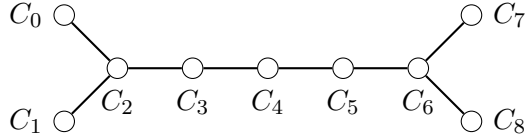
\begin{equation*}
\begin{gathered}
\begin{tikzpicture}
[node distance=1cm, font=\small]

\tikzstyle{vertex}=[circle, draw, fill=white, inner sep=0mm, minimum size=1.5ex]
\node[vertex]	(C2)  	at (1,0) 	[label=below:{$C_2$}] 		{};
\node[vertex]	(C3)			[right of=C2, label=below:{$C_3$}]	{};
\node[vertex]	(C4)			[right of=C3, label=below:{$C_4$}]	{};
\node[vertex]	(C5)			[right of=C4, label=below:{$C_5$}]	{};
\node[vertex]	(C6)			[right of=C5, label=below:{$C_6$}]	{};

\node[vertex]	(C7)			[above right of=C6, label=right:{$C_7$}]	{};
\node[vertex]	(C8)			[below right of=C6, label=right:{$C_8$}]	{};
\node[vertex]	(C0)			[above left  of=C2, label=left :{$C_0$}]	{};
\node[vertex]	(C1)			[below left  of=C2, label=left :{$C_1$}]	{};

\draw [thick] (C2)--(C3)--(C4)--(C5)--(C6);
\draw [thick] (C0)--(C2)--(C1);
\draw [thick] (C7)--(C6)--(C8);
\end{tikzpicture}
\end{gathered}
\end{equation*}
Now $\varphi^{-1}_\ind(\infty) = C_0+C_1+2(C_2+\ldots+C_6)+C_7+C_8$.
After renumeration, we may assume $(F\cdot C_0)=1$.
In turn, the connected curve $C_1+\ldots+C_8$ has dual graph $T_{2,2,6}$, and is necessarily contained
in $\psi^{-1}(\infty)$. This must have Kodaira symbol $\I_4^*$, as we already ruled out the symbol $\II^*$.
Let $R\subset \psi^{-1}(\infty)$ be the additional component, with $(R\cdot C_2)=1$.
We then have
$$
2\geq \varphi^{-1}_\ind(\infty)\cdot \psi^{-1}_\ind(\infty) = \varphi^{-1}_\ind(\infty)\cdot R = (C_0+2C_2)\cdot R=(C_0\cdot R)+2.
$$
Consequently $C_0\cap R=\varnothing$. In turn, $C_0+\ldots+C_3+R$ supports a curve of canonical type $\I_0^*$, 
in contradiction to Proposition \ref{possible configurations}.
So there are no multiple fibers with Kodaira symbol $\I_4^*$.

Suppose next that $\varphi^{-1}(\infty)$ is multiple with Kodaira symbol $\IV^*$.
The dual graph for $\varphi^{-1}(\infty)$ is:
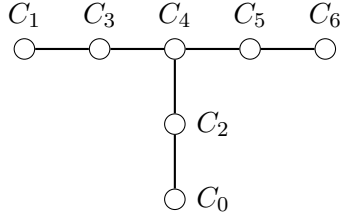
\begin{equation*}
\begin{gathered}
\begin{tikzpicture}
[node distance=1cm, font=\small]
\tikzstyle{vertex}=[circle, draw, fill=white, inner sep=0mm, minimum size=1.5ex]
\node[vertex]	(C1)  	at (1,0) 	[label=above:{$C_1$}] {};
\node[vertex]	(C3)			[right of=C1,  label=above:{$C_3$}]	{};
\node[vertex]	(C4)			[right of=C3,  label=above:{$C_4$}]	{};
\node[vertex]	(C2)			[below of=C4,  label=right:{$C_2$}]	{};
\node[vertex]	(C5)			[right of=C4,  label=above:{$C_5$}]	{};
\node[vertex]	(C6)			[right of=C5,  label=above:{$C_6$}]	{};
 \node[vertex]	(C0)			[below of=C2,  label=right:{$C_0$}]	{};

\draw[thick]    (C1)--(C3)--(C4)-- (C5)--(C6);
\draw[thick]    (C4)--(C2)--(C0);
\end{tikzpicture}
\end{gathered}
\end{equation*}
Now $\varphi^\inv_\ind(\infty)=(C_0+C_1+C_6)+2(C_2+C_3+C_5)+ 3C_4$.
After renumeration  we have $(F\cdot C_0)=1$. In turn, the connected curve $C_1+\ldots+C_6$ has dual graph $T_{2,3,3}$
and must be contained in $\psi^{-1}(\infty)$. 
Now recall that we already established that $\II^*$ does not appear in any genus-one fibration on $Y$.
Consequently, the only possible Kodaira symbols for
the fiber $\psi^{-1}(\infty)$ are $\IV^*$ and $\III^*$. 
In the former case, the fiber is multiple  by Proposition \ref{possible configurations},
and we write  $R\subset\psi^{-1}(\infty)$ for the additional component. 
Then $1=\varphi^{-1}_\ind(\infty)\cdot\psi^{-1}_\ind(\infty) = (C_0+2C_2)\cdot R\geq 2$, contradiction.
So the fiber $\psi^{-1}(\infty)$ has Kodaira symbol $\III^*$.
Let $R,R'\subset\psi^{-1}(\infty)$ be the two additional components, say with
$(R\cdot C_1)=(R'\cdot C_6)=1$. Then
$$
2\geq \varphi^{-1}_\ind(\infty)\cdot\psi^{-1}_\ind(\infty)  = (C_0+C_1+C_6)\cdot (R+R')= (C_0\cdot R) + (C_0\cdot R')+2.
$$
It follows that $\psi^{-1}(\infty)$ is simple. In light of Proposition \ref{possible configurations}, 
the configuration of fibers over the rational points  for $\psi:Y\ra\PP^1$
is $\III^*+E_4+E_4$. In turn, each $(-2)$-curve $C\subset Y$ intersects the fiber $\psi^{-1}(\infty)$.
By Proposition \ref{possible configurations}, 
we may assume that $\varphi^{-1}(0)=2(C_0'+\ldots+C_{r-1}')$ is multiple with  Kodaira symbol $\III$ or $\IV$
with respective integers $r=2$ or $r=3$.
Obviously  $C_0,\ldots,C_6$ are disjoint from  $\varphi^{-1}(0)$.
In turn, $R\cup R'$ intersects each component $C'_i\subset\varphi^\inv(0)$.
Moreover, $R$ is not contained in $\varphi^{-1}(0)$, because 
$R\cdot\varphi^{-1}(0) = R\cdot \varphi^{-1}(\infty)\geq 1$, and likewise for $R'$.
So $2=\varphi^{-1}_\ind(0) \cdot \psi^{-1}(\infty) $ equals
$$
(C_0'+\ldots+C'_{r-1})\cdot \psi^{-1}(\infty) = (C_0'+\ldots+C'_{r-1})\cdot(R+R') \geq r\geq 2.
$$
Thus we have equality at each step, hence $r=2$, and the fiber $\varphi^{-1}(0)$ has Kodaira symbol $\III$.
If $R,R'$ intersect the same component of $\varphi^{-1}(0)$, say $C'_0$, the curve $C_3+\ldots+C_6+R'+C'_0+R+ C_1$ is
a curve of canonical type $\I_8$, in contradiction to Proposition \ref{no I8}.
So after renumeration, we may assume $(R\cdot C'_0)=(R'\cdot C'_1)=1$ and $(R\cdot C'_1)=(R'\cdot C'_0)=0$.
In turn, $R\cap C'_0$ and $C'_1\cap C'_0$ yield two different rational points, and similarly for $C'_1$.
Now consider the   fibers $\psi^{-1}(0)$ and $\psi^{-1}(1)$, which are linearly equivalent curves.
Both are multiple, because $\psi^{-1}(\infty)$ is simple, and their reductions are smooth, as   observed above.
So $\psi^{-1}_\ind(0)$ and $\psi^{-1}_\ind(1)$ are disjoint but numerically equivalent integral curves.
After renumeration, we may assume that they intersect $C'_0$ but not $C'_1$. Then 
$$
R\cap C'_0\quadand C'_1\cap C'_0\quadand \psi^{-1}_\ind(0)\cap C'_0 \quadand \psi^{-1}_\ind(1)\cap C'_0
$$
define four distinct rational points on $C'_0=\PP^1$
contradicting $|\PP^1(\FF_2)|=3$.  Thus there is no fiber with Kodaira symbol $\IV^*$.

Let us now treat that case that  $\varphi^{-1}(\infty)$ is multiple with Kodaira symbol $\IV$.
By  Proposition \ref{possible configurations},  we may apply an automorphism of $\PP^1$ so that 
the  configuration of fibers over the rational points $t=\infty,0,1$ 
becomes $\IV+\I_1^*+E_5$. The former two are multiple, and the latter is simple.
The dual graph for $\varphi^{-1}(\infty)\cup\varphi^{-1}(0)$ is:
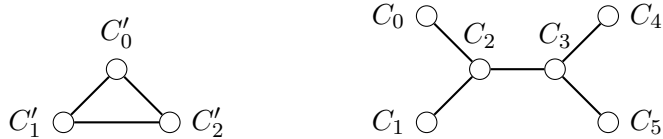
\begin{equation*}
\begin{gathered}
\begin{tikzpicture}
[node distance=1cm, font=\small]
\tikzstyle{vertex}=[circle, draw, fill=white, inner sep=0mm, minimum size=1.5ex]
\node[vertex]	(C2)  	at (1,0) [  label=above:{$C_2$}]{};
\node[vertex]	(C3)	[right of=C2,  label=above:{$C_3$}]{};
\node[vertex]	(C0)	[above left of=C2, label=left:{$C_0$}]{};
\node[vertex]	(C1)	[below left of=C2, label=left:{$C_1$}]{};
\node[vertex]	(C4)	[above right of=C3, label=right:{$C_4$}]{};
\node[vertex]	(C5)	[below right of=C3, label=right:{$C_5$}]{};

\node			(R)		[above left of =   C1, xshift=-1cm, label=below:{}]{};
\node			(R')	[above left of =   R,  xshift=-1cm, label=below:{}]{};

\node[vertex]	(C2')	[below left of = R, xshift=-1cm, label=right:{$C_2'$}]{};
\node[vertex]	(C0')	[above left of = C2', label=above:{$C_0'$}]{};
\node[vertex]	(C1')	[below left of = C0', label=left:{$C_1'$}]{};

\draw [thick] (C4)--(C3)--(C5);
\draw [thick] (C0)--(C2)--(C1);
\draw [thick] (C2)--(C3);

\draw[thick] (C0')--(C2')--(C1')--(C0');

\end{tikzpicture}
\end{gathered}
\end{equation*}
After renumeration, we may assume  that $(F\cdot C_0)=(F\cdot C'_0)=1$.
The   connected curve  $C_1+\ldots+C_5$ has dual graph $T_{2,2,3}$ and must be contained in   the fiber $\psi^{-1}(\infty)$.
Let $m\geq 2 $ be the multiplicity of $C_2\subset \psi^\inv_\ind(\infty)$.
From $2\geq \varphi^{-1}_\ind(\infty)\cdot \psi^{-1}(\infty)\geq C_0\cdot mC_2=m\geq 2$ we infer that the fiber $\psi^\inv(\infty)$
is simple. Likewise, $C'_1$ and $C'_2$ are $\psi$-vertical, say contained in $\psi_\ind^{-1}(0)$, with respective multiplicities $m_1,m_2\geq 1$.
Now
$$
2\geq \varphi^{-1}_\ind(\infty)\cdot \psi^{-1}(0)\geq C'_0\cdot (m_1C'_1+m_2C'_2) =m_1+m_2\geq 2,
$$
we see that  also $\psi^\inv(0)$ is simple, contradiction.
\qed

\medskip
Combining  Propositions  \ref{possible configurations} and \ref{no I8} and  \ref{no II*, IV*, IV}  we get the following key reduction:

\begin{proposition}
\mylabel{key reduction}
Suppose $Y$ is a non-exceptional Enriques surface over $k=\FF_2$ with constant Picard scheme.
Let $\varphi:Y\ra\PP^1$ be a genus-one fibration.
Up to  automorphisms of the projective line,
the possible configurations of fibers over the the rational points $t=0,1,\infty$ are given by
the following table:
$$
\begin{array}[c]{@{}lll}
\toprule
\text{elliptic}		& 		& \text{quasielliptic}\\
\midrule
\I_4^* +E_4+E_2		&		& \I_4^*+\II+\II\\
\midrule
\III^*+E_4+\I_2		&		& \III^*+\III+\II\\ 
\III^*+E_2+\tI_2\\
\III^*+E_4+E_4\\
\midrule
\I_1^*+E_4+\I_4	&		& \I_2^*+\III+\III\\
\bottomrule
\end{array}
$$
Moreover,   all $\I_1^*$-fibers are multiple, the fibers with Kodaira symbol  $\I_4^*,\I_n,\tI_2$ are simple,
and   there is another genus-one fibration
$\psi:Y\ra\PP^1$ with $\varphi^\inv(\infty)\cdot \psi^\inv(\infty)=4$.
\end{proposition}

\section{Elimination of \texorpdfstring{$\I_4^*$}{I4*}-fibers}
\mylabel{Elimination I4*}

Let $Y$ be a non-exceptional Enriques surface over the ground field $k=\FF_2$,  with constant Picard scheme.
The goal of this section is   to exclude the two $\I_4^*$-cases from the table in Proposition \ref{key reduction}.

\begin{proposition}
\mylabel{no I4*}
There is no genus-one fibration 
$\varphi:Y\ra\PP^1$ having a fiber  with Kodaira symbol $\I_4^*$. 
\end{proposition}

\proof
Seeking a contradiction, we assume that $\varphi^{-1}(\infty)$ has Kodaira symbol $\I_4^*$. Its dual graph is:
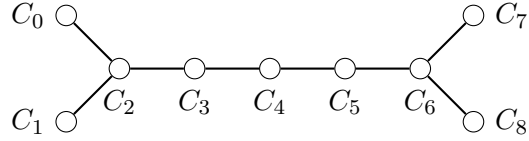
\begin{equation*}
\begin{gathered}
\begin{tikzpicture}
[node distance=1cm, font=\small]
\tikzstyle{vertex}=[circle, draw, fill=white, inner sep=0mm, minimum size=1.5ex]
\node[vertex]	(C2)  	at (1,0) 	[label=below:{$C_2$}] 		{};
\node[vertex]	(C3)			[right of=C2, label=below:{$C_3$}]	{};
\node[vertex]	(C4)			[right of=C3, label=below:{$C_4$}]	{};
\node[vertex]	(C5)			[right of=C4, label=below:{$C_5$}]	{};
\node[vertex]	(C6)			[right of=C5, label=below:{$C_6$}]	{};

\node[vertex]	(C7)			[above right of=C6, label=right:{$C_7$}]	{};
\node[vertex]	(C8)			[below right of=C6, label=right:{$C_8$}]	{};
\node[vertex]	(C0)			[above left  of=C2, label=left :{$C_0$}]	{};
\node[vertex]	(C1)			[below left  of=C2, label=left :{$C_1$}]	{};

\draw [thick] (C2)--(C3)--(C4)--(C5)--(C6);
\draw [thick] (C0)--(C2)--(C1);
\draw [thick] (C7)--(C6)--(C8);
\end{tikzpicture}
\end{gathered}
\end{equation*}
The fiber must be simple, according to   Proposition \ref{key reduction}, 
so the schematic fiber is  $\varphi^{-1}(\infty) = C_0+C_1+2(C_2+\ldots+C_6)+C_7+C_8$.
Moreover,  there is a genus-one fibration $\psi:Y\ra\PP^1$ so that the multiple fiber $2F$ have  $F\cdot\varphi^{-1}(\infty)=2$.
After renumeration, we have the following   possibilities: $(F\cdot C_i)=1$ with $i\in\{2,3,4\}$
or $(F\cdot C_0)=(F\cdot C_1)=1$ or $(F\cdot C_0)=(F\cdot C_8)=1$ or $(F\cdot C_0)=2$.
Our task is to rule out each of these six possibilities. 
Without restriction, we assume that   $\psi^{-1}(\infty)$ contains the largest number of irreducible components
among all fibers.

Suppose first that $(F\cdot C_4)=1$. Then the connected curves $C_0+\ldots+C_3$ and $C_5+\ldots+C_8$ both have
dual graph $T_{2,2,2}$ and are contained in some fibers of $\psi$. According to Proposition \ref{key reduction}, 
they both must be contained in  
$\psi^{-1}(\infty)$, which therefore has  Kodaira symbol $\I_4^*$ and must be simple.
Let $R\subset \psi^{-1}(\infty)$ the additional component, which has $(R\cdot C_i)=(R\cdot C_j)=1$ for some
$i\in\{0,1,3\}$ and $j\in\{5,7,8\}$. Let $1\leq m_i,m_j\leq 2$ be the multiplicities of $C_i,C_j\subset\varphi^{-1}(\infty)$.
From
$$
4=\varphi^{-1}(\infty)\cdot\psi^{-1}(\infty) = (m_iC_i+2C_4+m_jC_j)\cdot 2R = 2m_i+4(C_4\cdot R) + 2m_j 
$$
we infer that $R\cap C_4=\varnothing$ and $m_i=m_j=1$. After renumeration,
we may assume $i=1$ and $j=7$. Then $C_1+\ldots+C_7+R$ is a curve
of canonical type $\I_8$, in contradiction to Proposition \ref{key reduction}.

Suppose next that $(F\cdot C_3)=1$.  
The connected curve $C_4+\ldots+C_8$ has dual graph $T_{2,2,3}$ and must be contained in   $\psi^{-1}(\infty)$.
The curve $C_0+C_2+C_1$ is a chain and must be part of some fiber $\psi^{-1}(t)$ as well.
If $t\neq\infty$, 
we infer from Proposition \ref{key reduction} that $\psi^{-1}(t)$ has Kodaira symbol $\I_4$ 
and is simple.  Let $R\subset \psi^{-1}(t)$ be the additional component, with $(R\cdot C_0)=(R\cdot C_1)=1$.
This gives 
$$
4=\varphi^{-1}(\infty)\cdot \psi^{-1}(t) = (C_0+C_1+2C_3) \cdot R = 1+1+ 2(C_3\cdot R),
$$
and hence  $(C_3\cdot R)=1$. Thus $C_1+C_2+C_3+R$ is another curve of canonical type $\I_4$. Its intersection number with $C_4$ is one,
so the fiber is multiple, contradiction.
Summing up,  $C_0+C_1+C_2$ and $C_4+\ldots+C_8$ both belong to   $\psi^{-1}(\infty)$, which therefore
has Kodaira symbol $\I_4^*$ and must be simple, in light of Proposition \ref{key reduction}. Let $R$ be the additional component, which necessarily has
$(R\cdot C_2)=(R\cdot C_4)=1$. This gives 
$$
4=\varphi^{-1}(\infty)\cdot \psi^{-1}(\infty) = (2C_2+2C_3+2C_4)\cdot 2R = 4+4(C_3\cdot R) + 4\geq 8,
$$
contradiction.

Now suppose that $(F\cdot C_2)=1$. Then the connected curve $C_3+\ldots + C_8$ has dual graph $T_{2,2,4}$ 
and must be contained in $\psi^{-1}(\infty)$. Its Kodaira symbol is either $\II^*$, $\I_4^*$, $\III^*$ or $\I_2^*$.
The first alternative  is impossible by Proposition \ref{key reduction}. 
In the second alternative, the fiber must be simple.
The additional component $R\subset\psi^{-1}(\infty)$ has $(R\cdot C_i)=1$ for all $i\in\{0,1,3\}$,
which gives the contradiction
$$
4=\varphi^{-1}(\infty)\cdot\psi^{-1}(\infty) = (C_0+C_1+2C_2+2C_3)\cdot 2R= 2+2+4(C_2\cdot R) + 4\geq 8.
$$
In the third  alternative, we deduce
from Proposition \ref{key reduction} that one of the curves from $C_0+C_1$ also belongs to $\psi^{-1}(\infty)$.
Without restriction, we may assume that this is $C_1$. Let $R\subset\psi^{-1}(\infty)$ be the additional component,
say with  $(R\cdot C_7)=(R\cdot C_1)=1$. From
$$
4\geq \varphi^{-1}(\infty)\cdot\psi^{-1}_\ind(\infty) = (C_1+2C_2+C_7)\cdot 2R = 2+ 4(C_2\cdot R) + 2
$$
we get $R\cap  C_2=\varnothing$. But then $C_1+\ldots+C_7+R$ is a curve of canonical type $\I_8$, in contradiction
to Proposition \ref{key reduction}.
The only remaining possibility is that  $\psi^{-1}(\infty) $ has Kodaira symbol $\I_2^*$.
So the additional component $R\subset\psi^{-1}(\infty)$ has $(R\cdot C_4)=1$. From
$$
4\geq \varphi^{-1}(\infty)\cdot\psi^{-1}_\ind(\infty) = (2C_2+2C_4)\cdot R=   2(C_2\cdot R) + 2
$$
we deduce that $(C_2\cdot R)$ is either zero or one.
In the latter case we get a curve $D=C_2+C_3+C_4+R$ of canonical type $\I_4$.
The ensuing fiber is multiple, because $(D\cdot C_5)=1$, in contradiction to Proposition \ref{key reduction}.
The former case yields the following dual graph, in contradiction to  Proposition \ref{no I4*+R} below.
\begin{equation*}
\begin{gathered}
\begin{tikzpicture}
[node distance=1cm, font=\small]

\tikzstyle{vertex}=[circle, draw, fill=white, inner sep=0mm, minimum size=1.5ex]
\node[vertex]	(C2)  	at (1,0) 	[label=below:{$C_2$}] 		{};
\node[vertex]	(C3)			[right of=C2, label=below:{$C_3$}]	{};
\node[vertex]	(C4)			[right of=C3, label=below:{$C_4$}]	{};
\node[vertex]	(C5)			[right of=C4, label=below:{$C_5$}]	{};
\node[vertex]	(C6)			[right of=C5, label=below:{$C_6$}]	{};

\node[vertex]	(C7)			[above right of=C6, label=right:{$C_7$}]	{};
\node[vertex]	(C8)			[below right of=C6, label=right:{$C_8$}]	{};
\node[vertex]	(C0)			[above left  of=C2, label=left :{$C_0$}]	{};
\node[vertex]	(C1)			[below left  of=C2, label=left :{$C_1$}]	{};
\node[vertex]	(R)			[above of=C4, label=right :{$R$}]	{};

\draw [thick] (C2)--(C3)--(C4)--(C5)--(C6);
\draw [thick] (C0)--(C2)--(C1);
\draw [thick] (C7)--(C6)--(C8);
\draw [thick] (C4)--(R);
\end{tikzpicture}
\end{gathered}
\end{equation*}
Next, suppose that we are in the case $(F\cdot C_0)=(F\cdot C_1)=1$. Then the connected curve $C_2+\ldots+C_8$
has dual graph $T_{2,2,5}$  and must be contained in $\psi^{-1}(\infty)$.
The possible Kodaira symbols are $\II^*$ and $\I_4^*$. 
The former is impossible by Proposition \ref{key reduction}, so $\psi^\inv(\infty)$ is simple with Kodaira symbol $\I_4^*$.   
Let $R,R'\subset\psi^{-1}(\infty)$ be the two additional components,
with $(R\cdot C_2)=(R'\cdot C_2)=1$. Then $4=\varphi^{-1}(\infty)\cdot \psi^{-1}(\infty)$ becomes
$$
(C_0+C_1+2C_2)\cdot (R+R') = (C_0+C_1)\cdot (R+R') + 2+2, 
$$
and we infer that $R\cup R'$ is disjoint from $C_0\cup C_1$. In turn, $C_0+C_1+R+R'+C_2$ supports a curve of canonical type
$\I_0^*$, in  contradiction to Proposition \ref{key reduction}.

Now suppose $(F\cdot C_0)=(F\cdot C_8)=1$. Then $C_1+\ldots +C_7$ is a chain and contained in $\psi^{-1}(\infty)$.
The possible Kodaira symbols are $\I_4^*$ or $\III^*$, by Proposition \ref{key reduction}. 
The symbol $\I_4^*$ is handled as in the preceding paragraph. In the remaining case $\III^*$, let $R\subset \psi^{-1}(\infty)$ be the
additional component, which has $(R\cdot C_4)=1$. From
$$
4\geq \varphi^{-1}(\infty)\cdot \psi^{-1}_\ind(\infty) = (C_0+2C_4+C_8)\cdot 2R = 2(C_0\cdot R) + 4 + 2(C_8\cdot R)
$$
we infer that $R$ is disjoint from $C_0\cup C_8$. Now the configuration $C_0+\ldots+C_8+R$ contradicts Proposition \ref{no I4*+R} below.

It remains to cope with the last possibility   $(F\cdot C_0)=2$, which is the most difficult. 
Now the connected curve $C_1+\ldots+C_8$ has dual graph $T_{2,2,6}$
and must be contained in $\psi^{-1}(\infty)$. Its Kodaira symbol has to be $\I_4^*$, so the fiber is simple. The additional
component $R\subset\psi^{-1}(\infty)$ has $(R\cdot C_2)=1$, and from
$$
4=\varphi^\inv(\infty)\cdot \psi^\inv(\infty) = (C_0+2C_2)\cdot R = (C_0\cdot R) + 2
$$
we infer that the dual graph for $C_0+\ldots+C_8+R$ must be:
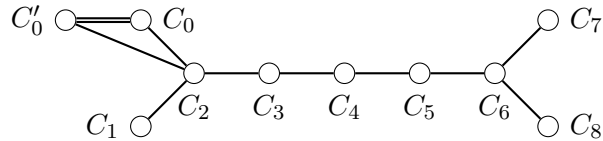
\begin{equation*}
\begin{gathered}
\begin{tikzpicture}
[node distance=1cm, font=\small]
\tikzstyle{vertex}=[circle, draw, fill=white, inner sep=0mm, minimum size=1.5ex]
\node[vertex]	(C2)  	at (1,0) 	[label=below:{$C_2$}] 		{};
\node[vertex]	(C3)			[right of=C2, label=below:{$C_3$}]	{};
\node[vertex]	(C4)			[right of=C3, label=below:{$C_4$}]	{};
\node[vertex]	(C5)			[right of=C4, label=below:{$C_5$}]	{};
\node[vertex]	(C6)			[right of=C5, label=below:{$C_6$}]	{};

\node[vertex]	(C7)			[above right of=C6, label=right:{$C_7$}]	{};
\node[vertex]	(C8)			[below right of=C6, label=right:{$C_8$}]	{};
\node[vertex]	(C0)			[above left  of=C2, label=right :{$C_0$}]	{};
\node[vertex]	(C1)			[below left  of=C2, label=left :{$C_1$}]	{};
\node[vertex]	(R)				[left of=C0, label=left :{$R$}]	{};

\draw [thick] (C2)--(C3)--(C4)--(C5)--(C6);
\draw [thick] (C0)--(C2)--(C1);
\draw [thick] (C7)--(C6)--(C8);
\draw [thick, double] (C0)--(R);
\draw [thick] (C2)--(R);
\end{tikzpicture}
\end{gathered}
\end{equation*}
But this configuration is impossible, by Proposition \ref{no I4*+R second} below.
\qed

\medskip
In the above arguments, we have used the following facts:

\begin{proposition}
\mylabel{no I4*+R}
There is no configuration $C_0+\ldots+C_8+R$ of ten $(-2)$-curves on $Y$ with simple normal crossings and the following dual graph:
\begin{equation}
\label{I4*+R}
\begin{gathered}
\begin{tikzpicture}
[node distance=1cm, font=\small]

\tikzstyle{vertex}=[circle, draw, fill=white, inner sep=0mm, minimum size=1.5ex]
\node[vertex]	(C2)  	at (1,0) 	[label=below:{$C_2$}] 		{};
\node[vertex]	(C3)			[right of=C2, label=below:{$C_3$}]	{};
\node[vertex]	(C4)			[right of=C3, label=below:{$C_4$}]	{};
\node[vertex]	(C5)			[right of=C4, label=below:{$C_5$}]	{};
\node[vertex]	(C6)			[right of=C5, label=below:{$C_6$}]	{};

\node[vertex]	(C7)			[above right of=C6, label=right:{$C_7$}]	{};
\node[vertex]	(C8)			[below right of=C6, label=right:{$C_8$}]	{};
\node[vertex]	(C0)			[above left  of=C2, label=left :{$C_0$}]	{};
\node[vertex]	(C1)			[below left  of=C2, label=left :{$C_1$}]	{};
\node[vertex]	(R)			[above of=C4, label=right :{$R$}]	{};

\draw [thick] (C2)--(C3)--(C4)--(C5)--(C6);
\draw [thick] (C0)--(C2)--(C1);
\draw [thick] (C7)--(C6)--(C8);
\draw [thick] (C4)--(R);
\end{tikzpicture}
\end{gathered}
\end{equation}
\end{proposition}

\proof
Seeking a contradiction, we assume that the configuration  exists. 
The subconfiguration  $C_0+\ldots+C_8$ supports a curve of canonical type $\I_4^*$. Let $\varphi:Y\ra\PP^1$ be the resulting
genus-one fibration.  Without restriction, we may assume that 
$\varphi^{-1}(\infty) = C_0+C_1+2(C_2+\ldots+C_6)+C_7+C_8$.
The fiber   indeed must be simple by Proposition \ref{key reduction}. Write $F\subset Y$ for some half-fiber, and let 
$S\subset\Num(Y)$ be the subgroup generated by $C_1,\ldots,C_8,R,F$.
To proceed, we exploit properties of the resulting Gram matrix $N\in\Mat_{10}(\ZZ)$, which can easily be checked
with computer algebra.
One computes in this way $\det(N)=-4$, whence $S\subset \Num(Y)$ if a subgroup of index two.
%
%
%
%

Nikulin's theory of discriminant forms \cite{Nikulin 1980} allows us to use the lattice $S$ to gain enough control of the slightly larger  $\Num(Y)$
and thus the geometry of  $Y$.
Set $S^*=\Hom(S,\ZZ)$, and consider the   injection $ S\ra S^*$ given by $x\mapsto (y\mapsto (x\cdot y))$.
This linear map is described, with respect to the given basis $C_1,\ldots,F\in S$ and the   dual basis $C_1^*,\ldots,F^*\in S^*$, 
by the matrix    $N\in\Mat_n(\ZZ)$.
With respect to this dual basis, the  induced $\QQ$-valued form   on $S^*$ has Gram matrix 
\begin{equation}
\label{gram matrix}
N^{-1}=
{\tiny
\begin{pmatrix} 
  -1   &-1   &-1   &-1   &-1   &-1 &-1/2 &-1/2   & 0   & 1\\
  -1   &-2   &-2   &-2   &-2   &-2  & -1 &  -1   & 0   & 2\\
  -1   &-2   &-3   &-3   &-3   &-3 &-3/2 &-3/2   & 0   & 3\\
  -1   &-2   &-3   &-4  & -4   &-4 &  -2 &  -2   & 0   & 4\\
  -1   &-2   &-3   &-4   &-5   &-5 &-5/2 &-5/2   & 0   & 4\\
  -1   &-2   &-3   &-4   &-5   &-6 &  -3 &  -3   & 0   & 4\\
-1/2   &-1 &-3/2   &-2 &-5/2   &-3 &  -2 &-3/2   & 0   & 2\\
-1/2   &-1 &-3/2   &-2 &-5/2   &-3 &-3/2 &  -2   & 0   & 2\\
   0   & 0   & 0   & 0  &  0   & 0 &   0 &   0   & 0   & 1\\
   1   & 2   & 3   & 4  & 4    &4  &  2  &  2    & 1   &-2
\end{pmatrix}
}
\end{equation}
and we have $S\subset\Num(Y)\subset S^*$.  Note that for each  column, the entries are the coordinates of the corresponding
dual basis vectors with respect to the basis vectors.
%
%
The  invariant factors  of $N$ are  computed to be $(2,2)$, so the \emph{discriminant group} $A_S=S^*/S$ is elementary abelian of order four,
hence a 2-dimensional vector space over $\FF_2$.
%
%
As explained in loc.\ cit., Section 3,
it comes with a  quadratic form $q:A_S\ra\QQ/2\ZZ$ whose associated bilinear form $A_S\times A_S\ra\QQ/\ZZ$ is non-degenerate.
Moreover, the generator for the extension $S\subset\Num(Y)$   corresponds  to a non-zero isotropic vector in $A_S$, by loc.\ cit., Section 4. 
Examining the columns in  \eqref{gram matrix},  we see that 
$$
C_1^*\equiv C_3^*\equiv C_5^*\quadand C_7^*\quadand C_8^* 
$$
are non-zero modulo $S$, whereas the other dual basis vectors vanish modulo  $S$. Moreover, one sees that $C_7^*+C_8^*\equiv C^*_1$ modulo $S$.
So each two of the above  form a basis of $A_S=S^*/S$, and the third one is their sum. Moreover,
$$
(C_1^*)^2=-1,\quad (C_7^*)^2=(C_8^*)^2=-2,\quad (C_1^*\cdot C_7^*)=(C_1^*\cdot C_8^*)=-1/2,
$$
hence the isotropic vectors in the discriminant group come from $C_7^*$ and $C_8^*$. 
 
\newcommand{\tR}{\tilde{R}}
Now back to our genus-one fibration $\varphi:Y\ra\PP^1$.
Consider the    genus-one curve $Y_K=\varphi^\inv(\eta)$ over the function field  $K=k(\PP^1)$.
The restriction map  $\Pic(Y)\ra\Pic(Y_K)$  is surjective. It factors over $\Num(Y)$  because
$\omega_Y$ can be written as the difference of half-fibers. 
The curve $Y_K$ contains no rational points,   the two-section $R$ gives a splitting $\Pic(Y_K)=\Pic^0(Y_K)\oplus 2\ZZ$,
and the curves $C_0,\ldots,C_8,F$ vanish in $\Pic(Y_K)$. This gives an identification $\Pic^0(Y_K)=\Num(Y)/S$, so this group has order two.
Write $\shN_K$ for the non-trivial invertible sheaf of degree zero on $Y_K$.
By Riemann--Roch we have $\shN_K(R)\simeq\O_{Y_K}(\tR)$ for another  two-section $\tR\subset Y$, and thus  $\Num(Y)=S+\ZZ\tR$.
Note that such a  two-section is not unique;
according to the Enriques Reducibility Theorem (\cite{Lang 1983}, Theorem 4.1), we may furthermore assume that $\tR^2=0$ or $\tR^2=-2$.
In light of the preceding paragraph, the element $\tR\in S^*$
is either congruent  to $C_7^*$ or $C_8^*$, modulo $S$.
Without loss of generality $\tR\equiv C_8^*$, which ensures $(\tR\cdot C_8)=1$.
Then $(\tR\cdot C_i)=0$ for $2\leq i\leq 6$, and $(\tR\cdot C_j)=1$ for exactly one $j\in\{0,1,7\}$ besides $j=8$. 
If $(\tR\cdot C_7)=1$ then $\tR\equiv C_7^*+C^*_8\equiv C_1^*$, contradiction.
So without loss of generality we may assume  $(\tR\cdot C_0)=1$. 
If $\tR^2=-2$ then $\tR=\PP^1$, and  $C_0+C_2+\ldots+C_6+C_8+\tR$ is of canonical type $\I_8$, which is impossible by Proposition \ref{key reduction}.
Thus $\tR^2=0$, and $\tR$ is a genus-one curve.

We can easily compute the intersection number    $n=(\tR\cdot R)$: As member of $S^*$ that is perpendicular to $C_1,\ldots,C_7$ we   have
$\tR = C_8^*+ nR^* + F^*$.  Glancing at  \eqref{gram matrix} again, one sees  
$$
(C_8^*)^2=(F^*)^2=-2,\quad (R^*)^2= (C_8^*\cdot R^*) = 0,\quad (C_8^*\cdot F^*) = 2,\quad (R^*\cdot F^*) = 1.
$$
From
$0=\tR^2=(C_8^*+ nR^* + F^*)^2 = (-2 + 0 -2) + 2(0+2+n)$ we get $n=0$, thus the curves $R,\tR\subset Y$ must be disjoint.
 
Let $\psi:Y\ra\PP^1$ be the   genus-one fibration with $\psi^{-1}(0)=2\tR$.  The fiber is indeed multiple,
because $(\tR\cdot C_0)=1$.  The curves $C_1+\ldots+C_7+R$ support a curve of canonical type $\III^*$ disjoint from $\tR$,
and we can assume that $\psi^{-1}(\infty)$ is the corresponding fiber. The latter is simple, because
$(C_0\cdot \tR)=1$ and $(C_0\cdot\psi^{-1}_\ind(\infty)) = (C_0\cdot 2C_2)=2$. The dual graph \eqref{I4*+R} gives 
$\psi^{-1}(\infty)\cdot \varphi^{-1}(\infty)=4$.
According to Proposition \ref{key reduction}, the configuration of fibers   for $\varphi$ is either $\I_4^*+E_4+E_2$ or $\I_4^*+\II+\II$.
In both cases we find some half fiber $\tilde{F}$ for $\varphi$ that is integral and whose  regular locus contains exactly two    rational points.
These rational points must be the intersections with the two half-fibers $\psi^{-1}_\red(0)$ and $\psi^\inv_\red(1)$.
However, we have $2=(\tilde{F}\cdot \psi^{-1}(\infty))=(\tilde{F}\cdot 2R)$, thus $\tilde{F}\cap R$ must be a third rational point in the regular locus of $\tilde{F}$,
contradiction.
\qed

\medskip
In a similar way, we establish:

\begin{proposition}
\mylabel{no I4*+R second}
There is no configuration $C_0+\ldots+C_8+R$ of ten $(-2)$-curves 
having    simple normal crossings, with   exception  $(C_0\cdot R)=2$, and the following dual graph:
\begin{equation}
\label{I4*+R second}
\begin{gathered}
\begin{tikzpicture}
[node distance=1cm, font=\small]
\tikzstyle{vertex}=[circle, draw, fill=white, inner sep=0mm, minimum size=1.5ex]
\node[vertex]	(C2)  	at (1,0) 	[label=below:{$C_2$}] 		{};
\node[vertex]	(C3)			[right of=C2, label=below:{$C_3$}]	{};
\node[vertex]	(C4)			[right of=C3, label=below:{$C_4$}]	{};
\node[vertex]	(C5)			[right of=C4, label=below:{$C_5$}]	{};
\node[vertex]	(C6)			[right of=C5, label=below:{$C_6$}]	{};

\node[vertex]	(C7)			[above right of=C6, label=right:{$C_7$}]	{};
\node[vertex]	(C8)			[below right of=C6, label=right:{$C_8$}]	{};
\node[vertex]	(C0)			[above left  of=C2, label=right :{$C_0$}]	{};
\node[vertex]	(C1)			[below left  of=C2, label=left :{$C_1$}]	{};
\node[vertex]	(R)			[left of=C0, label=left:{$R$}]	{};

\draw [thick] (C2)--(C3)--(C4)--(C5)--(C6);
\draw [thick] (C0)--(C2)--(C1);
\draw [thick] (C7)--(C6)--(C8);
\draw [thick,double] (C0)--(R);
\draw [thick] (C2)--(R);
\end{tikzpicture}
\end{gathered}
\end{equation}
\end{proposition}

\proof
Seeking a contradiction, we assume that such a configuration exists. 
The subconfigurations  $C_0+C_1+\ldots+C_8$ and $R+C_1+\ldots+C_8$ support
curves of canonical type $\I_4^*$. Let $\varphi:Y\ra \PP^1$ and $\varphi':Y\ra\PP^1$
be the resulting genus-one fibrations. Without loss of generality we may assume
\begin{gather*}
\varphi^\inv(\infty) = C_0+C_1+2(C_2+\ldots+C_6)+C_7+C_8,\\
\varphi'^\inv(\infty) = R+C_1+2(C_2+\ldots+C_6)+C_7+C_8.
\end{gather*}
Indeed, these fibers must be simple by Proposition \ref{key reduction}, and   $\varphi^\inv(\infty)\cdot \varphi'^\inv(\infty)=4$.
Choose some half-fibers $F,F'\subset Y$ for the respective fibrations, with intersection number  $(F\cdot F')=1$.
Let $S\subset \Num(Y)$ be the subgroup generated by $C_1,\ldots,C_8,F',F$. 
As in the preceding proof, we exploit   properties of the resulting Gram matrix $N\in\Mat_{10}(\ZZ)$.
Note that this is the lattice $U\oplus D_8$, as one referee pointed out. 
%
%
%
%
%
Again    $S\subset \Num(Y)$ has index two, and the discriminant   group $A_S=S^*/S$  is elementary abelian of order four.
Recall that the inverse
\begin{equation}
\label{gram matrix second}
N^{-1}=
{\tiny
\begin{pmatrix}  
-1   &-1   &-1   &-1   &-1   &-1 &-1/2 &-1/2    &0    &0\\
-1   &-2   &-2   &-2   &-2   &-2   &-1   &-1    &0    &0\\
-1   &-2   &-3   &-3   &-3   &-3 &-3/2 &-3/2    &0    &0\\
-1   &-2   &-3   &-4   &-4   &-4   &-2   &-2    &0    &0\\
-1   &-2   &-3   &-4   &-5   &-5 &-5/2 &-5/2    &0    &0\\
-1   &-2   &-3   &-4   &-5   &-6   &-3   &-3    &0    &0\\
-1/2 &-1 &-3/2   &-2 &-5/2   &-3   &-2 &-3/2    &0    &0\\
-1/2 &-1 &-3/2   &-2 &-5/2   &-3 &-3/2   &-2    &0    &0\\
0    &0    &0    &0    &0    &0    &0    &0     &0    &1\\
0    &0    &0    &0    &0    &0    &0    &0     &1    &0\\
\end{pmatrix}
}
\end{equation}
is the Gram matrix of the induced $\QQ$-valued form on $S^*$ with respect to the dual basis $C_1^*,\ldots C_8^*,F'^*, F^*$.
The generator for the extension $S\subset\Num(Y)$   corresponds  to a non-zero isotropic vector in $A_S$, with
respect to the quadratic form $q:A_R\ra\QQ/2\ZZ$.  
Glancing at \eqref{gram matrix second}, we see that 
$$
C_1^*\equiv C_3^*\equiv C_5^*\quadand C_7^*\quadand C_8^*. 
$$
are non-zero modulo $S$,     whereas the other dual basis vectors vanish modulo  $S$, and that
$(C_1^*)^2=-1$ and $(C_7^*)^2=(C_8^*)^2=-2$. So the isotropic vectors in the discriminant group come from $C_7^*$ and $C_8^*$.
Note that   $F^*=F'$ and $F'^*=F$ belong to the sublattice $S\subset S^*$.

By Proposition \ref{key reduction}, the half fiber  $F'$ is   integral. 
We have $F'\cdot \varphi^\inv(\infty) =2$, hence $F'$ is a two-section for  $\varphi:Y\ra\PP^1$.
Let  $Y_K=\varphi^{-1}(\eta)$ be the generic fiber.
As in the preceding proof, the two-section $F'$ gives a splitting $\Pic(Y_K)=\Pic^0(Y_K)\oplus 2\ZZ$,
we get an identification $\Pic(Y_K)=\Num(Y)/S$,
and there must be another integral curve $\tilde{R}\subset Y$ that is a two-section, 
with $\Num(Y)=S+\ZZ\tR$. Furthermore, we may may assume $\tR^2=0$ or $\tR^2=-2$.
We now examine the possible non-zero intersection numbers with the components of $\varphi^\inv(\infty)$:
Since the discriminant group is annihilated by two, 
the case $(\tR\cdot C_i)=2$ with $i\in \{0,1,7,8\}$ would imply that $\tR\in S^*$ belongs to $S$, contradiction.
Also $(\tR\cdot C_i)=1$ with $i\in \{2,\ldots,6\}$ is impossible, because then $\tR$ becomes zero or non-isotropic in $A_S=S^*/S$. 
Likewise, the two cases 
$$
 (\tR\cdot C_0) = (\tR\cdot C_1) =1\quadand  (\tR\cdot C_7)=(\tR\cdot C_8)=1
$$
are impossible, because $C_0^*+C_1^*\equiv F^*+C_1^*$ and $C_7^*+C_8^*$ are non-isotropic.
Up to renumeration, the only  possibilities  left are $(\tR\cdot C_0) = (\tR\cdot C_8)=1$ and $(\tR\cdot C_1) = (\tR\cdot C_8)=1$.
If $\tR^2=-2$ then $\tR=\PP^1$ and in both cases get  a curve of canonical type $\I_8$, in contradiction to Proposition \ref{key reduction}.
Thus $\tR^2=0$, and $\tR\subset Y$ is a genus-one curve.  But then the two cases where already discarded in the proof for Proposition \ref{no I4*}.
Summing up, the configuration \eqref{I4*+R second} also does not exist.
\qed

\section{Elimination of    \texorpdfstring{$\III^*$}{III*}-fibers}
\mylabel{Elimination III*}

Let $Y$ is a non-exceptional Enriques surface over the ground field $k=\FF_2$, with constant Picard scheme.
We now exclude the $\III^*$-cases from the table in Proposition \ref{key reduction}.
We proceed in two steps, treating the case of multiple fibers first:

\begin{proposition}
\mylabel{no multiple III*}
There is no genus-one fibration  $\varphi:Y\ra\PP^1$ having a multiple  fiber with Kodaira symbol $\III^*$. 
\end{proposition}

\proof
Seeking a contradiction, we assume   that $\varphi^{-1}(\infty)$
is multiple with Kodaira symbol $\III^*$. The dual graph is:
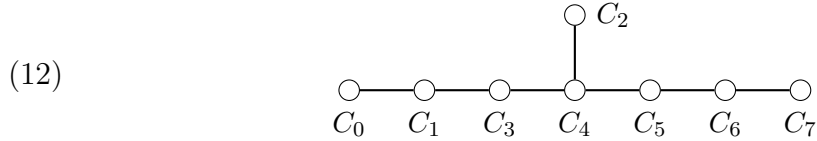
\begin{equation}
\label{dual graph III*}
\begin{gathered}
\begin{tikzpicture}
[node distance=1cm, font=\small]
\tikzstyle{vertex}=[circle, draw, fill=white, inner sep=0mm, minimum size=1.5ex]
\node[vertex]	(C0)  	at (1,0) [label=below:{$C_0$}] {};
\node[vertex]	(C1)	[right of=C0,  label=below:{$C_1$}]	{};
\node[vertex]	(C3)	[right of=C1,  label=below:{$C_3$}]	{};
\node[vertex]	(C4)	[right of=C3,  label=below:{$C_4$}]	{};
\node[vertex]	(C5)	[right of=C4,  label=below:{$C_5$}]	{};
\node[vertex]	(C6)	[right of=C5,  label=below:{$C_6$}]	{};
\node[vertex]	(C7)	[right of=C6,  label=below:{$C_7$}]	{};
\node[vertex]	(C2)	[above of=C4,  label=right:{$C_2$}]	{};
 
\draw[thick]    (C0)--(C1)--(C3)--(C4)-- (C5)--(C6)--(C7);
\draw[thick]    (C2)--(C4);
\end{tikzpicture}
\end{gathered}
\end{equation}
According to Proposition \ref{key reduction}, there is a genus-one fibration $\psi:Y\ra\PP^1$ having $\psi^{-1}(\infty)\cdot\varphi^{-1}(\infty)=4$.
We may assume that  $\psi^{-1}(\infty)$ has the largest number of irreducible components, among all fibers.
Choose a multiple fiber  $2F$  for $\psi$. After renumeration, we get $(F\cdot C_0)=1$.  
The connected curve $C_1+\ldots+C_7$ has dual graph $T_{2,3,4}$ and is thus contained in $\psi^{-1}(\infty)$. Its Kodaira
symbol must be $\III^*$, according to Proposition \ref{key reduction}. Let $C_0'\subset\psi^{-1}(\infty)$ be 
the additional component, which has $(C_0'\cdot C_1)=1$. From
$$
2\geq \varphi^{-1}_\ind(\infty)\cdot \psi_\ind(\infty) = (C_0+2C_1)\cdot C'_0 = (C_0\cdot C'_0) + 2
$$
we infer that $C_0$ and $C'_0$ are disjoint.
Thus  $C_0'+C_0+\ldots+ C_5$ supports a curve of canonical type $\I_2^*$.
Let $\tau:Y\ra\PP^1$ be the resulting genus-one fibration that has this curve as reduced fiber over $t=\infty$.
The fiber is multiple, because $\tau^{-1}_\ind(\infty)\cdot C_6=1$.
According to Proposition \ref{key reduction}, 
the fibers over $t=0,1$ have Kodaira symbol $\III$. Clearly, the curve $C_7$ belongs to one such fiber.
After applying an automorphism of $\PP^1$ we obtain
$$
\tau^{-1}_\ind(0)= C_7+C_7'\quadand \tau^{-1}_\ind(1) = R+R'.
$$
The dual graph for the configuration $C_0'+C_0+\ldots+C_7+C_7'+R+R'$ takes the following shape:
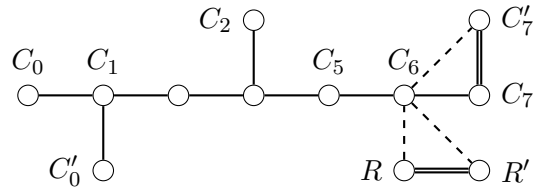
\begin{equation*}
\begin{gathered}
\begin{tikzpicture}
[node distance=1cm, font=\small]
\tikzstyle{vertex}=[circle, draw, fill=white, inner sep=0mm, minimum size=1.5ex]
\node[vertex]	(C0)  	at (1,0) 	[label=above:{$C_0$}] {};
\node[vertex]	(C1)			[right of=C0,  label=above:{$C_1$}]	{};
\node[vertex]	(C3)			[right of=C1,  label=above:{$$}]	{};
\node[vertex]	(C4)			[right of=C3,  label=above:{$$}]	{};
\node[vertex]	(C5)			[right of=C4,  label=above:{$C_5$}]	{};
\node[vertex]	(C6)			[right of=C5,  label=above:{$C_6$}]	{};
\node[vertex]	(C7)			[right of=C6,  label=right:{$C_7$}]	{};
\node[vertex]	(C2)			[above of=C4,  label=left:{$C_2$}]	{};
 
\node[vertex]	(C0')			[below of=C1,  label=left:{$C_0'$}]	{};
\node[vertex]	(C7')			[above of=C7,  label=right:{$C_7'$}]	{};
\node[vertex]	(R)			[below of=C6,  label=left:{$R$}]	{};
\node[vertex]	(R')			[below of=C7,  label=right:{$R'$}]	{};
 
\draw[thick]    (C0)--(C1)--(C3)--(C4)-- (C5)--(C6)--(C7);
\draw[thick]    (C2)--(C4);
\draw[thick]    (C1)--(C0');
\draw[thick,double]    (C7)--(C7');
\draw[thick,double]    (R)--(R');
\draw[thick,dashed]    (R)--(C6)--(R');
\draw[thick,dashed]     (C6)--(C7');
\end{tikzpicture}
\end{gathered}
\end{equation*}
The dashed edges indicate intersection numbers on which we can say the following:
First note that the intersection number $(R\cdot C_6)=1$   would lead to the configuration $C_0'+C_0+\ldots+C_7+R$, 
which contradicts   Proposition \ref{no I4*+R}.
For the same reason $(R'\cdot C_6)=1$  is impossible. 
On the other hand, we have 
$$
2= C_6\cdot\tau^{-1}(\infty) \geq C_6\cdot \tau^{-1}_\ind(1)=(C_6\cdot R)+(C_6\cdot R').
$$
After interchanging $R$ and $R'$ if necessary, we may assume that $(C_6\cdot R')=2$ and $(C_6\cdot R)=0$.
If follows that $R$ belongs to the fibers of  $\psi:Y\ra\PP^1$. 
According to Proposition \ref{possible configurations},   the    fiber $\psi^{-1}(\infty)$ 
with Kodaira symbol $\III^*$ must be multiple. 
Hence
$$
1=\varphi^{-1}_\ind(\infty)\cdot \psi^{-1}_\ind(\infty) = (2C_1\cdot C_0') =2,
$$
contradiction.
\qed

\begin{proposition}
\mylabel{no simple III*}
There is no genus-one fibration $\varphi:Y\ra\PP^1$ having a simple fiber with Kodaira symbol $\III^*$. 
\end{proposition}

\proof
Seeking a contradiction, we assume  that $\varphi^{-1}(\infty)$
is simple with Kodaira symbol $\III^*$, with dual graph as in \eqref{dual graph III*}.
According to Proposition \ref{key reduction}, 
there is another genus-one fibration $\psi:Y\ra\PP^1$ having $\psi^{-1}(\infty)\cdot\varphi^{-1}(\infty)=4$.
We may assume that  $\psi^{-1}(\infty)$ has the largest number of irreducible components, among all fibers.
Let $2F$ be one of the multiple fibers for $\psi$, such that $F\cdot \varphi^\inv(\infty)=2$.
After renumeration, we have  $(F\cdot C_2)=1$ or    $(F\cdot C_0)=(F\cdot C_7)=1$ or $(F\cdot C_1)=1$ or $(F\cdot C_0)=2$, with otherwise zero intersection numbers.
Our task is to rule out these four possibilities.

Suppose first that $(F\cdot C_2)=1$. Then $C_0+C_1+C_3+\ldots+C_7$ is a chain that is vertical with respect to $\psi$,
hence belongs to $\psi^{-1}(\infty)$. 
Its Kodaira symbol must be $\III^*$, by Propositions \ref{key reduction} and \ref{no I4*}.
We have $\psi^\inv_\red(\infty)= C_0+C_1+C_3+\ldots+C_7+R$,
with $(R\cdot C_4)=1$. This gives 
$$
4\geq \varphi^{-1}(\infty)\cdot \psi^\inv_\ind(\infty) = (4C_4+2C_2)\cdot 2R =8+ 4(C_2\cdot R)\geq 8,
$$
contradiction.

Next suppose that $(F\cdot C_0)=(F\cdot C_7)=1$. Then $C_1+\ldots+C_6$ has dual graph $T_{2,3,3}$ and is vertical with
respect to $\psi$, whence is contained in $\psi^\inv(\infty)$.
As above, the Kodaira symbol is $\III^*$, 
so $\psi^\inv_\red(\infty) = C_0'+C_1+\ldots+C_6+C_7'$, with $(C_0'\cdot C_1)=(C_7'\cdot C_6)=1$. Then
$$
4\geq \varphi^{-1}(\infty)\cdot \psi^\inv_\ind(\infty) = (C_0+2C_1+2C_6+C_7)\cdot(C_0'+C_7') \geq  4+  \sum (C_i\cdot C_j'),
$$
where the sum runs over $i,j\in \{0,7\}$. It follows that $C_0'\cup C_7'$ is disjoint from $C_0\cup C_7$.
Thus the curves $C_0'+C_0+\ldots C_7+C_7'$ forms a configuration as in Proposition \ref{no I4*+R}, contradiction.

We come to  the case $(F\cdot C_1)=1$, which is the most challenging.
The curve $C_2+\ldots+C_7$ has dual graph $T_{2,2,4}$ and is $\psi$-vertical, thus contained in $\psi^\inv(\infty)$.
Now the possible Kodaira symbols are  $\I_4^*$ or $\III^*$  or $\I_2^*$. 
The former was already discarded in  Proposition \ref{no I4*}. 
Suppose the symbol is $\III^*$.  It is simple by Proposition \ref{no multiple III*}, 
and with Proposition \ref{possible configurations} we infer that
also $C_0$ belongs to $\psi^\inv(\infty)$. Consequently  $\psi_\red^{-1}(\infty)=C_0+R+C_2+\ldots+C_7$ for some additional component $R$
with either $(R\cdot C_0)=(R\cdot C_3)=1$ or $(R\cdot C_0)=(R\cdot C_2)=1$.  In the former case
$$
4\geq \varphi^{-1}(\infty)\cdot \psi^\inv_\ind(\infty) = (C_0+2C_1+3C_3)\cdot 2R = 2 + 4(C_1\cdot R) + 6\geq 8,
$$
contradiction, and the latter case is treated analogously. So our fiber has Kodaira symbol $\I_2^*$,
and $\psi^\inv_\red(\infty)=C_2+\ldots+C_7+R$ for some additional component $R$ with $(R\cdot C_6)=1$. Now
$$
4\geq \varphi^{-1}(\infty)\cdot \psi^\inv_\ind(\infty) = (2C_1+ 2C_6)\cdot R = 2(C_1\cdot R) + 2.
$$
This gives $(C_1\cdot R)\leq 1$. But the case $(C_1\cdot R)=1$ yields a curve $R+C_1+C_3+\ldots+C_6$ of canonical type $\I_6$,
contradiction. Thus $C_1\cap R=\varnothing$.

According to Proposition \ref{key reduction}, $\psi$ has $\I_2^*+\III+\III$
as  configuration of fibers over rational points. The curve $C_0$ is $\psi$-vertical, and
we can assume $\psi^\inv_\ind(0)=C_0+C_0'$. Write $\psi^\inv_\ind(1)=D+D'$.
The dual graph for all our curves takes the following form:
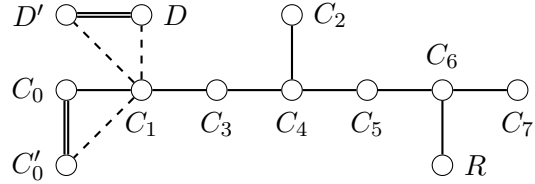
\begin{equation*}
\begin{gathered}
\begin{tikzpicture}
[node distance=1cm, font=\small]
\tikzstyle{vertex}=[circle, draw, fill=white, inner sep=0mm, minimum size=1.5ex]
\node[vertex]	(C0)  	at (1,0)  [		label=left:{$C_0$}] {};
\node[vertex]	(C1)	[right of=C0,	label=below:{$C_1$}] {};
\node[vertex]	(C3)	[right of=C1,	label=below:{$C_3$}] {};
\node[vertex]	(C4)	[right of=C3,	label=below:{$C_4$}] {};
\node[vertex]	(C5)	[right of=C4,	label=below:{$C_5$}] {};
\node[vertex]	(C6)	[right of=C5,	label=above:{$C_6$}] {};
\node[vertex]	(C7)	[right of=C6,	label=below:{$C_7$}] {};
\node[vertex]	(C2)	[above of=C4,	label=right:{$C_2$}] {};

\node[vertex]	(R)		[below of=C6, 	label=right:{$R$}] {};
\node[vertex]	(C0')	[below of=C0, 	label=left:{$C_0'$}] {};
\node[vertex]	(D)		[above of=C1, 	label=right:{$D$}] {};
\node[vertex]	(D')	[above of=C0, 	label=left:{$D'$}] {};

\draw[thick]    (C0)--(C1)--(C3)--(C4)-- (C5)--(C6)--(C7);
\draw[thick]    (C2)--(C4);
\draw[thick]    (R)--(C6);
\draw[thick,double]  (D)--(D');
\draw[thick,double]  (C0)--(C0');
\draw[thick,dashed]  (C1)--(C0');
\draw[thick,dashed]  (C1)--(D);
\draw[thick,dashed]  (C1)--(D');

\end{tikzpicture}
\end{gathered}
\end{equation*}
According to Proposition \ref{possible configurations}, the simple fiber $\varphi^\inv(\infty)$ with Kodaira symbol $\III^*$
must intersect each $(-2)$-curve on $Y$. It follows that  $C_1$ intersects $D$ and $D'$.
From
$4\geq \varphi^{-1}(\infty)\cdot \psi^\inv_\ind(1) = 2C_1\cdot (D+D')$
we infer $(C_1\cdot D)=(C_1\cdot D')=1$. Thus $D+C_0+\ldots+C_7+R$ forms a configuration as in Proposition \ref{no I4*+R}, 
again a contradiction.

It remains to deal with the case $(F\cdot C_0)=2$. Then $C_1+\ldots+C_7$ has dual graph $T_{2,3,4}$ and is vertical with respect to $\psi$,
whence contained in $\psi^\inv(\infty)$. Its Kodaira symbol must be $\III^*$, by Proposition \ref{key reduction}.
Write $\psi^\inv_\red(\infty)=C_0'+C_1+\ldots+C_7$ for some additional component $C_0'$ with $(C_0'\cdot C_1)=1$.
The fiber is simple, by Proposition \ref{no multiple III*}, so
$$
4=\varphi^\inv(\infty)\cdot \psi^\inv(\infty) = \varphi^\inv(\infty)\cdot C_0' = (C_0+ 2C_1)\cdot C_0' = (C_0\cdot C_0') +2.
$$
Whence $(C_0\cdot C_0')=2$, and $C_0+C_0'$ is a curve of canonical type, with Kodaira symbol $\I_2,\tI_2$ or $\III$. 
Let $\eta:Y\ra\PP^1$ be the
ensuing genus-one fibration. By Proposition \ref{key reduction}, we may assume that $\eta^{-1}(\infty)$ has 
symbol $\III^*$ or $\I_2^*$. In the former case Proposition \ref{no multiple III*} ensures that it is simple
so by  Proposition \ref{possible configurations} the configuration must be $\III^*+E_4+E_4$, contradiction.
In light of Proposition \ref{possible configurations}, the only remaining possible configuration is $\I_2^*+\III+\III$.
Without restriction $\eta_\ind^{-1}(0)=C_0+C_0'$ and $\eta_\ind^{-1}(1)=D+D'$. 
The curve  $C_2+\ldots+C_7$ has dual graph $T_{2,2,4}$ and thus belongs to $\eta^{-1}(\infty)$.
Let $R$ be the additional component.  The dual graph  of  our curves takes the following form:
\begin{equation*}
\begin{gathered}
\begin{tikzpicture}
[node distance=1cm, font=\small]
\tikzstyle{vertex}=[circle, draw, fill=white, inner sep=0mm, minimum size=1.5ex]
\node[vertex]	(C0)  	at (1,0)  [		label=left:{$C_0$}] {};
\node[vertex]	(C1)	[right of=C0,	label=below:{$C_1$}] {};
\node[vertex]	(C3)	[right of=C1,	label=above:{$C_3$}] {};
\node[vertex]	(C4)	[right of=C3,	label=below:{}] {};
\node[vertex]	(C5)	[right of=C4,	label=below:{}] {};
\node[vertex]	(C6)	[right of=C5,	label=above:{$C_6$}] {};
\node[vertex]	(C7)	[right of=C6,	label=below:{$C_7$}] {};
\node[vertex]	(C2)	[above of=C4,	label=right:{$C_2$}] {};

\node[vertex]	(R)		[below of=C6, 	label=right:{$R$}] {};
\node[vertex]	(C0')	[below of=C0, 	label=left:{$C_0'$}] {};
\node[vertex]	(D)		[above of=C1, 	label=right:{$D$}] {};
\node[vertex]	(D')	[above of=C0, 	label=left:{$D'$}] {};

\draw[thick]    (C0)--(C1)--(C3)--(C4)-- (C5)--(C6)--(C7);
\draw[thick]    (C2)--(C4);
\draw[thick]    (R)--(C6);
\draw[thick,double]  (D)--(D');
\draw[thick,double]  (C0)--(C0');
\draw[thick]  (C1)--(C0');
\draw[thick,dashed]  (C1)--(D);
\draw[thick,dashed]  (C1)--(D');
\draw[thick,dashed]  (C1)--(R);
\end{tikzpicture}
\end{gathered}
\end{equation*}
We have 
$\eta^{-1}_\ind(0)\cdot \varphi^{-1}(\infty) = (C_0+C_0')\cdot 2C_1 =   2+2=4$.
If $\eta^{-1}(0)$ is simple, each multiple fiber $2F'$ for $\eta:Y\ra\PP^1$ has $(F'\cdot C_1) =1$, in contradiction to the previous 
paragraph. Thus $\eta^{-1}(0)$ is multiple, such that  $\eta^{-1}(\infty)\cdot\varphi^{-1}(\infty)= 8$.
If the fiber $\eta^{-1}(\infty)$ is multiple as well, we  get
$$
4=\eta_\ind^{-1}(\infty) \cdot \varphi^{-1}(\infty) = R\cdot (2C_1+2C_6) = 2(R\cdot C_1) + 2,
$$
hence $C_3+\ldots+C_6+R+C_1$ is a curve of canonical type $\I_6$, contradiction.
Thus $\eta^{-1}(\infty)$ is simple, whereas $\eta^{-1}(1)$ must be the second multiple fiber. The latter gives
$$
4=\eta^{-1}_\ind(1)\cdot\varphi^{-1}(\infty) = (D+D')\cdot 2C_1 = 2(D\cdot C_1)+2(D'\cdot C_1).
$$
According to Proposition \ref{possible configurations}, the fiber $\varphi^{-1}(\infty)$ intersects both $D$ and $D'$,
so we have  $(D\cdot C_1)=(D'\cdot C_1)=1$.  
Thus $D+C_1+\ldots+C_7$ gives a curve of canonical type $\III^*$. This gives a further fibration on $Y$, in which 
  the   fiber of type $\III^*$  must be simple, according to Proposition \ref{no multiple III*}. This gives
$$
4\leq (D+2C_1+\ldots+2C_6+C_7)\cdot\varphi^{-1}(\infty)=D\cdot 2C_1=2,
$$
contradiction.
\qed

\section{Surfaces  with restricted configurations}
\mylabel{Restricted configurations}

Let $k$ be an algebraically closed ground field $k$ of arbitrary characteristic $p\geq 0$,
and $Y$ be a classical Enriques surface. The goal of this section is
to show that certain configurations of $(-2)$-curves cannot exist, once one imposes rather severe 
restrictions on the configuration of reducible fibers in genus-one fibrations.

\begin{proposition}
\mylabel{only two configurations}
There is no classical Enriques surface $Y$ such that every genus-one fibration $\varphi:Y\ra\PP^1$
satisfies one of the following properties:
\begin{enumerate}
\item
The configuration of reducible fibers is $\I_1^*+\I_4$, where the $\I_1^*$-fiber is multiple.
\item
The configuration  is $\I_2^*+\III+\III$, where two of these fibers are multiple.
\end{enumerate}
\end{proposition}

\proof 
Seeking a contradiction, we assume that such a surface exists. Choose a genus-one fibration $\varphi:Y\ra \PP^1$.
By \cite{Ekedahl; Shepherd-Barron 2004}, Theorem A the surface $Y$ is non-exceptional.
According to Proposition \ref{properties non-exceptional} there is another genus-one fibration $\psi:Y\ra\PP^1$,
with $\varphi^\inv(\infty)\cdot\psi^\inv(\infty)=4$. Choose a multiple fiber  $2F\subset Y$ for   $\psi$.
Without restriction, we may assume that the fibers of $\varphi$ and $\psi$ with Kodaira symbol $\I_n^*$ occur over
$t=\infty$. 

Suppose first that $\varphi^\inv(\infty)$ has Kodaira symbol $\I_1^*$, and hence must be multiple. The dual graph is:
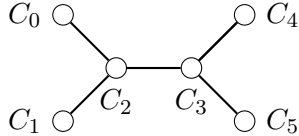
\begin{equation*}
\begin{gathered}
\begin{tikzpicture}
[node distance=1cm, font=\small]
\tikzstyle{vertex}=[circle, draw, fill=white, inner sep=0mm, minimum size=1.5ex]
\node[vertex]	(C2)  	at (1,0) [			label=below:{$C_2$}]	{};
\node[vertex]	(C3)	[right of=C2, 		label=below:{$C_3$}]	{};
\node[vertex]	(C4)	[above right of=C3, label=right:{$C_4$}]	{};
\node[vertex]	(C5)	[below right of=C3, label=right:{$C_5$}]	{};
\node[vertex]	(C0)	[above left of=C2, 	label=left :{$C_0$}]	{};
\node[vertex]	(C1)	[below left of=C2, 	label=left :{$C_1$}]	{};
 
\draw [thick] (C2)--(C3)--(C4);
\draw [thick] (C0)--(C2)--(C1);
\draw [thick] (C5)--(C3)--(C4);
\end{tikzpicture}
\end{gathered}
\end{equation*}
After renumeration, we may assume $(F\cdot C_0)=1$. The curve $C_1+\ldots +C_5$ has dual graph 
$T_{2,2,3}$ and is $\psi$-vertical, hence must belong to $\psi^\inv(\infty)$. Suppose the latter has Kodaira symbol
$\I_1^*$. Then $\psi^{-1}_\red(\infty)=R+C_1+\ldots+C_5$ for some additional component 
$R$ with $(R\cdot C_2)=1$, and the fiber is multiple. This gives the contradiction
$$
1=\varphi^\inv_\ind(\infty)\cdot\psi^\inv_\ind(\infty) = (C_0+2C_2)\cdot R = (C_0\cdot R) + 2\geq 2.
$$
Thus $\psi^\inv(\infty)$ has Kodaira symbol $\I_2^*$, and we write $\psi^\inv_\ind(\infty)=R+R'+C_1+\ldots+C_5$
with additional components $R$ and $R'$, which have  $(R\cdot C_1)=(R'\cdot C_1)=1$. Then
$$
2\geq \varphi^\inv_\ind(\infty)\cdot\psi^\inv_\ind(\infty) = (C_0+C_1)\cdot (R+R') = C_0\cdot(R+R')  + 1+1.
$$
Hence $R\cup R'$ is disjoint from $C_0$, and $\psi^\inv(\infty)$ must be simple. 
Without restriction, we may assume that $\psi^\inv(0)$ is multiple with Kodaira symbol $\III$.
Write $\psi^\inv_\ind(0)=D+D'$. Then
$$
2(D+D')\cdot C_0=\psi^\inv(\infty)\cdot C_0 =2C_2\cdot C_0 = 2.
$$
Without restriction, we may assume $(C_0\cdot D)=1$ and $(C_0\cdot D')=0$.
Then the dual graph for our curves becomes:
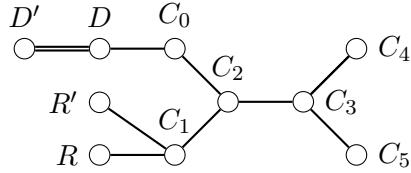
\begin{equation*}
\begin{gathered}
\begin{tikzpicture}
[node distance=1cm, font=\small]
\tikzstyle{vertex}=[circle, draw, fill=white, inner sep=0mm, minimum size=1.5ex]
\node[vertex]	(C2)  	at (1,0)	[		label=above:{$C_2$}]	{};
\node[vertex]	(C3)	[right of=C2, 		label=right:{$C_3$}]	{};
\node[vertex]	(C4)	[above right of=C3, label=right:{$C_4$}]	{};
\node[vertex]	(C5)	[below right of=C3, label=right:{$C_5$}]	{};
\node[vertex]	(C0)	[above left of=C2, 	label=above:{$C_0$}]	{};
\node[vertex]	(C1)	[below left of=C2, 	label=above   :{$C_1$}]	{};
 
\node[vertex]	(R)	[left of=C1, 		label=left:{$R$}]	{};
\node[vertex]	(R')	[above  of=R, yshift=-.3cm,	label=left:{$R'$}]	{};
\node[vertex]	(D)	[left of=C0, 		label=above:{$D$}]	{};
\node[vertex]	(D')	[left of=D, 		label=above:{$D'$}]	{};

\draw [thick] (C2)--(C3)--(C4);
\draw [thick] (C0)--(C2)--(C1);
\draw [thick] (C5)--(C3)--(C4);
\draw [thick] (C0)--(D);
\draw [thick, double] (D')--(D);

\draw [thick] (R)--(C1)--(R');

\end{tikzpicture}
\end{gathered}
\end{equation*}
Thus $R+D+C_0+\ldots+C_4$ is a curve of canonical type $\IV^*$, contradiction.

\emph{ Summing up, every  genus-one fibration on $Y$ 
 has   $\I_2^*+\III+\III$ as configuration of degenerate fibers, two of which are multiple.}
We tacitly assume that  the $\I_2^*$-fiber lies over $t=\infty$, whereas the
$\III$-fibers appear over $t=0,1$. 
In particular,  the  dual graph for the reducible fibers of our $\varphi:Y\ra\PP^1$ is:
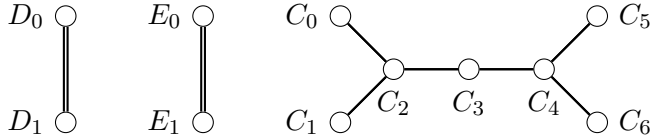
\begin{equation*}
\begin{gathered}
\begin{tikzpicture}
[node distance=1cm, font=\small]
\tikzstyle{vertex}=[circle, draw, fill=white, inner sep=0mm, minimum size=1.5ex]
\node[vertex]	(C2)  	at (1,0) [			label=below:{$C_2$}]	{};
\node[vertex]	(C3)	[right of=C2, 		label=below:{$C_3$}]	{};
\node[vertex]	(C4)	[right of=C3, 		label=below:{$C_4$}]	{};
\node[vertex]	(C5)	[above right of=C4, label=right:{$C_5$}]	{};
\node[vertex]	(C6)	[below right of=C4, label=right:{$C_6$}]	{};
\node[vertex]	(C0)	[above left of=C2, 	label=left :{$C_0$}]	{};
\node[vertex]	(C1)	[below left of=C2, 	label=left :{$C_1$}]	{};

\node[vertex]	(E0)			[left  of=C0, xshift=-2em, label=left :{$E_0$}]	{};
\node[vertex]	(E1)			[left  of=C1, xshift=-2em, label=left :{$E_1$}]	{};
\node[vertex]	(D0)			[left  of=E0, xshift=-2em, label=left :{$D_0$}]	{};
\node[vertex]	(D1)			[left  of=E1, xshift=-2em, label=left :{$D_1$}]	{};

\draw [thick] (C2)--(C3)--(C4);
\draw [thick] (C0)--(C2)--(C1);
\draw [thick] (C6)--(C4)--(C5);
\draw [thick,double] (D0)--(D1);
\draw [thick,double] (E0)--(E1);

\draw [thick] (C2)--(C3)--(C4);
\draw [thick] (C0)--(C2)--(C1);
\draw [thick] (C6)--(C4)--(C5);
\end{tikzpicture}
\end{gathered}
\end{equation*}
We first verify that $\varphi^{-1}(\infty)$ is simple. Seeking a contradiction, we assume that it is multiple,
such that $F\cdot \varphi^{-1}_\ind(\infty)=1$. After renumeration  we may assume $(F\cdot C_0)=1$.
The connected curve $C_1+\ldots+C_6$ has dual graph $T_{2,2,4}$, hence is contained in $\psi^{-1}(\infty)$.
Let $R\subset\psi^{-1}(\infty)$ be the additional component, which has $(R\cdot C_2)=1$.
From
$$
2\geq \varphi^{-1}_\ind(\infty)\cdot \psi^{-1}_\ind(\infty) = (C_0+2C_2)\cdot R = (C_0\cdot R) + 2
$$
we infer $R\cap C_0=\varnothing$. Thus $C_0+\ldots +C_3+R$ supports a curve of canonical type $\I_0^*$, 
contradiction.

\emph{Summing up, in  all genus-one fibrations the $\I_2^*$-fiber   is simple and the $\III$-fibers  
are multiple.} Recall that $2F$ is a multiple fiber for our fibration $\psi:Y\ra\PP^1$ with $\varphi^\inv(\infty)\cdot\psi^\inv(\infty)=4$.
After renumeration, we may assume that $(F\cdot D_0)=(F\cdot E_0)=1$ and $(F\cdot D_1)=(F\cdot E_1)=0$.
In particular, the curves $D_0$, $E_0$ are $\psi$-horizontal, whereas $D_1$, $E_1$ are $\psi$-vertical.

By symmetry, we now write  $\psi^{-1}_\ind(0)=R_0+R_1$ such that $R_0$ is $\varphi$-horizontal and  $R_1$ is $\varphi$-vertical.
We claim that $R_1$ is not contained in $\varphi^{-1}(0)$. 
For $R_1=D_1$ we get the contradiction
$$
1=(D_0+D_1)\cdot (R_0+R_1) = D_0\cdot (R_0+R_1)\geq (D_0\cdot R_1)=(D_0\cdot D_1)=2,
$$
whereas $R_1=D_0$ leads to  
$$
1=(D_0+D_1)\cdot (R_0+R_1) = D_1\cdot(R_0+R_1)    = 0,
$$
again a contradiction. Thus $R_1$ is disjoint from $\varphi^{-1}(0)$, which gives 
$$
(D_0\cdot R_0) = D_0\cdot (R_0+R_1) = (D_0+D_1)\cdot (R_0+R_1) = 1.
$$
By symmetry, we also get $(E_0\cdot R_0)=1$.
Likewise, we write $\psi^\inv_\ind(1)=S_0+S_1$ such that $S_0$ is $\varphi$-horizontal and  $S_1$ is $\varphi$-vertical,
and see $(D_0\cdot S_0)=(E_0\cdot S_0)=1$.
Thus   $D_0+ R_0+E_0+S_0$ is a curve of canonical type $\I_4$, contradiction.
\qed

\section{Proof of the main result}
\mylabel{Proof}

We are finally ready to prove Theorem \ref{no enriques over integers},  the main result of this paper.
The task is to show that the fiber category $\shM_\Enr(\ZZ)$ for the stack of Enriques surfaces is empty.
Seeking a contradiction, we assume that that there is a family of Enriques surfaces
$f:\foY\ra \Spec(\ZZ)$. According to Proposition \ref{constant picard scheme}, the Picard scheme $\Pic_{\foY/\ZZ}$ is constant.
So the closed fiber $Y=\foY\otimes\FF_2$ also has constant Picard scheme. 
It is  non-exceptional, by  \cite{Schroeer 2021b}, Theorem 7.2, because this Enriques surface
 comes with a deformation over $W_2(\FF_2)=\ZZ/4\ZZ$.
Therefore it suffices to establish:

\begin{theorem}
\mylabel{no non-exceptional}
There is no  Enriques surface $Y$ over the prime field $k=\FF_2$  that is  non-exceptional and  whose   Picard scheme
is constant.
\end{theorem}

\proof
Suppose such a $Y$ exists.
By Theorem \ref{family with constant pic} there is a genus-one fibration $\varphi:Y\ra\PP^1$.
For each geometric point $\ba:\Spec(\Omega)\ra \PP^1$ lying over a non-rational closed point  $a\in\PP^1$, the geometric fiber $Y_\ba$
is  irreducible, according to Proposition \ref{constant num over finite field},  thus the  Kodaira symbol is $\I_0, \I_1$ or $\II$.
Moreover, it must be simple by Theorem \ref{family with constant pic}.
The possible configurations of fibers over the rational points $t=0,1,\infty$ are tabulated
in Proposition \ref{key reduction}.  
 In Sections \ref{Elimination I4*} and \ref{Elimination III*} we saw that the Kodaira symbols
$\I_4^*$ and $\III^*$ actually do not occur.
For the classical Enriques surface $\bY=Y\otimes_kk^\alg$   this means that for every genus-one fibration, the configuration
of degenerate fibers is either $\I_1^*+\I_4$ or $\I_2^*+\III+\III$. Furthermore, in the former case the $\I_1^*$-fiber is multiple,
whereas in the latter case two of the three fibers are multiple.
 We showed in Proposition \ref{only two configurations} that
such an  Enriques surface does not exist. This gives the desired contradiction.
\qed


\end{document}